\documentclass[11pt,reqno]{amsart}

\usepackage{mathrsfs,amsmath}

\usepackage{hyperref} 
\usepackage{amsrefs} 

\DeclareMathOperator{\vMax}{\underline{Max}}

\DeclareMathOperator{\vword}{w-ord}

\newcommand{\calo}{{\mathcal {O}}}

\DeclareMathOperator{\Sing}{Sing}
\DeclareMathOperator{\Spec}{Spec}
\DeclareMathOperator{\ord}{ord}
\DeclareMathOperator{\Bl}{Bl}
\DeclareMathOperator{\DDiff}{\mathbb{D}iff}
\DeclareMathOperator{\diff}{Diff}
\renewcommand{\ord}{\operatorname{ord}}

\usepackage{latexsym, amsfonts, amsmath, amssymb, amscd, epsfig}

\usepackage[usenames]{color}
\usepackage{cancel}
\usepackage[normalem]{ulem}
\usepackage{dsfont}
\input xy             
\xyoption{all}    
\usepackage{amscd}

\newtheorem{Theorem}{Theorem}[section]
\newtheorem{Lemma}[Theorem]{Lemma}
\newtheorem{Proposition}[Theorem]{Proposition}

{\theoremstyle{definition}

\newtheorem{Definition}[Theorem]{Definition}
\newtheorem{Paragraph}[Theorem]{}
\newtheorem{Example}[Theorem]{Example}
\newtheorem{Remark}[Theorem]{Remark}
\newtheorem{Question}[Theorem]{Question}
}

\title[On Rees algebras and invariants for singularities\ldots]{On Rees algebras and invariants for singularities over 
perfect fields}

\author{A. Bravo, M.L. Garcia-Escamilla, O.E. Villamayor U.}

\address{Depto. Matem\'aticas,
Facultad de Ciencias, Universidad Aut\'onoma
de Madrid and Instituto de Ciencias Matematicas CSIC-UAM-UC3M-UCM, Canto Blanco 28049 Madrid, Spain}
\email{ana.bravo@uam.es, mariluz.garcia@uam.es, villamayor@uam.es}
\thanks{The authors were partially supported by MTM2009-07291.}


 \date{}

\begin{document}
\begin{abstract}
The purpose of this paper is to show how  Rees algebras can be applied in 
the study  of singularities embedded in smooth schemes over perfect fields.  In particular, we will study situations in which the multiplicity of a hypersurface is a multiple of the characteristic. 
As another application,  here we  indicate how the use  of these algebras  has trivialized  local-global questions 
in resolution of singularities over fields of characteristic zero.

\end{abstract}
\maketitle

\tableofcontents

\part{Introduction and statement of main results}
\section*{Introduction}
The purpose of this paper is twofold: on the one hand, we study some topics in commutative algebra concerning Rees algebras,  differential operators, and the theory of ideals on smooth schemes; on the other, 
we show how these results can we applied to different  problems in resolution of singularities.

\

\noindent {\bf Rees algebras, differential operators and integral closure.}

\

When working on smooth schemes over a perfect field, the action of differential operators on ideals 
leads to the definition of interesting invariants (see for instance the role of differential operators in connection with  test ideals in \cite{AMBL}, see also \cite[Remark 2.6]{BMS} or \cite{BMS2}). In addition it is very  frequent to come up with information that does not distinguish between an ideal and its integral closure,  an issue  that is quite natural to expect from the geometric point of view. 

\

Here we will be  concerned on certain properties of ideals that are preserved by integral closure, and 
that are not affected by considering the action of 
differential operators (in some sense  that will be  made precise along this manuscript). It is within this context that we find convenient to use Rees algebras as a suitable tool to 
combine both aspects.  Our results have applications in problems that arise in resolution of singularities.

\

\noindent {\bf Problems in Resolution of singularities.}

\

\begin{Paragraph}\label{Problem_1} {\bf Problem 1. On the multiplicity of hypersurfaces.}
\end{Paragraph}

Let $V$ be a smooth scheme over a perfect field $k$. 
It is quite natural to address the  resolution of singularities of a hypersurface $X\subset V$ 
  by paying attention to its maximum multiplicity locus, say $M_X\subset X$. In fact, a key point in  resolution  over fields of characteristic zero, 
    is the existence of the so called {\em hypersurfaces of maximal contact}. This means that, locally, in 
  a neighborhood of each $x\in M_X$, there is a smooth hypersurface containing $M_X$. Moreover, this containment 
  is preserved by suitably chosen monoidal transformations until the maximum multiplicity drops. This allows us  to 
  solve singularities by induction on the number of variables, since the original problem of lowering 
  the maximum multiplicity of $X$ is equivalent to another that can be formulated on a smooth scheme of lower dimension. 
  
\

As an example, let $V=\Spec k[z,x]$, and  
consider the affine curve  $C:=V(\langle z^2+x^3\rangle)$.  Its singular locus is  $(0,0)$, which is also  a twofold  point.  It can be checked, that, if $\operatorname{char}(k)\neq 2$, then 
$\{x=0\}$ is a hypersurface of maximal contact. Here the problem of lowering the maximum multiplicity of $C$ is equivalent to 
that of lowering the maximum order of the ideal $\langle x^3\rangle\subset k[x]$ below 2. We started with a problem 
in two variables  that,  in fact, is equivalent to a problem in the affine line. 

\

If $\operatorname{char}(k)=2$ we will show that  the problem of lowering the 
maximum multiplicity of $C$ cannot be reformulated as an equivalent problem in one dimension less. Thus, the curve 
$C$ is already   a counter-example to the possibility of using inductive arguments  in positive 
characteristic, at least if we intend to consider the multiplicity as the only invariant. 

\

\begin{Paragraph}\label{Problem_2} {\bf Problem 2. On constructive log-resolution of ideals.}
\end{Paragraph}

Let $V$ be a smooth scheme over a perfect field $k$, and let $J\subset {\mathcal O}_V$ 
be a non-zero sheaf of ideals. A {\em log-resolution of  $J$}   is a proper and birational morphism  of smooth schemes, $V\stackrel{\rho}{\longleftarrow} \mathfrak{V}$,   so that 
$J{\mathcal O}_{\mathfrak{V}}$ is locally invertible,  and supported on smooth hypersurfaces with normal crossings. 

\

A  {\em constructive 
log-resolution of  a sheaf of ideals}  $J$  consists on  
the construction of  a  finite sequence of blow ups at smooth centers, 
\begin{equation}
\label{LogRes}
\xymatrix@R=2pt@C=20pt{
V=V_0 & V_1 \ar[l]_{\hspace{25pt}\rho_{0}} &    \ldots  \ar[l]_{\rho_{1}} & V_{n} \ar[l]_{\rho_{n-1}}\\
J=J_0 &  J{\mathcal O}_{V_1} &  \ldots    & J{\mathcal O}_{V_n},
}
\end{equation}
so that the total  transform of $J$ in $V_n$, say $J{\mathcal O}_{V_n}$,  is of the form 
$$J{\mathcal O}_{V_n}=I(H_1)^{b_1}\cdots I(H_n)^{b_n},$$
where $H_1,\ldots,H_n\subset V_n$ are smooth hypersurfaces with normal crossings, and $b_i\in {\mathbb N}$ for $i=1, \ldots, n$. As an application, see \cite{EncVillaIberoamericana} for a procedure to 
construct a resolution of singularities of variety $X$ embedded in a smooth scheme $V$ over a field of characteristic zero,  via the constructive log-resolution of the ideal ${\mathcal I}(X)\subset {\mathcal O}_V$. 

\
 
Over fields of characteristic zero, constructive log-resolution is usually achieved by defining upper-semi continuous 
functions mapping to some well ordered set $(\Lambda, \geq)$, say  
\begin{equation}
\label{uppersemi}
f_i: V(J{\mathcal O}_{V_i})\to (\Lambda, \geq). 
\end{equation} 
These functions are so that if $\max f_i$ denotes the maximum value of $f_i$, then 
$$\vMax\ f_i:=\{ x\in V(J{\mathcal O}_{V_i}): f_i(x)=\max f_i\}$$
is smooth, and  determines the center of the monoidal transformation  
$$\xymatrix@R=2pt@C=20pt{V_i & \ar[l]_{\rho_{i}} V_{i+1}.}$$
Moreover if a sequence like (\ref{LogRes}) is defined by blowing up at $\vMax\ f_i$, for $i=0,1,\ldots, n-1$, then  
$$\max f_0>\max f_1>\ldots > \max f_{n-1}>\max f_n. $$
Here  $f_n$ being  constant on $V(J{\mathcal O}_{V_n})$ would mean  that $J{\mathcal O}_{V_n}$ is locally invertible and supported on smooth divisors with normal crossings. So, at each step of a 
resolution process,  the functions 
$f_i$  measure how far we are from  achieving  resolution.  

\

The   $f_i$ are 
usually referred to as {\em resolution functions}. They are traditionally constructed locally,  
in a neighborhood of each point, by successive    use of the so
 called {\em Hironaka's order function} (see \ref{Def_Pares} below), and {\em restriction to hypersurfaces 
of maximal contact}. These 
are all arguments of local nature, which indicate, locally, which is the center to be blown up. 
 It is highly non-trivial to prove that all these 
locally defined invariants patch so as to produce a global function,   hereby providing a 
(global)  
smooth center to blow up. This is usually  usually refer to as the {\em local-global problem in resolution of singularities}.  

\

We will show, that the locally defined invariants from above, in fact globalize because there is 
a {\em canonical way to construct them}. 

\

\begin{Paragraph}\label{Def_Pares}
{\bf The language of pairs.}  
\end{Paragraph}
 
In this paragraph we introduce the language of {\em pairs}. As we will see, Problems 1 and 2 
have a natural formulation in this context. 

\

Let $V$ be a smooth scheme over a perfect field $k$. A  {\em pair} is given by a couple $(J,b)$ where $J$ is a non-zero sheaf of 
ideals, and $b$ is a non-negative integer. The   {\em singular locus 
of a pair}  is the closed set 
\begin{subequations}
\begin{equation}
\label{SingHironaka}
\Sing(J,b):=\{x\in V: \nu_x(J)\geq b\},
\end{equation}
where $\nu_x(J)$ denotes the order of $J$ in the regular local 
ring ${\mathcal O}_{V,x}$. 

\

With this notation, {\em Hironaka's order 
function} is defined as 
\begin{equation}
\label{ordenHironaka}
\begin{array}{rrcl}
\ord_{(J,b)}:& \Sing(J,b) &  \longrightarrow & {\mathbb Q}_{\geq 0}\\
 & x & \longmapsto & \ord_{(J,b)}(x):=\frac{\nu_x(J)}{b}.\end{array}
\end{equation}

\

A {\em permissible center} is a  smooth closed subscheme  $Y\subset \Sing(J,b)$. The  {\em transform of $(J,b)$}   
after blowing up at a permissible center $Y$, 
$$\xymatrix@R=2pt@C=20pt{
V& \ar[l]_{\rho} V_1
}
$$
is defined as
a pair  
$(J_1,b)$ where 
\begin{equation}\label{trasformadoHironaka}J_1:=I(H_1)^{-b}J{\mathcal O}_{V_1} \end{equation}\end{subequations}
and $I(H_1)$ is the defining ideal of the exceptional divisor $\rho^{-1}(Y)$. 
Observe that  pairs have an easy-to-handle law of transformation under blow ups. A {\em resolution of 
a pair} is a finite sequence of blow ups at permissible centers 
$$
\xymatrix@R=2pt@C=20pt{
V=V_0 & V_1 \ar[l]_{\hspace{20pt} \rho_{0}} &    \ldots  \ar[l]_{\rho_{1}} & V_{n} \ar[l]_{\rho_{n-1}}\\
(J,b)=(J_0,b) & (J_1,b) & \ldots & (J_n,b),
}
$$
such that: 
\begin{itemize}
\item $\Sing(J_n,b)=\emptyset$, and 
\item the  exceptional locus of the composition $V\leftarrow V_n$ is a union of smooth hypersurfaces having only normal crossings in $V_n$.
\end{itemize}

\

\noindent {\bf Regarding to Problem \ref{Problem_1}.} With the same notation as in \ref{Problem_1}, observe that we can associate to  the maximum multiplicity locus of $C$ the pair $(\langle z^2+x^3\rangle,2)$. Lowering the maximum multiplicity of $C$ amounts to 
resolving the pair $(\langle z^2+x^3\rangle,2)$.  The question is:  can  be  found another pair, in one variable less, whose resolution induces a resolution of $(\langle z^2+x^3\rangle,2)$?

\

\noindent {\bf Regarding to Problem \ref{Problem_2}.} In characteristic zero, a constructive log-resolution of a sheaf of ideals $J$  can be achieved via constructive 
resolution 
of pairs. It is quite straightforward to observe that a resolution of $(J,1)$ provides a log-resolution of $J$. But a priori it may not be clear why to care about resolutions of pairs of the form $(J,b)$ with $b>1$, and this issue deserves some explanation.  Although a given log-resolution problem can be stated as the resolution of a pair of the form $(J,1)$, the procedure followed in constructive resolution requires the use of some form of induction and stratifying functions that immediately forces the appearance of pairs of more general form. These are arguments of local nature, and it is not obvious, a priori, that they  lead to the definition of global invariants and resolution functions. 

\begin{Paragraph}\label{Integral_Differential} {\bf Local sequences and equivalence of pairs.}
\end{Paragraph}

As indicated in \ref{Def_Pares}, we can reformulate Problems 1 and 2 in terms of pairs.  However,  the assignation of   a pair to a given problem is not unique. To fix ideas, assume that $B$ is a smooth algebra of finite type over a perfect field $k$, and let $J\subset B$ be an ideal. Let us consider two different situations:  

\

On the one hand, observe that   a resolution of $(J,1)$ is also a resolution of $(J^2,2)$. More generally, if two pairs, $(J,b)$ and $(I,c)$ are such that the ideals  $J^c$ and $I^b$ have the same integral closure, then one expects that  they will undergo the same resolution. It is not hard to see that both pairs 
have the same singular locus: note that $\nu_x(J^c)=\nu_x(I^b)$ for all $x\in \Spec B$ (this follows, for instance, from  \cite[Proposition 6.8.10]{IrenaHuneke}). Hence 
$\Sing(J,b)=\Sing(I, c)$. Also,  a   sequence of permissible transformations for $(J,b)$ is also    permissible for $(I,c)$, and a resolution of  one  of them induces a resolution of the other (see \ref{TRess_Tpairs_integral} and \ref{Resol_Rees}).

\

On the other hand, if $\diff_{B|k}^{r}$ denotes the 
 (locally free)  $B$-module of $k$-linear differential operators of order at most $r$,  and if 
$$J\subset \diff_{B|k}^{r}(J):=\langle  D(f): f\in J; D\in \diff_{B|k}^{s}, s=0,1,\ldots, r \rangle,$$  then one can check that 
$$\Sing (J,b)=\Sing(\diff_{B|k}^{1}(J),b-1)=\ldots=\Sing(\diff_{B|k}^{b-1}(J),1).$$

As in the previous situation, a sequence of monoidal transformations is permissible for $(J,b)$ if and only if   is also permissible for  $\Sing(\diff_{B|k}^{i}(J),b-i)$ 
for   $i=0,1,\ldots,b-1$ (see Section  \ref{extensions}).  

\

Thus, given a pair $(J,b)$ we have seen how to  find others that {\em codify the same information from the resolution point of view}.  
The question is: is there an effective criterion  to identify all  pairs with this property? This questions leads us to consider other morphisms, that also play a role in resolution problems  
 (in addition to  the permissible monoidal transformations that we have already defined). 

\

 If  $\operatorname{Spec} B=V \stackrel{\rho}{\longleftarrow} \mathfrak{V}$ is the constructive resolution of $(J,b)$ obtained by the use of the resolution functions from (\ref{uppersemi}), and if  $\varphi: \mathbb{V}\to V$ is a smooth morphism,  then it is desirable    that $\mathbb{V} \leftarrow \mathfrak{V}\times \mathbb{V}$ be also  the  constructive resolution of $(\varphi^*(J),b)$ obtained via the same resolution functions (for instance, passing to an open \'etale neighborhood of a point should not affect the output of our arguments).

\

This discussion leads us to the  consideration of   local sequences. A {\em local sequence} over $V=\operatorname{Spec}B$ is a sequence of morphisms 
\begin{equation}
\label{Local_Sequ_Ideals}
\xymatrix@R=2pt@C=45pt{
V=V_0 & V_1 \ar[l]_{\hspace{15pt} \pi_{0}} &    \ldots  \ar[l]_{\pi_{1}} & V_{m} \ar[l]_{\pi_{m-1}}
}
\end{equation}
where each 
$V_i \stackrel{\pi_i}{\longleftarrow} V_{i+1}$  is either the blow up at a smooth center or 
a smooth morphism. A local sequence like (\ref{Local_Sequ_Ideals}) is {\em $(J,b)$-local}  if whenever
 $\pi_i$  is a monoidal transformation with smooth center $Y_i$, then 
$Y_i\subset \Sing (J_i,b)$, and  $J_{i+1}$ is the transform of  $J_i$ as in (\ref{trasformadoHironaka}), and 
whenever 
$\pi_i$  is 
 a smooth morphism the ideal $J_{i+1}$ is the pull-back of $J_i$ in $V_{i+1}$.   

\

We  declare two pairs $(J,b)$ and $(I, c)$ to be {\em equivalent}  if  any local sequence for one is also a local sequence for the other. Moreover, if 

\begin{equation}
\label{Local_Sequ_IdealsB}
\xymatrix@R=2pt@C=25pt{
V=V_0 & V_1 \ar[l]_{\hspace{15pt} \pi_{0}} &    \ldots  \ar[l]_{\pi_{1}} & V_{m} \ar[l]_{\pi_{m-1}}\\
(J,b)=(J_0,b)  & (J_1,b)  &   \ldots    & (J_m,b)\\
(I,c)=(I_0,c)  & (I_1,c)  &   \ldots    & (I_m,c),
}
\end{equation}
is $(J,b)$-local (hence $(I, c)$-local), 
one has that $\Sing(J_i,b)=\Sing(I_i, c)$ for $i=0,1,\ldots, m$. 

\

Once we have defined this equivalence relation among pairs, we formulate our problem:

\begin{quote}
{\bf Goal 1: Canonicity principle.} {\em Our goal is to  provide a canonical representative among all pairs $(I,c)$ that define the same singular locus as $(J,b)$ under  sequences of $(J,b)$-local transformations as (\ref{Local_Sequ_IdealsB}). 
  By exhibiting a canonical representative,   
we mean to give a criterion to be able to 
construct, without ambiguity, a particular pair among all pairs defining the same 
local sequences.}
\end{quote}

\begin{quote}
{\bf Goal 2: Applications.} {\em In combination 
with the theory of elimination of variables   (see  \cite{hpositive}), the canonicity principle will provide a 
useful tool to address Problems \ref{Problem_1} and  \ref{Problem_2}. }
\end{quote}

\

\noindent {\bf The role of  Rees algebras.} 

\

The most natural way to associate a Rees algebra to a given ideal is to 
consider the Rees ring generated by it (i.e., the direct sum of powers of the given ideal). But, when passing from a Rees ring to its integral closure we end up with a Rees algebra that, in general,  is not a Rees ring any more.  This forces us to 
work in a wider context and   consider arbitrary 
(finitely generated) Rees algebras instead. 

\

Assume that  $B$ is  as above. A $B$-{\em Rees algebra} is a finitely generated graded ring say ${\mathcal G}=\bigoplus_{n\in {\mathbb N}}I_nW^n$, where $I_0=B$, $I_n\subset B$ is an ideal for $n\in {\mathbb N}_{\geq 1}$, and $I_k\cdot I_l\subset I_{k+l}$ for $k,l\in{\mathbb N}$. Here $W$ is just a variable that helps us keeping track of the grading. 
It is the grading of ${\mathcal G}$ that   enables us to define an action of the differential operators in a natural way. 

\

In our arguments, we will associate a Rees algebra to a given ideal in a suitable way, so that  
if two ideals have the same integral closure (in the sense of ideals), then the 
Rees algebras associated to them will have the same integral closure as well (in the sense of 
integral closure of Rees algebras). For instance,  given $(J,b)$ as before, 
   we will be interested in the Rees algebra generated by $J$ in degree $b$. If   $(I,c)$ is such that 
 $I^b$ has the same integral closure as $J^c$, then the Rees algebras associated to both pairs will also 
 have the same  integral closure. This assignation defines a  map from pairs $(J,b)$ to Rees algebras: 
 $$\begin{array}{rrcl} 
  \mu: &  \text{Pairs on } B & \longrightarrow &   B\operatorname{-Rees \ algebras} \\
  &  (J,b) &   \longmapsto  & B[JW^b]. 
 \end{array}$$
This map has an inverse (at least if we agree not to distinguish between two Rees algebras 
if they share the same integral closure). In fact, it can be shown that any (finitely generated) Rees algebra is, up to integral closure, the Rees algebra generated by some ideal in some weight (see Lemma \ref{N_adecuado}). This provides a (natural) dictionary    from the class of $B$-Rees algebras to 
the set  of pairs on $B$. 

\

The {\em singular locus of a Rees algebra} $ {\mathcal G}=\bigoplus_{n\in {\mathbb N}}I_nW^n$ is defined as 
$$\Sing{\mathcal G}:=\{x\in \Spec B: \nu_x(I_n)\geq n, \operatorname{ for } n\in {\mathbb N}\}.$$
So $\Sing(J,b)=\Sing B[JW^b]$. 

\

Permissible monoidal transformations and transforms 
can be defined for Rees algebras in the same fashion as for pairs (see Section \ref{RessVSPairs}). In addition, 
${\mathcal G}$-local sequences can be defined as well (see Section \ref{WeakMTheorem}).  This allows us to reformulate  our question about pairs, in the language of 
Rees algebras instead. Denote by ${\mathscr C}_B({\mathcal G})$ the set of all $B$-Rees algebras that 
have the same singular locus of $\mathcal{G}$, and this condition is preserved 
after considering any local sequence. 

\begin{quote}
{\bf Goal.} {\em Characterize all Rees algebras in ${\mathscr C}_B({\mathcal G})$,   and exhibit a canonical representative in this set,  i.e., give a criterion in order to construct, without ambiguity, a unique Rees algebra within the class.}
\end{quote}

 The 
 advantage of working with Rees algebras (instead of pairs) is that there is a natural way to saturate a given Rees algebra ${\mathcal G}$ 
via  the action of differential operators. We refer to this saturated algebra as $\DDiff ({\mathcal G})$.  This saturation is compatible with taking integral closure. In other words, if two Rees algebras have the same integral closure, then 
so do their differential saturations (see \cite[Theorem 6.13]{integraldifferential}). 

\

Given a Rees algebra ${\mathcal G}$ we will denote by $\overline{\mathcal G}$ its integral closure in $B[W]$. With this notation we will prove the following theorem: 

\begin{quote}  {\bf Theorem \ref{CanonicalChoice}. Canonicity Principle.}  {\em Let 
$V$ be a smooth scheme of finite type over a perfect field $k$. 
Let ${\mathcal K}$ and ${\mathcal G}$ be  ${\mathcal O}_V$-Rees algebras. Then  ${\mathcal K} \in {\mathscr C}_V({\mathcal G})$ if and only if 
 $\overline{\DDiff({\mathcal K})}=\overline{\DDiff({\mathcal G})}$. Therefore the 
 differential Rees algebra $\overline{\DDiff({\mathcal G})}$ is the canonical representative of  ${\mathscr C}_V({\mathcal G})$. }
 \end{quote}

This result is closely related to the Finite Presentation Theorem in  \cite{Hironaka05}, where geometric methods are used. However, here we present a proof that entirely relies on techniques coming from commutative algebra. Our proof  is  directed to the Canonicity Principle (already used in \cite{BrV}), and to its applications. 

\

Regarding to Problem \ref{Problem_1},  using Theorem \ref{CanonicalChoice} it is proven that, already for plane curves, 
it is impossible to describe  
the maximum multiplicity locus of a plane curve in one dimension less. This is  a pathology  intrinsic  to the  
positive  characteristic, since in characteristic zero this always can be done using hypersurfaces 
of maximal contact.  In other words, the multiplicity is not an appropriate invariant for inductive arguments in resolution of singularities, and the obstruction appears already in dimension one (see Section \ref{Non_Maximal}). 

\

As for   Problem \ref{Problem_2},  an application of Theorem \ref{CanonicalChoice}   trivializes the  local-global problems that arise in  constructive resolution of singularities.  More precisely,  a  corollary of Theorem \ref{CanonicalChoice} is that the local invariants 
defined in resolution come from  the canonically defined Rees  algebras. Therefore the local-global problem has a (natural) trivial solution. These ideas will be made  precise in the next paragraphs, and will be fully explained in Part \ref{applications}.

\

The paper is organized as follows. Section \ref{ReesAlgebras} is devoted to presenting 
the basics of Rees algebras, while  in Section \ref{RessVSPairs}  the dictionary between Rees algebras and pairs is 
properly stated. The notion of     {\em weak equivalence } on Rees algebras is presented in 
Section \ref{WeakMTheorem}.  Theorems \ref{Canonical} and \ref{CanonicalChoice} are stated in the same section. Section \ref{extensions} is dedicated to studying   the action of differential operators on Rees algebras.  In  Sections \ref{technical_integral_closure} and \ref{technical_closed}  we state some technical results needed for the proof of Theorem \ref{Canonical}, which is finally given in Section \ref{proof_canonical}. 
The last part of the paper, Part \ref{applications}, is devoted to 
describe some applications of   Theorem \ref{Canonical}:
\begin{enumerate}
\item[(i)] The (canonical) assignment of pairs to the maximum stratum of the Hilbert-Samuel function of a given variety (Section \ref{HilbertSamuel}); 
\item[(ii)] The definition of invariants for resolution of singularities (Section \ref{Ejemplos_Invariantes}); this together with 
(i) settles the local-global problem; 
\item[(iii)] We study a feature  of positive characteristic  related to the  the impossibility of using inductive arguments to address  resolution  using  the multiplicity as the unique invariant(see Section \ref{Non_Maximal}).
\end{enumerate}

\

{\em Acknowledgments.} We would like to thank the referee for careful reading this manuscript. His/her useful comments and suggestions have 
helped to improve the presentation of the paper.

\section{Rees algebras} 
\label{ReesAlgebras}
\begin{Definition} \label{Reesalg}{\rm Let $B$ be a Noetherian ring, and let
$\{I_n\}_{n\in {\mathbb N}}$ be a sequence of ideals in $B$ satisfying the
following conditions:
\begin{enumerate}
\item[i.] $I_0=B$;
\item[ii.] $I_k\cdot I_l\subset I_{k+l}$.
\end{enumerate}
The graded subring ${\mathcal G}=\bigoplus_{n\geq 0}I_nW^n$ of
the polynomial ring $B[W]$ is said to be a {\em $B$-Rees algebra}, or a  Rees algebra over $B$,   if
it is a finitely generated $B$-algebra. }
\end{Definition}
A Rees algebra can be described by
giving a finite set of   generators, say
$\{f_{1}W^{n_1},\ldots,f_{s}W^{n_s}\}$,   
$${\mathcal G}=B[f_{1}W^{n_1},\ldots,f_{s}W^{n_s}]\subset B[W]$$ 
with $f_{i}\in B$ for $i=1\ldots,s$.  An element
$g\in I_n$ will be of the form $g=F_n(f_{1},\ldots,f_{s})$ for
some weighted homogeneous polynomial of degree $n$ in $s$-variables,
$F_n(Y_1,\ldots,Y_s)$, where $Y_i$ has weight $n_i$ for
$i=1,\ldots,s$.

\begin{Example} {\rm The  typical  example of a Rees algebra is the  {\em Rees
ring of an ideal}  $J\subset B$,  say   ${\mathcal
G}=B[JW]=\oplus_nJ^nW^n$.   
 }
\end{Example}

\begin{Example}
{\rm Following a slightly different pattern, also the graded algebra ${\mathcal  G}=B[JW^b]=B\oplus 0W\oplus\ldots \oplus 0W^{b-1}\oplus JW^b\oplus 0 W^{b+1}\oplus\ldots \oplus 0W^{2b-1}\oplus J^2W^{2b}\oplus 0W^{2b+1}\oplus \ldots$  is a Rees algebra, which we will refer to as an {\em almost-Rees ring},  and will be denoted by ${\mathcal G}_{(J,b)}$.   Almost-Rees rings will play a central role in our arguments. On the one hand,  as we will see, they are the natural bridge between almost-Rees rings  and pairs (see the Introduction and also Section \ref{RessVSPairs}). 
On the other, any Rees algebra is, up to integral closure, an almost-Rees ring (see \cite[Remark 1.3]{EV}, or  \cite[2.3]{integraldifferential}), and indeed this almost-Rees ring can be chosen in a particular way (see Lemma \ref{N_adecuado}).
So, philosophically, from our point of view, the study of Rees algebras reduces to the understanding of the theory of almost-Rees rings.   }
\end{Example}

The notion of Rees algebra extends to   schemes  in the obvious manner:  a sequence of sheaves of ideals  $\{I_n\}_{n\geq 0}$ on a scheme $V$,  defines a sheaf of  Rees algebras ${\mathcal G}$ over $V$ if  $I_0={\mathcal O}_V$, and  $I_k\cdot I_l\subset I_{k+l}$
for all non-negative integers  $k,l$,  and if there is an affine
open cover $\{U_i\}$  of $V$,  such that ${\mathcal G}(U_i)\subset
{\mathcal O}_V(U_i)[W]$ is an ${\mathcal O}_V(U_i)$-Rees  algebra
in the sense of   Definition \ref{Reesalg}.

\begin{Paragraph}\label{integral_closure_Rees}{\bf Rees algebras and integral closure.} As indicated in the Introduction, from the point of view of resolution, it seems quite natural not to distinguish between two Rees algebras if they 
have the same integral closure. Since integral closure is 
a concept  of local nature, we may assume  to be working on a smooth ring $B$ over a perfect field $k$,  with quotient field $K(B)$. We will be interested in studying 
 the integral closure of 
a  $B$-Rees algebra  ${\mathcal G}\subset B[W]$ in $K(B)[W]$, and will denote it by $\overline{\mathcal G}$. 

The integral closure of a 
$B$-Rees algebra is a $B$-Rees algebra again: on the one hand, 
  \cite[Theorem 2.3.2]{IrenaHuneke} ensures that the integral closure of a 
Rees algebra is a graded ring; 
on the other, the fact that Rees algebras are, by definition, finitely generated  
over an excellent ring,   guarantees  that their integral 
closure is finitely generated too  
 (see \cite[7.8.3.ii), vi)]{EGAIV}). 
\end{Paragraph}

\begin{Paragraph}\label{Veronese}{\bf The Veronese action on a Rees algebra 
\cite[2.3]{integraldifferential}, \cite[4.3.1]{hpositive}.} Given a natural number $M$, the  {\em $M$-th Veronese} action on   a Rees algebra 
${\mathcal G}=B\oplus I_1W\oplus I_2W^2\oplus\ldots\oplus I_nW^n\oplus \ldots$, 
is defined as 
$${\mathbb V}_M({\mathcal G}):=\bigoplus_{k\geq 0}I_{Mk}W^{Mk}.$$ 
Since ${\mathcal G}$ is finitely generated,   ${\mathbb V}_M({\mathcal G})\subset {\mathcal G}$ is a finite extension 
for any choice of  $M$. 
\end{Paragraph}

\begin{Remark}\label{integral_in_WM} Observe that 
${\mathbb V}_M(\overline{{\mathcal G}})$ is the integral closure of ${\mathbb V}_M({{\mathcal G}})$ in
$B[W^M]$. On the one hand, the inclusion 
${\mathbb V}_M(\overline{{\mathbb V}_M({\mathcal G})})\subset {\mathbb V}_M(\overline{\mathcal G})$ is clear. To check 
the other inclusion, we argue as follows.  Let  $f\in K(B)$,   
and let $n$ be a positive integer. If $fW^n\in \overline{\mathcal G}$ then it satisfies an integral relation of the form $$(fW^n)^\ell+a_1W^{n}(fW^n)^{\ell-1}+\dots +a_\ell W^{n\ell}=0$$
were each $a_i W^{ni} (\in \mathcal G)$ is homogeneous of degree $n \cdot i$. In particular, if $n$ is a multiple of $M$, then $fW^n\in B[W^M]$ and   $a_i\in 
{\mathbb V}_M({{\mathcal G}})$ for $i=1,\ldots,\ell$. Therefore 
$fW^n$ is in the integral closure of $ {\mathbb V}_M({\mathcal G})$  and also  in $B[W^M]$,  so 
${\mathbb V}_M(\overline{\mathcal G})\subset  {\mathbb V}_M(\overline{{\mathbb V}_M({\mathcal G})})$. Therefore, ${\mathbb V}_M(\overline{\mathcal G})= {\mathbb V}_M(\overline{{\mathbb V}_M({\mathcal G})})=\overline{{\mathbb V}_M({\mathcal G})}\cap B[W^M]$, and hence ${\mathbb V}_M(\overline{{\mathcal G}})$ is the integral closure of ${\mathbb V}_M({{\mathcal G}})$ in
$B[W^M]$. In what follows, and unless 
otherwise indicated,  we will always be 
considering the integral closure of a Rees algebra in $B[W]$. 
\end{Remark}

\begin{Example}
\label{IntegralReesRing} 
\cite[Proposition 5.2.1]{IrenaHuneke} 
If  ${\mathcal G}$ is the Rees ring of an ideal 
$J$, then its integral closure,  $\overline{\mathcal G}$, is ${B}\oplus \overline{J}W\oplus\overline{J^2}W^2\oplus\ldots\oplus 
\overline{J^n}W^n\oplus \ldots$,
where $\overline{J}$ denotes the integral closure 
of the ideal $J$ in $B$. 
\end{Example}

\begin{Example}\label{InClosure_ring_almost} If ${\mathcal G} =B[JW^b]$ is an almost-Rees ring,  then $\overline{\mathcal G}=B\oplus I_1W\oplus 
\ldots \oplus I_bW^b\oplus\ldots$, where each $I_n$ is integrally closed and, by Remark \ref{integral_in_WM} and Example \ref{IntegralReesRing},  
$I_{bk}=\overline{J^k}$ for all $k\in {\mathbb N}$.
\end{Example}

\begin{Lemma}\label{N_adecuado} For suitable choices of $N$, the Rees algebras  $\mathbb{V}_N(\mathcal{G})$ are    almost-Rees rings. In particular, any Rees algebra is finite over an almost-Rees ring. 
\end{Lemma}

\noindent {\em Proof:} 
Assume that ${\mathcal G}$  is generated by  $f_1W^{n_1}, \ldots, f_sW^{n_s}$, and let 
 $M$ be a common multiple of all $n_i$. Since $f_i^{M/n_i}W^M\in I_MW^M$,  we obtain the following finite extension of (graded) algebras:
$$B[I_MW^M]=\bigoplus_{n\geq 0}(I_M)^nW^{nM}\subset \mathbb{V}_{M}(\mathcal{G})=\bigoplus_{n\geq 0}I_{nM}W^{nM} \subset \mathcal{G}.$$
If $I_M ^n= I_{nM}$   for all $n\geq 1$  the proof is  complete. 
Otherwise, we can assume that $\mathbb{V}_M(\mathcal{G})$ is generated as $B[I_MW^M]$-module by a finite set of homogeneous elements. Moreover, these generators can be 
chosen with degrees at most $AM$, for some non-negative integer $A$.  Then,  $I_{nM}=(I_M)^{(n-A)}I_{AM}$ for all $n \geq A$  (see \cite[Proposition 5.2.5]{IrenaHuneke}). Now,  replacing $n$ 
 by 
$kA$ in the previous equality,  we obtain
$$\left(I_{AM}\right)^k \subset I_{kAM} = (I_M)^{A(k-1)} \cdot I_{AM} \subset \left(I_{AM}\right)^{k-1}\cdot I_{AM}=\left(I_{AM}\right)^{k}
$$
Thus  $(I_{AM})^k=I_{kAM}$.  This shows that $\mathbb{V}_{AM}(\mathcal{G})$ is an almost-Rees ring. The same holds for all $\mathbb{V}_N(\mathcal{G})$ if $N$ is a multiple of $AM$.\qed

\begin{Remark}\label{potencias_enteras} By Lemma \ref{N_adecuado}, for any integrally closed Rees algebra ${\mathcal G}$ there are  infinitely many choices of $N$ for which ${\mathbb V}_N({\mathcal G})={\mathcal O}_V[I_NW^N]$  is such that in  addition all powers  of $I_N$ are integrally closed. 
\end{Remark}

Let  ${\mathcal G}$ and ${\mathcal K}$ be two Rees algebras, and let $N$ be a non-negative integer so that 
both $\mathbb{V}_N({\mathcal G})$ and $\mathbb{V}_N({\mathcal K})$ are almost-Rees rings.  Then,  by \ref{Veronese}, $\overline{\mathcal G}=\overline{\mathcal K}$ if and only if  $\overline{\mathbb{V}_N({\mathcal G})}=\overline{\mathbb{V}_N({\mathcal K})}$.  The following lemma gives a useful criterion to compare the integral closure of two almost-Rees rings.

\begin{Lemma}\label{almostIntegral} Let  $B[JW^b]$ and $B[IW^c]$ be two almost-Rees rings. Then 
$B[JW^b]\subset \overline{B[IW^c]}$ if and only if  $J^c\subset \overline{I^b}$. Thus, $\overline{B[JW^b]}= \overline{B[IW^c]}$
 if and 
only if $\overline{J^c}=\overline{I^b}$. 

\end{Lemma}

\noindent{\em Proof:}  The ``only if'' part is easy, since just by looking at the degree $bc$ piece  in both $B[JW^b]$ and $\overline{B[IW^c]}$ we obtain that  $J^c  \subset \overline{I^b}$.
For the converse, since $B[JW^b]$ and $\mathbb{V}_{bc}(B[JW^b])$ have the same integral closure in B[W] we only need to prove that $\mathbb{V}_{bc}(B[JW^b])\subset \overline{B[IW^c]}$, which 
is clear since  $J^c\subset \overline{I^b}$ implies that  $J^{ck} \subset (\overline{I^b})^k \subset \overline{I^{bk}}$. \qed

\

In  this paper we will come across  some information, provided by a Rees algebra, that will be preserved 
up to integral closure.  This is the case of the  concepts to be defined below such as the  
{\em singular locus of an algebra},   the {\em order of a Rees algebra  at a point}, and the 
{\em zero  set of a Rees  algebra}.  We will see that if two algebras have the same integral closure, then they share 
the same singular locus, the order at a given point is the same, and they both have identical zero set (see \ref{SingOrdIntegral}). 
Considering 
Rees algebras up to integral closure  will allow us to 
define the normalized blow up of a Rees algebra (see \ref{normalblowup}), which plays a role in the 
proof of Theorem \ref{Canonical} (see \ref{ProofTheo}).

\begin{Paragraph} \label{singularlocus}{\bf The singular locus  of a Rees algebra.}  \cite[1.2]{positive}
 Let $V$ be a smooth scheme over a perfect field $k$,   and let
${\mathcal G}=\oplus_nI_nW^n$ be a sheaf of ${\mathcal O}_V$-Rees 
algebras.  Then  the  {\em singular locus of ${\mathcal
G}$},  $\Sing{\mathcal G}$,  is 
$$\Sing{\mathcal G}:=\bigcap_{n\in {\mathbb N}_{>0}}\{x\in V: \nu_x(I_n)\geq n\},$$
where   
$\nu_x(I_n)$ denotes the usual  order of $I_n$ in the regular local
ring ${\mathcal O}_{V,x}$. Observe that $\Sing{\mathcal G}$  is a closed subset in 
  $V$. If ${\mathcal G}$ is generated by $f_1W^{n_1},\ldots,f_sW^{n_s}$ on an affine open set $U\subset V$, then it can be shown that 
\begin{equation} \label{sing_generators}\Sing{\mathcal G}\cap U =\bigcap_{i=1}^s\{x\in V: \nu_x(f_i)\geq n_i\} \end{equation}
(see \cite[Proposition 1.4]{EV}).   
\end{Paragraph}

\begin{Paragraph} \label{orderRees}{\bf The order of a Rees algebra at a point.}
\cite[6.3]{EV} Let ${\mathcal
G}=\bigoplus_{n\geq 0}I_nW^n$ be a Rees algebra on a smooth scheme $V$,  let $x\in \Sing{\mathcal
G}$, and assume that $fW^n\in I_nW^n$ in an open neighborhood of $x$. Then set
$$\ord_{f
W^n 
}(x)=\frac{\nu_x(f)}{n}\in {\mathbb Q},$$
where, as before, $\nu_x(f)$ denotes the order of $f$ in the regular local ring
${\mathcal O}_{V,x}$. Notice that $\ord_x(f)\geq 1$ since
$x\in \Sing{\mathcal G}$. Now define
$$\ord_{\mathcal
G}(x)=\inf_{n\geq 1} \left\{\frac{\nu_x(I_n)}{n} \right\}.$$ If
${\mathcal G}$ is generated by
$\{f_{1}W^{n_1},\ldots,f_{s}W^{n_s}\}$ on an affine neighborhood of $x$ then it can be shown that 
\begin{equation}\label{order_generators}\ord_{\mathcal G}(x)=\min\{\ord_{f_{i}
W^{n_i}}(x):
i=1,\ldots,s\},\end{equation} and therefore, since
 $x\in \Sing{\mathcal
G}$, $\ord_{\mathcal G}(x)$ is a rational number that is greater
than or equal to one.

\end{Paragraph}

\begin{Paragraph}{\bf The zero set  of a Rees algebra.}\label{zeros} Given an  
$\mathcal{O}_V$-Rees algebra, $\mathcal{G}=\bigoplus_{n\geq 0} I_nW^n$, we define the 
\emph{zero set of $\mathcal{G}$}, $V(\mathcal{G})$,  as 
  $$V(\mathcal{G}):= \bigcap_{n\geq 1}V(I_n).$$
This closed set satisfies the following   properties:
\begin{itemize}
\item If $\mathcal{G}$ is an almost-Rees algebra, say 
$\mathcal{O}_V[JW^b]$, then $V(\mathcal{G})=V(J)$.
 \item  In general, $\Sing\mathcal{G} \subsetneqq V(\mathcal{G}).$ 
\end{itemize}
\end{Paragraph}

\begin{Example} \label{ExampleA} 
Let $\mathbb{H}\subset V$ be a hypersurface, and let $b$ be a non-negative 
integer. Then the singular locus of the Rees algebra generated by 
$I(\mathbb{H})$ in degree $b$, i.e., the singular locus of $\mathcal{G}_{(I(\mathbb{H}),b)}=\mathcal{O}_V[I(\mathbb{H})W^b](\subset \mathcal{O}_V[W])$,
  is the closed set of points of multiplicity at
least $b$ of $\mathbb{H}$ (which may be empty). The order of $\mathcal{G}_{(I(\mathbb{H}),b)}$ at a point in the singular locus is the multiplicity of $\mathbb{H}$ divided by $b$, and the zero set of ${\mathcal G}$ is $\mathbb{H}$. 
\end{Example}
\begin{Example}
\label{ExampleB}
In the same manner, if 
$J\subset {\mathcal O}_V$ is an arbitrary non-zero sheaf of ideals, and $b$ is a
non-negative integer,    then 
$\Sing \mathcal{G}_{(J,b)}$ consists of the points of $V$ where the order of $J$ is at least $b$,
$\ord_{\mathcal{G}_{(J,b)}}(x)= \frac{\nu_x(J)}{b}$   for all $x\in \Sing \mathcal{G}_{(J,b)}$, 
 and $V(\mathcal{G})=V(J)$. 
\end{Example}

\begin{Paragraph} \label{SingOrdIntegral}
{\bf Singular locus, order, zero set, and integral closure.}
Two Rees algebras with the same integral closure have the same singular locus \cite[Proposition
4.4 (1)]{integraldifferential}, the same order at a point \cite[Proposition 6.4 (2)]{EV}, and the same zero set. 
In the following lines we sketch a proof of these facts for self-containment.

\end{Paragraph}

\begin{proof}[Sketch of the proof.]
We will argue locally on an affine open set, say $U=\Spec(B)$.
Assume that $\mathcal{G}$ and $\mathcal{K}$ are two Rees algebras with the same integral closure. By Lemma \ref{N_adecuado}, one can choose an appropriated non-negative integer $N$ so that  both, $\mathbb{V}_N(\mathcal{G})$ and $\mathbb{V}_N(\mathcal{K})$, are  almost-Rees rings. Note that the  four algebras have the same integral closure. Then we only need to consider the following  two cases:

\

\noindent {\bf Case 1.}  Let $\mathcal{G}=\bigoplus_nI_nW^n=B[f_1W^{n_1}, \ldots, f_sW^{n_s}]$ and let $\mathbb{V}_N(\mathcal{G})=B[I_NW^N]=\bigoplus_{n} I_{nN}W^{nN}$:

Noticing that for any generator $f_iW^{n_i}$ of $\mathcal{G}$, $\frac{\nu_x(f_i)}{n_i}=\frac{\nu_x(f^{N/n_i}_i)}{N}$, one has that
\begin{equation*}\Sing{\mathcal G}=\bigcap_{i=1}^s\{x\in U: \nu_x(f_i^{N / n_i})\geq N\}=\{x\in U: \nu_x(I_N)\geq N\}=\Sing \mathbb{V}_N(\mathcal{G}), \end{equation*}
and for $x\in \Sing{\mathcal G}=\Sing \mathbb{V}_N(\mathcal{G})$,
\begin{equation*}\ord_{\mathcal G}(x)=\min\{\ord_{f_{i}^{N /n_i}W^N}(x):
i=1,\ldots,s\}=\frac{\nu_x(I_N)}{N}=\ord_{\mathbb{V}_N(\mathcal{G})}(x).\end{equation*}
Moreover, by
the definition of Rees algebras in \ref{Reesalg} and by the choice of $N$, note that $I_n^N W^{nN}\subset I_{nN}W^{nN}=[\mathbb{V}_N(\mathcal{G})]_{nN}$ for all $n\geq 1$, and hence
\begin{equation*}V(\mathcal{G})=\bigcap_{n\geq 1}V(I_n^N)\supset \bigcap_{n \geq 1} V(I_{nN}) = V(\mathbb{V}_N(\mathcal{G})) \supset V(\mathcal{G}). \end{equation*}

\

\noindent {\bf Case 2.} Let $I, J\subset B$ be two ideals, and assume  that $B[JW^n]$ and $B[IW^n]$ have the same integral closure. Note  that Lemma \ref{almostIntegral}  (or Example \ref{InClosure_ring_almost}) guarantees that $\overline{I}=\overline{J}$. Hence $\nu_x(I)=\nu_x(J)$ for all $x\in U$ by \cite[Proposition 6.8.10]{IrenaHuneke}, and $V(J)=V(I)$ (since $I\subset \overline{I}=\overline{J}\subset \sqrt{I}$).

\end{proof}

\begin{Paragraph}\label{normalblowup}{\bf The normalized blow up of a Rees algebra.}
We can associate a blow up to any   Rees algebra ${\mathcal G}=\oplus_nI_nW^n$ 
in the following way. First choose an appropriate $N$ so that the $N$-th Veronese action on  ${\mathcal G}$  
 is an almost-Rees ring, say  ${\mathbb V}_N({\mathcal G})=B[I_NW^{N}]$ (see Lemma \ref{N_adecuado}). Next   
 consider the blow up of $\Spec B$ at $I_N$, i.e.,
$$\xymatrix@R=2pt@C=20pt{
\Spec B &\ar[l] \operatorname{Proj}(B[I_NW^N])
}$$
and define $$\Bl({\mathcal G}):=\operatorname{Proj}(B[I_NW^N]).$$ 
It turns out that this construction 
is independent of the choice of $N$: if $M$ is so that 
 ${\mathbb V}_{M}({\mathcal G})$ is also an almost-Rees ring, say ${\mathbb V}_{M}({\mathcal G})=B[I_{M}W^{M}]$, then 
$(I_N)^{M}=(I_{M})^N$, and therefore the blow up of 
$\Spec B$ at $I_N$  coincides with 
 that  of $\Spec B$ at $I_{M}$. In this way we attach to $\mathcal{G}$ another invariant: the {\em normalized blow up of $\Spec B$ with respect to ${\mathcal G}$},
$$
\xymatrix@R=2pt@C=20pt{
\Spec B & \ar[l]  \Bl(\mathcal{G}) & \ar[l] \overline{\Bl(\mathcal{G})}:=\overline{\operatorname{Proj}(B[I_NW^N])}.
}
$$
Using  Remark \ref{integral_in_WM}, Example \ref{InClosure_ring_almost}, and the previous discussion, it can be checked that $\overline{\Bl({\mathcal G})}=\Bl(\overline{\mathcal G}).$
\end{Paragraph}

\section{Rees algebras 
 and 
pairs}
\label{RessVSPairs}
In this section we will explore the connection between Rees algebras and pairs. Our ultimate  interest   is to  reformulate the problem of constructive log-resolution of ideals  in terms of  Rees algebras.
 
\begin{Paragraph}\label{ordersingintegral}{\bf The bridge between Rees algebras and pairs.} \cite[Proposition
4.4 (1)]{integraldifferential}, \cite[Proposition 6.4]{EV}
Rees algebras can be related to pairs in the   following way.  On the one hand, we can attach an almost-Rees ring to a given 
 pair $(J,b)$, say ${\mathcal G}_{(J,b)}=\mathcal{O}_V[JW^b]$. On the other hand,  
  as indicated in Lemma \ref{N_adecuado},  up to integral closure, any Rees algebra   is 
an almost-Rees ring. In other words, for every Rees algebra ${\mathcal G}=\bigoplus_{n\geq 0}I_nW^n$ there is some $N$ such that 
$\mathbb{V}_N({\mathcal G})=\mathcal{O}_V[I_NW^N]=\mathcal{G}_{(I_N,N)}$ is an almost-Rees ring, and it  therefore  
 can be interpreted as the Rees algebra associated to the pair $(I_N,N)$. Moreover, by   \ref{SingOrdIntegral} and   
 the definition of the singular locus of a pair  (see (\ref{SingHironaka}), and Example \ref{ExampleB}),  one has that 
 $$\Sing{\mathcal G}=\Sing{\mathbb V}_N({\mathcal G})=\Sing(I_N,N).$$ 
In addition,   
by \ref{SingOrdIntegral} and by the definition of  Hironaka's order function for a pair  
(see (\ref{ordenHironaka})), one has that for any $x \in \Sing \mathcal{G}$
$$\ord_{\mathcal G}(x)=\ord_{{\mathbb V}_N({\mathcal G})}(x)=
\ord_{(I_N,N)}(x).$$

\end{Paragraph}

This parallelism 
between Rees algebras and pairs can be carried over one step further, since resolution 
 can also be formulated for Rees algebras as explained in the following paragraphs.

\begin{Paragraph} \label{weaktransforms} {\bf Transforms of Rees algebras by blow ups.}
A smooth closed subscheme $Y\subset V$ is said to be {\em permissible} for ${\mathcal G}=\bigoplus_{n}J_nW^n\subset{\mathcal O}_{V}[W]$ if     $Y\subset
\Sing {\mathcal G}$.  A {\em permissible monoidal transformation}  is the blow up at a permissible center,  $V\leftarrow V_1$. If 
$H_1\subset V_1$ denotes  the exceptional divisor, then for each
$n\in {\mathbb N}$,
$$J_n{\mathcal O}_{V_1}=I(H_1)^nJ_{n,1}$$
for some  sheaf of ideals $J_{n,1}\subset {\mathcal
O}_{V_1}$.  The   {\em  transform of ${\mathcal
G}$}  in $V_1$ is then defined as: 
$${\mathcal G}_1:=\bigoplus_{n}J_{n,1}W^n;$$
and for a given homogeneous element $fW^m\in {\mathcal G}$, a  
{\em weighted 
transform}, $f_1W^m\in {\mathcal G}_1$,  is defined by choosing any generator of the principal ideal 
$$I(H_1)^{-m}\cdot \langle f\rangle {\mathcal O}_{V_1}.$$

\end{Paragraph}

The next proposition gives a local description of the 
transform of a Rees algebra  after a permissible
monoidal transformation.
\begin{Proposition}\label{localt}\cite[Proposition 1.6]{EV}
 Let ${\mathcal G}=\bigoplus_nI_nW^n$ be a
Rees algebra on a smooth scheme $V$   over a field $k$, and let
$V\leftarrow V_1$ be a permissible monoidal transformation. Assume, for simplicity, that 
$V$ is affine. If  
${\mathcal G}$ is generated by
$\{f_{1}W^{n_1},\ldots,f_{s}W^{n_s}\}$, then its  transform ${\mathcal
G}_1$ is generated by
$\{f_{1,1}W^{n_1},\ldots,f_{s,1}W^{n_s}\}$,
where   $f_{i,1}$ denotes a weighted transform of $f_{i}$ in $V_1$ 
for  $i=1,\ldots,s$.
\end{Proposition}

\begin{Paragraph}{\bf Transforms of Rees algebras, 
transforms of pairs, and integral extensions.}\label{TRess_Tpairs_integral}
We emphasize that the suitable integer, $N$, that links Rees algebras and pairs passing through an almost-Rees ring, is preserved by permissible monoidal transformations (see \ref{ordersingintegral} above). More precisely, let  $\mathcal{G}=\bigoplus_nJ_nW^n$ be a Rees algebra. If $\mathbb{V}_N(\mathcal{G})$ is an almost-Rees ring, 
i.e.,  if $\mathbb{V}_N(\mathcal{G})=\mathcal{G}_{(J_N,N)}$ for some pair $(J_N,N)$, then observe that 
\begin{equation}\label{transform_veronese}  \mathbb{V}_N(\mathcal{G})_1=\mathbb{V}_N(\mathcal{G}_1). \end{equation}
Hence  
 $\mathbb{V}_N(\mathcal{G}_1)$  is also an almost-Rees ring, and,  moreover,  
  $\mathbb{V}_N(\mathcal{G}_1)=\mathcal{G}_{(J_{N,1},N)}$ where $(J_{N,1},N)$ is the transform of the pair $(J_N,N)$ 
  by the permissible transformation $V\leftarrow V_1$ (see (\ref{trasformadoHironaka}) in the Introduction). This shows that, if  $\mathbb{V}_N(\mathcal{G})$ is an almost Rees ring, then 
  so is $\mathbb{V}_N(\mathcal{G}_1)$ for the same $N$. 
  
\

Using a similar  argument,   it is easy to prove that the transforms of two Rees algebras with the same integral closure also have the same integral closure: suppose that $\mathcal{G}=\bigoplus_n J_nW^n$
    and $\mathcal{K}=\bigoplus_n I_nW^n$  have the same 
 integral closure. By Lemma  \ref{N_adecuado}  there is  a suitable $N$ such that both  $\mathbb{V}_N(\mathcal{G})=\mathcal{G}_{(J_N,N)}$ and $\mathbb{V}_N(\mathcal{K})=\mathcal{G}_{(I_N,N)}$ are almost-Rees rings. So, by Lemma \ref{almostIntegral}, $\overline{J_N}=\overline{I_N}$ as sheaves of  ideals in $V$, and hence    $\overline{J_N\mathcal{O}_{V_1}}=\overline{I_N\mathcal{O}_{V_1}}$. One can now check that also  $\overline{J_{N,1}} = \overline{I_{N,1}}$. To conclude, note that by Lemma \ref{almostIntegral} and (\ref{transform_veronese}), the Rees algebras $\mathcal{G}_1$ and $\mathcal{K}_1$ have the same integral closure.
    
\end{Paragraph}

\begin{Paragraph}\label{Resol_Rees}{\bf Resolution of Rees algebras
 and 
resolution of pairs.}  {\rm A {\em resolution} of a Rees algebra  ${\mathcal G}$  on a smooth scheme $V$ is a finite sequence of blowing ups, 
\begin{equation}
\label{ResolReesAlg}
\xymatrix@R=2pt@C=20pt{
V=V_0 & V_1 \ar[l]_{\hspace{20pt} \rho_{0}} &    \ldots  \ar[l]_{\rho_{1}} & V_{n} \ar[l]_{\rho_{n-1}}\\
{\mathcal G}={\mathcal G}_0  & {\mathcal G}_1 & \ldots & {\mathcal G}_n
}
\end{equation}
at permissible centers $Y_i\subset \Sing{\mathcal G}_i$ for $i=0,1,\ldots,n-1$, so that 
\begin{itemize}
\item $\Sing{\mathcal G}_n=\emptyset$, and
\item the 
 exceptional locus of the composition $V\leftarrow V_n$ is a union of smooth hypersurfaces having only normal crossings in $V_n$.
\end{itemize}

\

 Recall  that the   transformation law under 
permissible transformations is compatible  for both, Rees algebras and pairs  (see \ref{weaktransforms}). 
Hence,  if $\mathbb{V}_N({\mathcal G})={\mathcal G}_{(J,b)}$ for some pair $(J,b)$, then 
a resolution of ${\mathcal G}$ as in (\ref{ResolReesAlg}) gives a resolution  of 
$(J,b)$ (see the Introduction for details on resolutions of pairs).

\

Note, in addition, that the definition of an arbitrary sequence of    transforms  over $\mathcal{G}$ is equivalent to a sequence of transforms of the pair $(J,b)$, 

\begin{equation}
\label{ResolPairs}
\xymatrix@R=2pt@C=20pt{
V=V_0 & V_1 \ar[l]_{\hspace{20pt} \pi_{0}} &    \ldots  \ar[l]_{\pi_{1}} & V_{n} \ar[l]_{\pi_{n-1}}\\
{\mathcal G}={\mathcal G}_0  & {\mathcal G}_1 & \ldots & {\mathcal G}_n\\
(J,b)=(J_0,b) & (J_1,b) & \ldots & (J_n,b).
}
\end{equation}
and vice versa. This is due to the fact that $\mathbb{V}_N({\mathcal G}_i)=\mathbb{V}_N({\mathcal G})_i={\mathcal G}_{(J_i,b)}$ for $i=0,1,\ldots,n$.} 
\end{Paragraph}

\section{Weak equivalence and Main Theorem}
\label{WeakMTheorem}
Paralleling the resemblance between Rees algebras and pairs, 
here we introduce an equivalence relation among 
Rees algebras:  {\em weak equivalence} (see the discussion in Introduction  for a notion of equivalence within pairs). This relation works  in such a way that algebras associated to equivalent pairs, in the sense of Hironaka, 
will also be weakly equivalent. Two Rees algebras  that are weakly equivalent   will  have the same resolution invariants, and 
hence the same constructive resolutions.   

\

Weak equivalence is  defined by taking into 
account  a  {\em  tree of closed sets}  determined by the singular locus of a Rees algebra, ${\mathcal G}$,   and  the singular loci of transforms of ${\mathcal G}$ under suitable morphisms (see Definition \ref{local_sequenceC} below).

\begin{Definition}
\label{local_sequenceA} 
{\rm  Let $V$ be a smooth scheme over a perfect field $k$. A {\em local sequence over   $V$}  is a sequence of the form 
$$\xymatrix@R=2pt@C=30pt{
V=V_0 & V_1 \ar[l]_{\hspace{20pt} \pi_{0}} &    \ldots  \ar[l]_{\pi_{1}} & V_{m} \ar[l]_{\pi_{m-1}}
}$$
where  for $i=0,1,\ldots,m-1$, each 
 $\pi_i$  is either 
the blow up at a smooth closed subscheme, or a smooth morphism. }
\end{Definition}

\begin{Definition} 
\label{local_sequenceB} 
{\rm If ${\mathcal G}$ is an ${\mathcal O}_{V}$-Rees algebra,  a 
  {\em ${\mathcal G}$-local sequence over $V$}   is a local sequence over $V$, 
\begin{equation}\label{ABCranso1}
\xymatrix@R=2pt@C=30pt{
(V=V_0,\mathcal{G}=\mathcal{G}_0) & \ar[l]_{\hspace{20pt} \pi_{0}}  (V_1,\mathcal{G}_1) & \ar[l]_{\hspace{20pt} \pi_{1}}  \cdots & \ar[l]_{\pi_{m-1}} (V_m,\mathcal{G}_m),
}
\end{equation}
where  for $i=0,1,\ldots,m-1$  each
 $\pi_i$ 
 is either a permissible monoidal transformation for 
${\mathcal G}_i\subset {\mathcal O}_{V_i}[W]$ (and then ${\mathcal G}_{i+1}$ is the transform of ${\mathcal G}_{i}$ 
in the sense of  \ref{weaktransforms}), 
or a smooth morphism (and then ${\mathcal G}_{i+1}$ is   the pull-back of 
${\mathcal G}_i$ in $V_{i+1}$).   }
\end{Definition}

\begin{Definition}
{\rm Let ${\mathcal G}$ be an ${\mathcal O}_{V}$-Rees algebra, and let   
\begin{equation}\label{ABCranso12}
\xymatrix@R=2pt@C=30pt{
(V=V_0,\mathcal{G}=\mathcal{G}_0) & \ar[l]_{\hspace{20pt} \pi_{0}}  (V_1,\mathcal{G}_1) & \ar[l]_{\hspace{20pt} \pi_{1}}  \cdots & \ar[l]_{\pi_{m-1}} (V_m,\mathcal{G}_m),
}
\end{equation}
be a  ${\mathcal G}$-local sequence over $V$. Then the collection of closed subsets 
$$\begin{array}{cccc}
\Sing{\mathcal G}_0\subset V_0, & \Sing{\mathcal G}_1\subset V_1, & \ldots, &\Sing{\mathcal G}_N\subset V_N
\end{array}$$
determined by the ${\mathcal G}$-local sequence (\ref{ABCranso12}) is a {\em branch of closed 
subsets over $V$ determined by ${\mathcal G}$}. The union of all branches of closed subsets  
   by considering all ${\mathcal G}$-local sequences over $V$ is 
the {\em tree of closed subsets over $V$ determined by ${\mathcal G}$}. We will   denote it by ${\mathscr F}_V({\mathcal G})$.}
\end{Definition}

\begin{Definition}
\label{local_sequenceC} 
{\rm If ${\mathcal G}$ and ${\mathcal K}$ are 
 two ${\mathcal O}_{V}$-Rees algebras, then a  ${\mathcal G}$-${\mathcal K}$-{\em local sequence over $V$} is a local 
 sequence over $V$ that is both ${\mathcal G}$-local and ${\mathcal K}$-local. 
}
\end{Definition}

\begin{Definition} \label{Defweak} {\rm 
Two Rees algebras, $\mathcal{G}$ and $\mathcal{K}$,  
  are said to be {\em
weakly equivalent} if: 

(i) $\Sing\mathcal{G}=\Sing \mathcal{K}$;  

(ii) Any local ${\mathcal G}$-local sequence over $V$ induces a  ${\mathcal K}$-local sequence over $V$, 
 and any   ${\mathcal K}$-local sequence over $V$ induces a  ${\mathcal G}$-local sequence over $V$; 

(iii) For  any ${\mathcal G}$-local  sequence   over   $V$,  
\begin{equation*}
\xymatrix@R=2pt@C=30pt{
(V, \mathcal{G})= (V_0,\mathcal{G}_{0}) & \ar[l]_{\hspace{20pt} \pi_{0}}  (V_1,\mathcal{G}_{1}) & \ar[l]_{\hspace{20pt} \pi_{1}}  \cdots & \ar[l]_{\hspace{-10pt} \pi_{m-1}} (V_m,\mathcal{G}_{m}),
}
\end{equation*}
and the corresponding  induced ${\mathcal K}$-local over $V$,
\begin{equation*}
\xymatrix@R=2pt@C=30pt{
(V, \mathcal{G})= (V_0,\mathcal{K}_{0}) & \ar[l]_{\hspace{20pt} \pi_{0}}  (V_1,\mathcal{K}_{1}) & \ar[l]_{\hspace{20pt} \pi_{1}}  \cdots & \ar[l]_{\hspace{-10pt} \pi_{m-1}} (V_m,\mathcal{K}_{m}),
}
\end{equation*}
there is an equality of closed sets,  
$\Sing (\mathcal{G}_{j})=\Sing (\mathcal{K}_{j})$
for $0\leq j \leq m$; and vice versa.  
 
}
\end{Definition}

\begin{Remark}
Weak equivalence is an equivalence relation within the class of 
Rees algebras defined over $V$. We denote by ${\mathscr C}_V({\mathcal G})$ the equivalence class 
of an ${\mathcal O}_V$-Rees algebra ${\mathcal G}$. By definition  two Rees algebras are weakly equivalent when 
 they determine the same tree of closed subsets over $V$; i.e.,  two Rees algebras over $V$, say ${\mathcal G}$ and ${\mathcal K}$, are weakly equivalent when ${\mathscr F}_V({\mathcal G})={\mathscr F}_V({\mathcal K})$.
 \end{Remark}

 \begin{Example} 
 \label{localintegral} 
 If two ${\mathcal O}_V$-Rees algebras  ${\mathcal G}$ and ${\mathcal K}$ have the same integral closure, then 
 by \ref{SingOrdIntegral},  $\Sing{\mathcal G}=\Sing{\mathcal K}$. If $\varphi: V_1\rightarrow  V$ is a 
smooth morphism, then  $\varphi^*({\mathcal G})$ and 
$\varphi^*({\mathcal K})$  also  have the same integral closure in $V_1$, and therefore 
$\Sing \varphi^*({\mathcal G})=\Sing \varphi^*({\mathcal K})$.   
Moreover,  if   $V\leftarrow V_1$ is a permissible 
transformation with   center   $Y\subset \Sing{\mathcal G}=\Sing
{\mathcal K}$, then by \ref{TRess_Tpairs_integral}
    the   transforms  of ${\mathcal G}$ and ${\mathcal K}$ in $V_1$, say ${\mathcal
G}_1$ and $ {\mathcal K}_1$,  also have   same integral closure and therefore 
$\Sing{\mathcal
G}_1=\Sing{\mathcal
K}_1$. This already shows that  two Rees algebras 
with  the same integral closure  are weakly equivalent. 
\end{Example}

Given two Rees algebras $\mathcal{G}$ and $\mathcal{K}$ we
 can make sense of the expression ${\mathscr F}_V({\mathcal G})\subset {\mathscr F}_V({\mathcal K})$ in a natural way:

\begin{Definition}
\label{definition_inclusion}
{\rm Let $ {\mathcal K}$ and  ${\mathcal G}$ be two Rees algebras on $V$. 
We will say that 
$${\mathscr F}_V({\mathcal K})\subset {\mathscr F}_V({\mathcal G})$$ if 
$\Sing {\mathcal K}\subset \Sing {\mathcal G}$, and 
any  $ {\mathcal K}$-local sequence over $V$,
\begin{equation*}
\xymatrix@R=2pt@C=20pt{
(V,\mathcal{K})=(V_0,\mathcal{K}_{0}) & \ar[l]_{\hspace{20pt} \pi_{0}}
(V_1,\mathcal{K}_{1}) & \ar[l]_{\hspace{20pt} \pi_{1}} \cdots & \ar[l]_{\hspace{-10pt}\pi_{m-1}} 
(V_m,\mathcal{K}_{m}),
}
\end{equation*}
induces a ${\mathcal G}$-local sequence over $V$,
\begin{equation*}
\xymatrix@R=2pt@C=20pt{
(V,\mathcal{G})=(V_0,\mathcal{G}_{0}) & \ar[l]_{\hspace{20pt} \pi_{0}}
(V_1,\mathcal{G}_{1}) & \ar[l]_{\hspace{20pt} \pi_{1}} \cdots & \ar[l]_{\hspace{-10pt}\pi_{m-1}} 
(V_m,\mathcal{G}_{m}),
}
\end{equation*}
with $\Sing{\mathcal K}_i\subset \Sing{\mathcal G}_i$ for $i=0,\ldots,m$. We will say that 
$${\mathscr F}_V({\mathcal K})= {\mathscr F}_V({\mathcal G})$$
if ${\mathscr F}_V({\mathcal K})\subset {\mathscr F}_V({\mathcal G})$  and 
${\mathscr F}_V({\mathcal G})\subset {\mathscr F}_V({\mathcal K})$.}
\end{Definition}

\begin{Remark}
\label{treesInclusion}
Observe that if $\mathcal{G}\subset \mathcal{K}$ is an inclusion of graded rings,  then ${\mathscr F}_V(\mathcal{K})\subset {\mathscr F}_V(\mathcal{G})$. Moreover, ${\mathscr F}_V({\mathcal G})= {\mathscr F}_V({\mathcal K})$ if and only if 
${\mathscr C}_V({\mathcal G})= {\mathscr C}_V({\mathcal K})$ 
for any two ${\mathcal O}_V$-Rees algebras ${\mathcal G}$ and ${\mathcal K}$. 
\end{Remark}

There are  two different (natural) ways in which a Rees algebra ${\mathcal G}$ 
can be extended   and still remain in the same equivalence class: taking  its integral closure, say $\overline{\mathcal G}$  (see Example \ref{localintegral}) 
 or 
extending it  by the action of differential operators, say  $\DDiff({\mathcal G})$ (this 
is the smallest algebra containing $\mathcal{G}$ which is saturated (in some way) by differential operators
and it 
will  be discussed in Section \ref{extensions}). We will  show that, in fact, a combination of both operations  
leads to a complete characterization of each equivalence class.  More precisely, we will prove the following  
theorems: 

\begin{Theorem}{\bf Duality.}
\label{Canonical}
 Let $V$ be a smooth scheme over a perfect field $k$, and let ${\mathcal G}$ and 
${\mathcal K}$ be   Rees algebras. Then  
 ${\mathscr F}_V({\mathcal K})\subset  {\mathscr F}_V({\mathcal G})$ if and only if 
$\overline{\DDiff({\mathcal G})}\subset \overline{\DDiff({\mathcal K})}
$.   
\end{Theorem}

\begin{Theorem} {\bf Canonicity.}
\label{CanonicalChoice}
 Let $V$ be a smooth scheme over a perfect field $k$, and let ${\mathcal G}$ be a Rees algebra. Then the differential Rees algebra $\overline{\DDiff({\mathcal G})}$ is the canonical representative of  ${\mathscr C}_V({\mathcal G})$. 
  
\end{Theorem}

As indicated in the Introduction, these statements are related to Hironaka's 
Finite presentation 
Theorem  in \cite{Hironaka05}. However, we address here the proof of the theorems using the techniques coming from commutative algebra (already presented in the previous sections). 

Theorem \ref{Canonical}  asserts that  given two Rees algebras, ${\mathcal G}$ and ${\mathcal K}$, 
there are canonical representatives for both ${\mathscr C}_V({\mathcal G})$ and ${\mathscr C}_V({\mathcal K})$, 
 namely $\overline{\DDiff({\mathcal G})}$  and $\overline{\DDiff({\mathcal K})}
$,   in such a way that ${\mathscr F}_V({\mathcal K})\subset  {\mathscr F}_V({\mathcal G})$ 
if and only if there is an inclusion between the canonical representatives, i.e., if and only if 
$\overline{\DDiff({\mathcal G})}\subset \overline{\DDiff({\mathcal K})}
$.  

\

Theorem \ref{CanonicalChoice}  asserts that in a given equivalence class, say ${\mathscr C}_V({\mathcal G})$,  the element $\overline{\DDiff({\mathcal G})}$  is
 the 
maximum within its class with respect to inclusion   of Rees algebras. This maximum element is thus a canonical representative within the class.

\

\noindent {\bf On the strategy of the proof of Theorems \ref{Canonical} and \ref{CanonicalChoice}.} 
Notice that Theorem   \ref{CanonicalChoice} follows from Theorem  \ref{Canonical}.    Theorem \ref{Canonical} will be proved  in Section  \ref{proof_canonical} using  algebraic methods,  like  the ones developed in   Section  \ref{technical_integral_closure} and the normalized blow up of a Rees algebra introduced in  \ref{normalblowup}. A key point that follows from   Theorem \ref{CanonicalChoice} is the fact that the extension of a Rees algebra to the Differential Rees algebra that it generates is canonical (see Remark \ref{localgenDif}). This issue, and more generally, the action of differential operators on Rees algebras will be treated in Section \ref{extensions}. In 
Section  \ref{technical_integral_closure} we will study questions regarding to the inclusion $\overline{\DDiff({\mathcal G})}\subset \overline{\DDiff({\mathcal K})}
$ 
 that appears in the statement of Theorem \ref{Canonical}. 
In particular we establish a suitable 
criterion  to characterize the    integral  closure of a  Rees algebra   (here we make use of ideas that appear in   \cite{kaw}).  Section  \ref{technical_closed}  is devoted to the study of the inclusion  ${\mathscr F}_V({\mathcal K})\subset  {\mathscr F}_V({\mathcal G})$ of Theorem \ref{Canonical}. 
As indicated above, the 
 proof of Theorem \ref{Canonical}  is finally  given in  Section 
\ref{proof_canonical}.

\begin{Remark}
There may be other (natural) notions of ${\mathcal G}$-{\em local sequences} 
over smooth schemes, that in principle could lead to different 
equivalence relations on Rees algebras. This will be discussed in  Section 
 \ref{other_weak_equivalence}.
\end{Remark}

\section{Differential Rees algebras and Giraud's Lemma}
\label{extensions}
We have already proved that a Rees algebra and its  integral closure are weakly equivalent (see Example \ref{localintegral}). In  this section we will see that a Rees algebra and the {\em differential Rees algebra} (see Definition \ref{DiffAlg} below)  expanded by it are also 
weakly equivalent.  This is essentially a result of Giraud (see Lemma \ref{difintegral}). 
See also \cite{WLL}  and \cite{kaw} for other results in this line.  

\

Let  $V$  be a smooth scheme  over a  perfect field $k$. For any 
non-negative integer $r$,  denote by $\diff^r_{V|k}$  the (locally free) sheaf of $k$-linear differential operators
of order at most $r$.

\begin{Definition}\label{DiffAlg} A   Rees algebra ${\mathcal G}=\bigoplus_nI_nW^n$
is said to be a {\em differential Rees algebra}, if the
following condition holds:
\begin{enumerate}
\item[] There is an affine open covering of $V$, $\{U_i\}$, such
that for any $D\in \diff^r_{V|k}(U_i)$ and any $h\in I_n(U_i)$ we have
that $D(h)\in I_{n-r}(U_i)$ provided that $n\geq r$.
\end{enumerate} 
In particular, $I_{n+1}\subset I_n$, since $\diff^0_{V|k}(U_i)\subset \diff^1_{V|k}(U_i)$.
\end{Definition}

\begin{Remark}
\label{localgenDif}
{\rm Given any  Rees algebra  ${\mathcal G}$ 
 on a smooth scheme $V$, there is a
natural way to construct  the smallest differential Rees algebra containing it:  the \emph{differential 
Rees  algebra generated 
by ${\mathcal G}$},   $\DDiff(\mathcal{G})$ (see \cite[Theorem 3.4]{integraldifferential}). 
 More precisely,  if ${\mathcal G}$ is locally generated by $\{f_{1}W^{n_1},\ldots,f_{s}W^{n_s}\}$ 
 on an affine open set $U$, 
 then it can be shown that
 $\DDiff(\mathcal{G}(U))$ is generated by the elements 
\begin{equation}
\label{generatorsDiff}
\{D(f_{i})W^{n_i-r}: D\in \diff^r_{V|k}(U),  0\leq r< n_i, i=1,\ldots,s\}.
\end{equation} 
Note that $\diff^r_{V|k}(U)\subset \diff^{\ell}_{V|k}(U)$ if $r \leq \ell$. Thus, if $D\in 
\diff^r_{V|k}(U)\subset \diff^{\ell}_{V|k}(U)$ is a differential operator, then $D(f_i)W^{n_i-r}$ is in (\ref{generatorsDiff}) and  also $D(f_i)W^{n_i-\ell}$ is in (\ref{generatorsDiff}), as long as $r\leq \ell < n_i$. Moreover, it suffices to take $D$ as part of a   finite system of generators   of $\diff^r_{V|k}(U)$ with $r$ being  strictly smaller than $n_i$   which in particular implies that the differential algebra generated by ${\mathcal G}$ is a (finitely generated) Rees algebra 
(cf. \cite[Proof of Theorem 3.4]{integraldifferential}). 
}
\end{Remark}

\begin{Paragraph}\label{Local_Gene_Diff_Op}
{\bf Local generators for the sheaf of differential operators.} Along these lines we give a local 
description of the generators of the locally free sheaf $\diff^r_{V|k}$. Let $x\in V$ be a closed point, and 
let  $\{z_1,z_2,\ldots,z_d\}\subset {\mathcal O}_{V,x}$ be a regular system of parameters.  On $\widehat{{\mathcal O}_{V,x}}\simeq k^{\prime}[[z_1,\ldots,z_d]] $, where $k^{\prime}$ is the residue field at $x$,  consider the Taylor expansion  
\begin{equation}
\label{Formula_Taylor}
\begin{array}{r@{\hspace{0pt}}rcl} 
\operatorname{Tay} :& k^{\prime}[[z_1,\ldots,z_d]]  & \longrightarrow & k^{\prime}[[z_1,\ldots,z_d,T_1,\ldots,T_d]]\\
 & f(z_1,z_2,\ldots,z_d) & \longmapsto & f(z_1+T_1,z_2+T_2,\ldots,z_d+T_d)=\sum_{\alpha\in {\mathbb N}^d}\Delta^{\alpha} (f)T^{\alpha}
\end{array}
\end{equation}
as in \cite[Definition 1.2]{hpositive}
and in \cite[Theorem 3.4]{integraldifferential}. Then  for each $\alpha=(\alpha_1,\alpha_2,\ldots,\alpha_d)\in {\mathbb N}^d$, 
\begin{equation}\label{DeltaComp1}
\begin{array}{rrcl}
\Delta^{(\alpha_1,\alpha_2,\ldots,\alpha_d)}: & \widehat{{\mathcal O}_{V,x}} & \longrightarrow  &  \widehat{{\mathcal O}_{V,x}} \\
 \ & f & \longmapsto &
\Delta^{(\alpha_1,\alpha_2,\ldots,\alpha_d)}(f) 
\end{array}
\end{equation}
is a differential operator of order $|\alpha|=\alpha_1+\alpha_2+\ldots+\alpha_d$, which defines by restriction 
a differential operator $D^{\alpha}: {\mathcal O}_V(U)\to {\mathcal O}_V(U)$ in some neighborhood $U$ of $x$, since $V$ is smooth over the perfect field $k$.  Moreover, 
the sheaf of differential operators  up to order $r$, say $\diff_{V|k}^r$, is locally generated by the $D^{\alpha}$ with 
$|\alpha|\leq r$
(at $U$). 
\end{Paragraph}

\begin{Paragraph}\label{difsingular}{\bf Differential Rees algebras, order,  singular locus and zero set.}
Let $J\subset {\mathcal O}_V$ be a non-zero sheaf of ideals. Then, for a non-negative 
integer $r$, define 
$$\diff^{r}_{V|k}(J):=\langle D(f):D\in \diff^{r}_{V|k} \operatorname{ and } f\in J\rangle.$$
Let $x\in V$, and let $b$ be a non-negative integer. Since $V$ is smooth over a perfect field $k$, using  Taylor expansions 
as in \ref{Local_Gene_Diff_Op},   note that 
$$\nu_x( J)\geq b \Leftrightarrow x\in V(\diff^{b-1}_{V|k}(J)).$$ 
Therefore,  if
  ${\mathcal
G}=\bigoplus_nI_nW^n$, then, 
$$\Sing{\mathcal G}=\bigcap_{n\geq 1}V(\diff^{n-1}_{V|k}(I_n)),$$
(see \cite[Definition 4.2 and Proposition 4.4]{integraldifferential}). In particular, 
$$\Sing{\mathcal
G}=\Sing\DDiff(\mathcal{G});$$ and moreover, if $x\in
\Sing{\mathcal G}=\Sing (\DDiff ({\mathcal G}))$ then
$$
 \ord_{\mathcal
G}(x)=\ord_{\DDiff(\mathcal{G})}(x) 
$$ (cf.
\cite[Proposition 6.4  
(3)]{EV}). Furthermore, if  ${\mathcal G}$ is a
differential Rees algebra, then $\Sing{\mathcal G}=V(I_n)$ for any
positive integer $n$ (see \cite[Proposition
4.4
(5)]{integraldifferential}), and  $\Sing\mathcal{G}=V(\mathcal{G})$ (see \ref{zeros}).
\end{Paragraph}

\begin{Paragraph} {\bf Local sequences and differential extensions.} \label{localExtDiff} The pull-back of a Differential Rees algebra by a smooth morphism is a 
Differential Rees algebra again. So, to completely understand the behavior of 
Differential Rees algebras when considering local sequences on $V$ we need to 
study the case of permissible transformations. A main consequence of Giraud's Lemma \ref{difintegral}   (see part (ii)) is that the   transforms of $\mathcal{G}$ and $\DDiff(\mathcal{G})$ have the same invariants studied in \ref{difsingular}. 

\begin{Lemma} \label{difintegral}  {\bf Giraud's Lemma 
 \cite[Theorem 4.1]{EV}}. Let ${\mathcal
G} \subset {\mathcal K}\subset {\mathcal R}$ be an inclusion
of Rees algebras, such that 
$\mathcal{R}=\DDiff(\mathcal{G})$,
and let $V\leftarrow
V_1$ be a permissible monoidal transformation with center
$Y\subset \Sing {\mathcal R}$($=\Sing {\mathcal G}=
\Sing {\mathcal K}$). Then:

(i) There is an inclusion of  transforms
$${\mathcal G}_1 \subset {\mathcal K}_1\subset
{\mathcal R}_1.$$

(ii) Even if ${\mathcal R}_1$ is not a differential Rees algebra
 over $V_1$,  the three algebras ${\mathcal G}_1 \subset {\mathcal K}_1\subset
{\mathcal R}_1$ span the same differential Rees algebra, and therefore
$$\Sing{\mathcal G}_1=\Sing{\mathcal K}_1=\Sing
{\mathcal R}_1.$$
\end{Lemma}
\end{Paragraph}

\noindent {\bf {Summarizing: who is in ${\mathscr C}_V({\mathcal G})$?}}
So far we conclude that    any   Rees algebra ${\mathcal G}$ is weakly equivalent 
to both $\overline{\mathcal G}$ and  $\DDiff({\mathcal
G})$.  Theorem \ref{CanonicalChoice} asserts that combining both operators 
we  get a characterization  all  the Rees algebras in  ${\mathscr C}_V({\mathcal G})$, since $\overline{\DDiff({\mathcal K})}=
\overline{\DDiff({\mathcal G})},$ for all ${\mathcal K}\in {\mathscr C}_V({\mathcal G})$.

\section{Testing integral closure on Rees algebras}\label{technical_integral_closure}
The purpose of this section is to develop suitable tools to check   
whether  a given Rees algebra is contained in the integral closure of another.   
 In Section \ref{ReesAlgebras} we showed that  this question can be formulated in 
terms of   integral closure of ideals   (see \ref{Veronese} 
and  Lemmas \ref{N_adecuado} and 
 \ref{almostIntegral}).  Our main goal here is to  prove  Proposition \ref{CasoHyper}, which is a particular case of Theorem \ref{Canonical}. In fact the proof of Theorem \ref{Canonical} will be   reduced to this special  case.

\begin{Remark}\label{IntegCriterBlowup}  Recall  that if  $I$ and $J$  are  ideals in a normal domain $B$, and  if 
$$\xymatrix@R=2pt@C=20pt{
\Spec B & \ar[l]_{\hspace{15pt}\Theta}  \mathscr{B}
}
$$
is  the normalized blow up of $\Spec B$ at $I$, 
then $J {\mathcal O}_{\mathscr{B}}\subset I {\mathcal O}_{\mathscr{B}}$
 if and only if $J\subset \overline{I}$ in $B$. Moreover, since $\mathscr{B}$ is normal, and $I{\mathcal O}_{\mathscr{B}}$ is locally invertible, it follows that $J {\mathcal O}_{\mathscr{B}}\subset I {\mathcal O}_{\mathscr{B}}$ if and only if 
$J {\mathcal O}_{\mathscr{B},\mathbb{H}}\subset I {\mathcal O}_{\mathscr{B},{\mathbb{H}}}$ for all irreducible hypersurfaces   $\mathbb{H}$ of $\mathscr{B}$.
\end{Remark}

  Remark \ref{IntegCriterBlowup}
motivates our interest in grasping the integral closure 
of a particular class of Rees algebras, namely those  generated by
 a locally principal ideal  in some weight.

\begin{Lemma} \label{integralhyp} Let $\mathbb{H}$ be an irreducible and reduced  hypersurface on a smooth scheme  $V$, let ${\mathcal H}={\mathcal O}_V[ I(\mathbb{H})^NW^n]$, and let ${\mathcal G}=$ $\mathcal{O}_V[I(\mathbb{H})^QKW^q]$ for some  ideal sheaf $K \not\subset I(\mathbb{H})$. Then ${\mathcal G}\subset \overline{\mathcal H}$ if and only if $\frac{Q}{q}\geq \frac{N}{n}$. 
\end{Lemma} 

\noindent  {\em Proof:}  We first show that if $\mathcal{G}\subset \overline{\mathcal{H}}$ then $\tfrac{Q}{q}\geq \tfrac{N}{n}$. By Example \ref{InClosure_ring_almost},  looking at the degree $nq$ piece in both algebras we find that $I(\mathbb{H})^{Qn}K^n \subset \overline{I(\mathbb{H})^{Nq}}$.  And by
Remark \ref{IntegCriterBlowup},    $I(\mathbb{H})^{Qn}\mathcal{O}_{V,\mathbb{H}}=I(\mathbb{H})^{Qn}K^n\mathcal{O}_{V,\mathbb{H}} \subset I(\mathbb{H})^{Nq}\mathcal{O}_{V,\mathbb{H}}$,  so $Qn\geq Nq$.
The converse follows from Lemma \ref{almostIntegral} since $Qn\geq Nq$ implies that $I(\mathbb{H})^{Qn}K^n\subset I(\mathbb{H})^{Nq}= \overline{I(\mathbb{H})^{Nq}}$.\qed

\begin{Proposition} \label{CasoHyper} Let $\mathbb{H}$ be an irreducible and reduced  hypersurface on a smooth scheme  $V$, let ${\mathcal H}={\mathcal O}_V[I(\mathbb{H})^NW^n]$ and let 
${\mathcal G}={\mathcal O}_{V}[JW^b]$ for some sheaf of ideals
 $J\subset {\mathcal O}_V$. Assume $N\geq n$. Then 
${\mathscr F}_V({\mathcal H})\subset {\mathscr F}_V({\mathcal G})$ if and only if $\mathcal G\subset 
\overline{\mathcal H}$.
\end{Proposition}

\noindent  {\em Proof:}   Thanks to the finiteness of the Veronese action (see \ref{Veronese}) and the fact that two Rees algebras that have the same integral closure are weakly equivalent (see  
Example 
\ref{localintegral}) we can assume ${\mathcal H}={\mathcal O}_{V}[I(\mathbb{H})^{Nb}W^{nb}]$   and that  ${\mathcal G} ={\mathcal O}_{V}[J^nW^{nb}]$. Writing $J^n=I(\mathbb{H})^QJ_1$ with $J_1\not\subset I(\mathbb{H})$, and by Lemma \ref{integralhyp}, we note that $\mathcal{G}\subset \overline{\mathcal{H}}$ if and only if $\frac{Q}{nb} \geq \frac{N}{n}$ i.e., if and only if  $Q \geq Nb$.

Suppose first  that  $Q \geq Nb$. This means  that $\mathcal{G}\subset \mathcal{H}$ as graded rings, and by 
Remark 
\ref{treesInclusion},  $\mathscr{F}_V(\mathcal{H})\subset \mathscr{F}_V(\mathcal{G})$.

To prove the converse,  we will assume that   ${\mathscr F}_V({\mathcal H})\subset {\mathscr F}_V({\mathcal G})$ to conclude that $Q\geq Nb$. 
If $N=n$, the conclusion  is clear  since $\mathbb{H}=\Sing\mathcal{H} \subset \Sing\mathcal{G}= \{x \in V, \nu_x(I(\mathbb{H})^QJ_1)\geq nb\}$, so $Q\geq nb=Nb$. 

If $N>n$, we will 
use   Hinonaka's trick, and  construct 
a suitable $\mathcal{H}$-local sequence (which by hypothesis will also be a $\mathcal{G}$-local sequence).   First notice  that if $\mathbb{H}\subset V$ is not smooth we can replace $V$ with an open set $U$ (by  removing a closed set of codimension at least two) so that $\mathbb{H}|_U$ is smooth. Since the restriction $U\rightarrow V$ is a smooth morphism, one has that  $\mathscr{F}_U(\mathcal{H})\subset \mathscr{F}_U(\mathcal{G})$. Now consider    the product of $U$ by an affine line, say $V^{\prime}_0=U\times {\mathbb A}_k^1$, and the corresponding pull-backs 
${\mathcal H}^{\prime}_0$ and ${\mathcal G}^{\prime}_0$ of ${\mathcal H}$ and ${\mathcal G}$, respectively, to $V^{\prime}_0$. 
Then $\mathbb{H}^{\prime}_0=\mathbb{H}|_U\times {\mathbb A}^1_k=\Sing{\mathcal H}^{\prime}_0\subset \Sing{\mathcal G}^{\prime}_0$. 
Let $x\in {\mathbb A}^1_k$ be a closed point and consider the permissible monoidal transformation with  center 
 $\mathbb{H}|_U\times \{x\}$,  say $V^{\prime}_0\leftarrow V^{\prime}_1$. Let $H_1$ be the exceptional divisor, and let $\mathbb{H}_1^{\prime}$ denote the 
 strict transform of  $\mathbb{H}^{\prime}_0$ in $V^{\prime}_1$. Now   the transforms ${\mathcal H}^{\prime}_1$ and 
${\mathcal G}^{\prime}_1$ of ${\mathcal H}^{\prime}_0$ and ${\mathcal G}^{\prime}_0$ respectively, factor as follows:
$${\mathcal H}^{\prime}_1={\mathcal O}_{V^{\prime}_1}[I(\mathbb{H}_1^{\prime})^{Nb}I(H_1)^{(N-n)b}W^{nb}]; \ \ 
{\mathcal G}^{\prime}_1={\mathcal O}_{V^{\prime}_1}[I(\mathbb{H}_1^{\prime})^QI(H_1)^{Q-nb}J_1{\mathcal O}_{V_1^{\prime}}W^{nb}].$$ 

Blow up at the permissible center 
$\mathbb{H}_1^{\prime}\cap H_1$, and continue blowing up at   the successive intersections of the strict transforms of $\mathbb{H}_1^{\prime}$ and the new exceptional divisors
so as to obtain the sequence
\begin{equation}\label{sucesion_sin_cruces_normales}\begin{xy} 
(0,6)*+{V}; (12,6)*+{U}; (35,6)*+{U\times \mathbb{A}^1_k=V'_0}; (70,6)*+{V'_1}; (90,6)*+{V'_2};(108,6)*+{\ldots}; (128,6)*+{V'_{m-1}};(148,6)*+{V'_m};%
(35,0)*+{\mathbb{H}|_{U}\times \{x\}}; 
(70,0)*+{\mathbb{H}_1^{\prime}\cap H_{1}}; %
(128,0)*+{\mathbb{H}_{m
-1
}^{\prime}\cap H_{m-1}.};
{\ar (9,6);(3,6) } ; 
{\ar (21,6);(15,6) };
{\ar (67,6);(49,6)\POS?*^+!D{\scriptstyle\rho_0}}; 
{\ar (86,6);(73,6)\POS?*^+!D{\scriptstyle\rho_1}}; 
{\ar (104,6);(94,6)}; 
{\ar (122,6);(112,6)};
{\ar (144,6);(134,6)\POS?*^+!D{\scriptstyle\rho_{m-1}} }; 
\end{xy}\end{equation}

The corresponding transforms of ${\mathcal H}^{\prime}_1$ and  ${\mathcal G}^{\prime}_1$ can be expressed as:
$${\mathcal H}^{\prime}_m={\mathcal O}_{V^{\prime}_m}[I(\mathbb{H}^{\prime}_m)^{Nb}I(H_1)^{(N-n)b}\cdots I(H_m)^{m(N-n)b}W^{nb}]; $$
$${\mathcal G}^{\prime}_m={\mathcal O}_{V^{\prime}_m}[I(\mathbb{H}^{\prime}_m)^QI(H_1)^{Q-nb}\cdots I(H_m)^{m(Q-nb)}J_1{\mathcal O}_{V^{\prime}_m}W^{nb}].$$ 

To conclude observe that for $m$ sufficiently large, $H_m$ is a permissible center for ${\mathcal H}^{\prime}_m$
which can be blown up
at most $\left[\frac{m(N-n)b}{nb}\right]$ times. Since by hypothesis ${\mathscr F}_{V^{\prime}_m}({\mathcal H}_m^{\prime})\subset {\mathscr F}_{V^{\prime}_m}({\mathcal G}_m^{\prime})$, 
these blow ups are permissible for $\mathcal{G}_m^{\prime}$ too. Thus $\left[\frac{m(Q-nb)}{nb}\right]\geq \left[\frac{m(N-n)b}{nb}\right]$, 
and since $m \gg 0$ necessarily we have
$\frac{Q-nb}{nb}\geq \frac{Nb-nb}{nb}$, and hence $Q\geq Nb$.  \qed

\section{On restrictions of   trees    to smooth closed subschemes} 
\label{technical_closed}
 Recall that  a Rees algebra  ${\mathcal G}$  over $V$ determines a tree of closed sets, say  ${\mathscr F}_V({\mathcal G})$, by 
considering ${\mathcal G}$-local sequences on $V$. Indeed,  the class ${\mathscr C}_V({\mathcal G})$   
is completely determined by  ${\mathscr F}_V({\mathcal G})$ (see Remark \ref{treesInclusion}). Roughly speaking, in 
this section we study the restriction of the trees ${\mathscr F}_V({\mathcal G})$ to smooth closed subschemes  $X\subset V$.  The central result   is Proposition \ref{restriction} which  plays a key  role 
in the proof of Theorem \ref{Canonical}.   In Proposition \ref{restriction} we treat the notion of restriction of differential operators, which is one of the most significant properties of Giraud's use of  higher differential operators. This also appears as the 
{\em restriction property} in \cite{Hironaka05}.

\begin{Definition}
\label{defInterseccion}
{\rm Let $ {\mathcal G}$,  ${\mathcal K}$ and ${\mathcal R}$   be Rees algebras on $V$. 
We will say that 
$${\mathscr F}_V({\mathcal R})={\mathscr F}_V({\mathcal G})
\cap {\mathscr F}_V({\mathcal K})$$ if: } 
\begin{enumerate}
\item[(i)] $\Sing {\mathcal R} = \Sing {\mathcal G}\cap \Sing {\mathcal K}$; 
\item[(ii)]  Any   ${\mathcal G}$-${\mathcal K}$-local sequence  over $V$ induces an  $ {\mathcal R}$-local sequence over $V$, and any $ {\mathcal R}$-local sequence over $V$ induces a $ {\mathcal G}$-${\mathcal K}$-local sequence  over $V$. 
\item[(iii)]  Given a  $ {\mathcal G}$-${\mathcal K}$-${\mathcal R}$-local sequence  over $V$,
$$
\xymatrix@R=2pt@C=20pt{
(V,\mathcal{G}, \mathcal{K}, \mathcal{R})=(V_0,\mathcal{G}_{0}, \mathcal{K}_{0},\mathcal{R}_{0}) & \ar[l] 
(V_1,\mathcal{G}_{1},\mathcal{K}_{1}, \mathcal{R}_{1}) & \ar[l] \cdots  & \ar[l] (V_m,\mathcal{G}_{m}, \mathcal{K}_{m}, \mathcal{R}_{m}),
}$$
there is an equality of closed sets  $$\Sing{\mathcal R}_i=\Sing{\mathcal G}_i\cap \Sing{\mathcal K}_i \operatorname{ \ for }  i=0,\ldots,m.$$
\end{enumerate}
\end{Definition}

\begin{Definition}\label{amalgamagenerada}
Let ${\mathcal G}$ and ${\mathcal K}$ be two Rees algebras on $V$. We 
will denote by ${\mathcal G}\odot {\mathcal K}$ the smallest Rees algebra containing both ${\mathcal G}$ and ${\mathcal K}$. Suppose that on an open 
 affine  set $U\subset V$,  
${\mathcal G}$ is locally generated by 
$f_1W^{n_1},\ldots,f_sW^{n_s}$ and ${\mathcal K}$  is locally  generated by 
$g_1W^{m_1},\ldots, g_rW^{m_r}$.     
 Then $({\mathcal G}\odot {\mathcal K})(U)$ is generated 
by $f_1W^{n_1},\ldots,f_sW^{n_s}, g_1W^{m_1},\ldots, g_rW^{m_r}$. 
\end{Definition}
\begin{Lemma}
\label{amalgama}
Let $ {\mathcal G}$,  ${\mathcal K}$ be two ${\mathcal O}_V$-Rees algebras. Then 
$${\mathscr F}_V({\mathcal G}\odot {\mathcal K})={\mathscr F}_V({\mathcal G})
\cap {\mathscr F}_V({\mathcal K}).$$
\end{Lemma} 

\noindent {\em Proof:} By  the local description of Definition \ref{amalgamagenerada}, the lemma is a straight consequence of  \ref{singularlocus} and
Proposition   
\ref{localt}. \qed

\begin{Remark}
 By Lemma  \ref{amalgama}, the tree ${\mathscr F}_V({\mathcal G})\cap {\mathscr F}_V({\mathcal K})$  is in fact a tree defined by a Rees algebra. In addition, with the notation of Definition \ref{definition_inclusion}, we notice that the tree $\mathscr{F}_V(\mathcal{R})$ is the biggest  defined by a Rees algebra satisfying both conditions,  ${\mathscr F}_V({\mathcal R})\subset {\mathscr F}_V({\mathcal G})$, and ${\mathscr F}_V({\mathcal R})\subset {\mathscr F}_V({\mathcal K})$. In other words, if there is another Rees algebra, say ${\mathcal S}$, such that ${\mathscr F}_V({\mathcal S})\subset {\mathscr F}_V({\mathcal G})\cap {\mathscr F}_V({\mathcal K})$, then ${\mathscr F}_V({\mathcal S})\subset {\mathscr F}_V({\mathcal R})$. 
\end{Remark}

\begin{Paragraph} \label{discussion} 
{\bf Local sequences over smooth closed subschemes.} 
Intersections of trees of closed subsets over $V$, as   in Definition \ref{defInterseccion}, are of 
 special interest in the following particular case. Assume that ${\mathcal G}$ is 
any ${\mathcal O}_V$-Rees algebra,  set $X\subset V$ be a smooth closed subscheme, and 
let ${\mathcal X}={\mathcal O}_V[I(X)W]$. {Then} any ${\mathcal G}$-${\mathcal X}$-local sequence  over $V$, 
\begin{equation}
\label{localintersection}
\xymatrix@R=2pt@C=20pt{
(V,\mathcal{G}, \mathcal{X})=(V_0,\mathcal{G}_{0}, \mathcal{X}_{0}) & \ar[l]_{\hspace{30pt}\pi_0}
(V_1,\mathcal{G}_{1},\mathcal{X}_{1}) & \ar[l]_{\hspace{20pt}\pi_1}\cdots & \ar[l]_{\hspace{-30pt}\pi_{m-1}}
(V_m,\mathcal{G}_{m}, \mathcal{X}_{m}),
}
\end{equation}
induces, at the same time, a local sequence over  $X$ (where   $X_{i+1}=\Sing{\mathcal X}_{i+1}$ is either  the strict transform of $X_i$ in  
$V_{i+1}$,  
if $\pi_i: V_i\leftarrow V_{i+1}$ is a permissible monoidal transformation, or the preimage  of $X_i$ by a smooth 
morphism otherwise), 
\begin{equation*}
(X,\Sing{\mathcal G} 
\cap X)=(X_0,\Sing \mathcal{G}_{0}\cap 
X_0) \stackrel{\pi_0|_{X_0}}{\longleftarrow}
(X_1,\Sing \mathcal{G}_{1}\cap 
 {X}_{1})\stackrel{\pi_1|_{X_1}}{\longleftarrow}\cdots \stackrel{\pi_{m-1}|_{X_{m-1}}}{\longleftarrow}
(X_m,\Sing \mathcal{G}_{m}\cap 
{X}_{m}),
\end{equation*}
$$\xymatrix@R=2pt@C=40pt{
X=X_0 & \ar[l]_{\hspace{20pt}\pi_0|_{X_0}}
X_1 & \ar[l]_{\pi_1|_{X_1}} \cdots & \ar[l]_{\pi_{m-1}|_{X_{m-1}}}
X_m,
}$$
and  each $\Sing \mathcal{G}_{i}\cap 
{X}_{i}$ can be naturally interpreted 
as a closed subset of $X_i$ for $i=0,1,\ldots,m$. Here the tree ${\mathscr F}_V({\mathcal G})$ 
of  
closed sets over $V$ induces a tree of closed sets over $X$,  and 
our purpose is to  
find out  if this is the tree of some  $\mathcal{O}_X$-Rees algebra
(see Question \ref{comparar_arriba_abajo}). To  face this question, we will have  to compare trees on $X$ and trees on $V$, so we will start by showing that  
   local sequences over $X$ can be lifted, at least locally,  to local sequences over $V$. 
   
Observe  that a local sequence   
   over $X$ includes, for example, a smooth morphism say $X\stackrel{\varphi}{\longleftarrow}  U$.  Since $X$ is closed in $V$ there is, a priori, no natural way to define a smooth morphism, say 
  $V\stackrel{{\varphi}^{\prime}}{\longleftarrow}  U'$, so that $X\stackrel{\varphi}{\longleftarrow}  U$ can be identified with the restriction
  of the former over the pull-back of $X (\subset V)$.  The existence of local retractions (shown in \ref{Existence_Local_Retra}) will enable us to overcome this difficulty. We shall first show that   for  $X\subset V$ as above, one can 
 define local retractions $X\stackrel{p}{\longleftarrow}  V$. Here we may have to replace $V$ by an open cover of $V$ in \'etale  topology.
 Once this point is settled,  we show that given $X\subset V$  together with a retraction $X\stackrel{p}{\longleftarrow}  V$,
  any local sequence over $X$ can be lifted to a local sequence over $V$. This will suffice to address the question formulated above (see \ref{Levantar_local}).

\begin{Paragraph}\label{Existence_Local_Retra} {\bf On the existence of local retractions.}   Let 
$x\in X\subset V$ be  a closed point. 
A regular system of parameters $\{z_1, \dots , z_d\}$ in $\mathcal{O}_{V,x}$ defines an inclusion of a polynomial ring in $d$ variables, say 
$k'[Z_1, \dots , Z_d]\subset \mathcal{O}_{V,x}$, where $k'$ is the residue field at $x$. 
This in turn says that $(V,x)$ is an \'etale neighborhood of $(\mathbb{A}^{d},\mathbb{O})$. Assume  the closed smooth 
subscheme $X$ has dimension $e\leq d$. Then the regular system of parameters above can be chosen so that   $z_{e+1}=0,\dots, z_d=0$ defines $X$ locally in a neighborhood of $x$, and $(X,x)$ is an \'etale 
neighborhood of $(\mathbb{A}^{e},\mathbb{O})$.   The projection of  $\mathbb{A}^{d}\to  \mathbb{A}^{e}$ over the first $e$-coordinates leads to a diagram 
$$
\xymatrix{
(\mathbb{A}^{d},\mathbb{O}) \ar[d] & (V, x) \ar[l]\\
(\mathbb{A}^{e},\mathbb{O})  & (X, x), \ar[l]
}
$$ where the horizontal arrows are \'etale, and the vertical one is smooth. 
Taking fiber products we get   
$$
\xymatrix@C=10pt{
V \ar[d] & **[r]V'=V \times X \ar[l] \ar[d]^{\rho}\\
\mathbb{A}^{e} & X. \ar[l]
}
$$
with vertical smooth arrows and horizontal \'etale arrows in a neighborhood of $x$.
 The closed immersion $X\subset V$ defines a diagonal map, $X \rightarrow V$. This, in turn, indicates that $\rho$ admits a  section, say $s:X\to V^{\prime}$, which 
defines  a local retraction in a neighborhood of $x$, $\rho: V^{\prime} \to X$. 
\end{Paragraph}

\begin{Paragraph}\label{Levantar_local} {\bf Lifting local sequences.} 

\noindent {-\em Lifting smooth morphisms.} If $X_1\rightarrow X$ is a smooth morphism, then  
$$
\xymatrix@C=10pt{
V' \ar[d]^{\varphi}& **[r]V_1=V' \times X_1 \ar[d]^{\varphi_1} \ar[l] \\
X & X_1. \ar[l]}
$$
induces  a smooth morphism over $V^{\prime}$,    and a commutative diagram. The closed immersion $X\subset V^{\prime}$ guarantees the existence of a diagonal map from $X_1$ to $V^{\prime}$. So, \ as before,  there is a section $s_1: X_1\to V_1$ which together with $\varphi_1: V_1\to X_1$ defines a retraction. 

\noindent {-\em Lifting blow ups.} Observe, in addition, that if $Y\subset X$ is a smooth closed subscheme, then the blow up with center $Y$ induces a 
commutative diagram of blow ups and closed immersions, 
$$
\xymatrix{
V & V_1 \ar[l] \\
X \ar[u] & X_1 \ar[u] \ar[l]
}
$$
where $X_1$ can be identified with the strict transform of $X$. Moreover, the diagram also indicates that a local retraction $\varphi: V\to X$ can be naturally lifted to a local retraction $\varphi_1: V_1\to X_1$ at least locally in a neighborhood of a point. 
\end{Paragraph}

\begin{quote}
{\bf Conclusion 1.}  Thus, using the retraction $\varphi: V^{\prime}\to X$, any local sequence over $X$ can be lifted to a local sequence over 
$V$ (in the sense of Definition \ref{local_sequenceA}). 
\end{quote}
\begin{quote}
{\bf Conclusion 2.}  If a retraction $\varphi: V \to X$ is given, then 
given   an ${\mathcal O}_X$-Rees  algebra  ${\mathcal S}$, 
the tree of closed subsets 
${\mathscr F}_X({\mathcal S})$ can be identified with a tree of closed subsets over $V$.
\end{quote}
\end{Paragraph}

\begin{Question} 
\label{comparar_arriba_abajo} Let  ${\mathcal G}$ be  
any ${\mathcal O}_V$-Rees algebra,  let $X\subset V$ be a smooth closed subscheme, and 
let ${\mathcal X}={\mathcal O}_V[I(X)W]$. Is there   an ${\mathcal O}_X$-Rees algebra ${\mathcal S}$ 
such that 
$${\mathscr F}_X({\mathcal S})={\mathscr F}_V({\mathcal G})\cap {\mathscr F}_V
({\mathcal X})?$$
Here the inclusion  ``$\subset$''     has to be interpreted viewing $\mathscr{F}_X(\mathcal{S})$ as a tree of closed sets over $V$, and the equality of trees should hold  for {\em any} local 
  retraction $\varphi: V\to X$ (since, a priori,   no local retraction $V\to X$ is given, and 
infinitely many can be constructed). 

The first obvious  guess, setting ${\mathcal S}$ as ${\mathcal G}|_X$, 
will not work since in  general the inclusion 
$$\Sing{\mathcal G}\cap X \subset \Sing({\mathcal G}|_X)$$
is strict. 
However, a more careful choice of representative within the class ${\mathscr C}_V({\mathcal G})$ will do the job,  as  it is shown in the next proposition. 

\end{Question}

\begin{Proposition} \label{restriction} Let $X\subset V$ be a smooth closed subscheme. Set  ${\mathcal X}={\mathcal O}_V[I(X)W]$ and 
let  ${\mathcal G}$ be an arbitrary  Rees algebra. Then 
 $${\mathscr F}_X((\DDiff(\mathcal{G}))|_X)
={\mathscr F}_V({\mathcal G})\cap {\mathscr F}_V({\mathcal X}).$$

Moreover, 
$(\DDiff(\mathcal{G}))|_X$
is an ${\mathcal O}_X$-differential Rees algebra. 

\end{Proposition}

\noindent{\em Proof:} Lemma \ref{amalgama} asserts that 
 $${\mathscr F}_V\left({\mathcal G}\odot {\mathcal X} \right)=
{\mathscr F}_V({\mathcal G})\cap {\mathscr F}_V({\mathcal X}).$$
 It suffices to prove the proposition in the particular case in which ${\mathcal G}={\mathcal O}_V[fW^b]$ and 
$X$ is a smooth hypersurface, since the general case would then follow 
from this using an inductive argument. 

Let $x\in (\Sing{\mathcal G})\cap X$ be a closed point, and fix a regular system of parameters $\{z_1,\ldots,z_d\}\subset {\mathcal O}_{V,x}$ 
with $z_1$ being a local equation for $X$. Consider the expansion of $f$ in the completion $\widehat{{\mathcal O}_{V,x}}$, 
\begin{equation}
\label{desarrollo}
f=h_0(z_2,\ldots,z_d)+z_1h_1(z_2,\ldots,z_d)+z_1^2h_2(z_2,\ldots,z_d)+\ldots=\sum_{i=0}^\infty z_1^{i}h_i(z_2,\ldots, z_d).
\end{equation}
Notice that there is a  formal retraction defined in algebraic 
terms   by the  natural inclusion  
\begin{equation}
\label{retractionComple}
\rho: \widehat{{\mathcal O}_{X,x}}\simeq k^{\prime}[[z_2,\ldots,z_d]]\to \widehat{{\mathcal O}_{V,x}}\simeq k^{\prime}[[z_1,z_2,\ldots,z_d]],
\end{equation} where $k^\prime$ is the residue field at $x$. This inclusion 
allows us to view the  $h_i$ as elements in  both $\widehat{{\mathcal O}_{V,x}}$ and $\widehat{{\mathcal O}_{X,x}}$ for   $i=0,1,\ldots$. 

Now observe  that $x\in \Sing({\mathcal G}\odot  {\mathcal X})$ if and only if 
\begin{equation}
\label{CondicionHi}
\nu_x(h_0)\geq b,  \ \  \nu_x(h_1)\geq b-1,\ldots,
   \nu_x(h_{b-1})\geq 1,  \ \ \nu_x(z_1)\geq 1;
  \end{equation}
where $\nu_x$ denotes the order function in the regular local ring $\widehat{{\mathcal O}_{V,x}}$.  This,  in turn,  is equivalent to asking that $\nu_{x}(h_i)\geq b-i$ for $i=0,1,\ldots, b-1$ in the regular local ring 
$\widehat{{\mathcal O}_{X,x}}$.  We would like to express these conditions in terms of elements of ${\mathcal O}_{V,x}$ and  $\mathcal{O}_{X,x}$. 
To this end, we will make use of differential operators. 

Following the arguments in \ref{Local_Gene_Diff_Op}, consider   
$$\begin{array}{r@{\hspace{0pt}}rcl} 
\operatorname{Tay}_{z_1}:& k^{\prime}[[z_1,\ldots,z_d]]  & \longrightarrow & k^{\prime}[[z_1,\ldots,z_d,T]]\\
 & f(z_1,z_2,\ldots,z_d) & \longmapsto & f(z_1+T,z_2,\ldots,z_d)=\sum_j\Delta^j_{z_1}(f)T^j.
\end{array}$$
Then 
$$
\operatorname{Tay}_{z_1}(f)=\sum_{i=0}^\infty (z_1+T)^{i}h_i =  \sum_{i=0}^\infty\left( \sum_{j=0}^{i}{\scriptstyle \left(\hspace{-4pt}\begin{array}{c}i \\[-2pt] j\end{array}\hspace{-4pt}\right)} z_1^{i-j}T^{j}\right)h_i   = \sum_{j=0}^{\infty} \left(\sum_{i=j}^{\infty} {\scriptstyle \left(\hspace{-4pt}\begin{array}{c} i\\[-2pt] j\end{array}\hspace{-4pt}\right)}  z_1^{i-j}h_i\right)T^{j}.$$

And recall that for $j=0,1,\ldots,b-1$, 
\begin{equation}\label{DeltaComp}
\begin{array}{rrcl}
\Delta^j_{z_1}: & \widehat{{\mathcal O}_{V,x}} & \longrightarrow  &  \widehat{{\mathcal O}_{V,x}} \\
 \ & f & \longmapsto &
 h_j+{\scriptstyle \left(\hspace{-4pt}\begin{array}{c} j+1\\[-2pt] j\end{array}\hspace{-4pt}\right)}z_1h_{j+1}+{\scriptstyle \left(\hspace{-4pt}\begin{array}{c} j+2\\[-2pt] j\end{array}\hspace{-4pt}\right)}z_1^{2}h_{j+2}+\ldots 
\end{array}
\end{equation}
is a differential operator of order $j$.

Now, using (\ref{DeltaComp}), 
\begin{equation}
\label{operadorRest}
\begin{array}{l}
\widehat{\mathcal{O}_{V,x}}[\Delta^{b-1}_{z_1}(f)W, \Delta^{b-2}_{z_1}(f)W^{2},\Delta^{b-3}_{z_1}(f)W^{3},\ldots, 
\Delta^1_{z_1}(f)W^{b-1},fW^b]\odot \widehat{{\mathcal O}_{V,x}}[z_1W]= \\[10pt]
=\widehat{\mathcal{O}_{V,x}}  [z_1W, 
h_{b-1}W, 
\Delta^{b-2}_{z_1}(f)W^{2},\Delta^{b-3}_{z_1}(f)W^{3}, \ldots,
\Delta^1_{z_1}(f)W^{b-1},fW^b] =\\[10pt]
=\widehat{\mathcal{O}_{V,x}} [z_1W,   h_{b-1}W,
h_{b-2}W^{2},  
\Delta^{b-3}_{z_1}(f)W^{3}, \ldots, 
\Delta^1_{z_1}(f)W^{b-1},fW^b] =\\[5pt]
=\ldots=\\[5pt]
= \widehat{\mathcal{O}_{V,x}}  [z_1W,   h_{b-1}W, h_{b-2}W^{2},  
\ldots, h_{1}W^{b-1},
h_0W^b]. 
\end{array}
\end{equation}

So,  by  (\ref{operadorRest})  and by the double nature of 
$h_j$ (see (\ref{retractionComple})), one has that  
 $x\in \Sing({\mathcal G}\odot  {\mathcal X})$ if and only if, 
$$\begin{array}{c@{\hspace{3pt}}c@{\hspace{3pt}}l}x&\in &\Sing\left(\widehat{\mathcal{O}_{X,x}}[h_{b-1}W, h_{b-2}W^{2},
 \ldots,  h_1W^{b-1},h_0W^b]\right)\\[7pt]
&=& \Sing\left(\widehat{\mathcal{O}_{X,x}}[(\Delta^{b-1}_{z_1}(f))|_{_X}W,   (\Delta^{b-2}_{z_1}(f))|_{_X}W^{2}, \ldots, 
(\Delta^1_{z_1}(f))|_{_X}W^{b-1}, f|_{_X}W^b]\right).\end{array}$$
Since $V$ is smooth over a perfect field $k$,  each differential operator $\Delta^j_{z_1}$,  
originally defined on $\widehat{{\mathcal O}_{V,x}}$,   
defines, by restriction,  a differential operator of order $j$,   say 
$D^j: {\mathcal O}_{V}(U)\to {\mathcal O}_{V}(U)$,   in a neighborhood $U$ of $x$
(since $V$ is smooth over a perfect field $k$, if $z_1,\ldots,z_d\in {\mathcal O}_{V,x}$ is a regular 
system of parameters, then $dz_1,\ldots,dz_d$ form a basis of $\Omega_{{\mathcal O}_{V}|k}^1$  in a neighborhood $U$ of $x$, now use \cite[Theorem 16.11.2]{EGAIV}). Thus 
$$x\in \Sing({\mathcal G}\odot  {\mathcal X})\Leftrightarrow$$ 
$$\Leftrightarrow x\in \Sing({\mathcal O}_{X,x}[D^{b-1}(f)|_{_X}W, 
D^{b-2}(f)|_{_X}W^2, \ldots
D^1(f)|_{_X}W^{b-1},  f|_{_X}W^b]). $$
From here it follows that 
$$\Sing({\mathcal G}\odot  {\mathcal X})=\Sing
\left((\DDiff(\mathcal{G}))|_{X}\right).$$
 
It is worthwhile mentioning here that, via a local retraction, the differential operators $D^j$ can be interpreted 
as relative differential operators via the  local inclusion ${\mathcal O}_{X,x}\subset {\mathcal O}_{V,x}$ (here we may need to work on an \'etale neighborhood of $x$). Thus, from the 
previous discussion it follows that, in a neighborhood $U\subset V$ of $x$,  
$$f=D^0(f)+z_1D^1(f)+\ldots+z_1^sD^s(f) \mod \langle z_1^{s+1}\rangle. $$

Observe that the equality of singular loci is preserved by any $\mathcal{G}$-$\mathcal{X}$-local sequence.  
On the one hand,  it is clearly preserved by pull-backs 
under smooth morphisms. On the other, notice that the equality  is also preserved after 
permissible monoidal transformations. This follows from   Proposition \ref{localt} and the fact that  retractions can be lifted 
after a permissible monoidal transformation as proved in  \ref{Levantar_local} (so  the expansion in (\ref{desarrollo}) can also be lifted after a monoidal transformation).

\

Now we  show that 
$(\DDiff(\mathcal{G}))|_{X}$ is an ${\mathcal O}_X$-Differential Rees algebra. 
We will argue locally, in a neighborhood of  a point $x\in X\subset V$, and assume, as before, that there is a regular system of parameters $\{z_1,z_2,\ldots,z_d\}\subset {\mathcal O}_{V,x}$ such that $z_1$ is a local equation defining $X$.

\

As in \ref{Local_Gene_Diff_Op}, consider the Taylor expansion on  $\widehat{{\mathcal O}_{V,x}}\simeq k^{\prime}[[z_1,\ldots,z_d]] $,  
$$\begin{array}{r@{\hspace{0pt}}rcl} 
\operatorname{Tay} :& k^{\prime}[[z_1,\ldots,z_d]]  & \longrightarrow & k^{\prime}[[z_1,\ldots,z_d,T_1,\ldots,T_d]]\\
 & f(z_1,z_2,\ldots,z_d) & \longmapsto & f(z_1+T_1,z_2+T_2,\ldots,z_d+T_d)=\sum_{\alpha\in {\mathbb N}^d}\Delta^{\alpha} (f)T^{\alpha},
\end{array}$$
and recall that  for each $\alpha=(\alpha_1,\alpha_2,\ldots,\alpha_d)\in {\mathbb N}^d$, $\Delta^{(\alpha_1,\alpha_2,\ldots,\alpha_d)}$ 
is a differential operator of order $|\alpha|=\alpha_1+\alpha_2+\ldots+\alpha_d$, which defines by restriction 
a differential operator $D^{\alpha}: {\mathcal O}_V(U)\to {\mathcal O}_V(U)$ in some neighborhood $U$ of $x$.  In fact,   $\diff_{V|k}^r$, is locally generated by the $D^{\alpha}$ with 
$|\alpha|\leq r$. 

Using the formal retraction in (\ref{retractionComple}), observe 
that for $\alpha_2+\ldots+\alpha_d\leq r$, 
\begin{equation}\label{DeltaComp2}
\begin{array}{rrcl}
\Delta^{(0,\alpha_2,\ldots,\alpha_d)}: & \widehat{{\mathcal O}_{V,x}} & \longrightarrow  &  \widehat{{\mathcal O}_{V,x}} \\
 \ & f & \longmapsto &
\Delta^{(0,\alpha_2,\ldots,\alpha_d)}(f) 
\end{array}
\end{equation}
generate the $\widehat{{\mathcal O}_{X,x}}$-module of differential operators $\diff^r_{\widehat{{\mathcal O}_{X,x}}|k^{\prime}}$. 
Since $X$ is smooth over a perfect field $k$, these differential operators define, by restriction to ${\mathcal O}_{X,x}$ 
differential operators $D^{\alpha}: {\mathcal O}_X(U)\to {\mathcal O}_X(U)$ in some neighborhood $U$ of $x$. Moreover they  
generate the ${\mathcal O}_{X(U)}$-module $\diff_{X|k}^r$. 

To conclude, it is enough to observe that using the formal retraction in (\ref{retractionComple}), for each 
$\alpha=(0,\alpha_2,\ldots,\alpha_d)\in {\mathbb N}^d$, 
$$(\Delta^{\alpha}f)|_{\{z_1=0\}}=\Delta^{\alpha}(f|_{\{z_1=0\}}),$$
from where it follows that 
$(\DDiff(\mathcal{G}))|_{X}=\DDiff(\mathcal{G}|_{X})$.
 \qed

\section{Proof of Theorem \ref{Canonical}}\label{proof_canonical}

\begin{Paragraph}
\label{ProofTheo} {{\em Proof of Theorem  \ref{Canonical}:}}  Let ${\mathcal G}$ and ${\mathcal K}$ be two  Rees algebras. 
We want to show that  ${\mathscr F}_V({\mathcal K})\subset {\mathscr F}_V({\mathcal G})$ if and only if $\overline{\DDiff(\mathcal{G})}\subset\overline{\DDiff(\mathcal{K})}$.  

The ``if" part  follows easily from 
Example \ref{localintegral}, \ref{localExtDiff}, and Lemma 
\ref{difintegral}, so we only have to prove the ``only if"  part and assume that ${\mathscr F}_V({\mathcal K})\subset {\mathscr F}_V({\mathcal G})$. Let $\DDiff(\mathcal{K})=\bigoplus_nI_nW^n$ and let $\DDiff(\mathcal{G})=\bigoplus_nJ_nW^n$.  

To check whether  a given  Rees algebra is contained in the integral closure of another,  we start by choosing a suitable Veronese action, say $\mathbb{V}_N$, so that both  ${\mathbb V}_N(\DDiff(\mathcal{K}))={\mathcal O}_V[I_NW^N]$ and ${\mathbb V}_N(\DDiff(\mathcal{G}))={\mathcal O}_V[J_NW^N]$ 
are almost-Rees rings (see Lemma \ref{N_adecuado}). By Lemma \ref{almostIntegral}, to show that $\DDiff(\mathcal{G})\subset \overline{\DDiff(\mathcal{K})}$ it is 
 enough to prove that  $J_N\subset \overline{I_N}$. By Remark \ref{IntegCriterBlowup}, 
this can be tested in the normalized blow up of $V$ at $I_N$,  say $V\stackrel{\Theta}{\longleftarrow} \mathscr{B}=\overline{\Bl_{I_N}(V)}$ (observe that 
 $\mathscr{B}$ is the normalized blow up of $V$ with respect to $\DDiff(\mathcal{K})$, see  \ref{normalblowup}). 

Then,   
$\Theta^*({\mathbb V}_N(\DDiff(\mathcal{K})))={\mathcal O}_{\mathscr{B}}[I_N{\mathcal O}_{\mathscr{B}}W^{N}]$ where  $I_N{\mathcal O}_{\mathscr{B}}$ is a  locally invertible sheaf of ideals on $\mathscr{B}$, i.e., 
\begin{equation}
\label{factorizacion}
I_N{\mathcal O}_{\mathscr{B}}W^N=I(\mathbb{H}_1)^{N_1}\cdots I(\mathbb{H}_s)^{N_s}W^N
\end{equation}
for some irreducible and reduced hypersurfaces $\mathbb{H}_1,\ldots,\mathbb{H}_s\subset  \mathscr{B}$.   Thus, again by 
Remark
$\ref{IntegCriterBlowup}$, it is enough to prove that $J_N\mathcal{O}_{\mathscr{B},\mathbb{H}_i}\subset I(\mathbb{H}_i)^{N_i}\mathcal{O}_{\mathscr{B},\mathbb{H}_i}$, for $i=1\ldots s$. 

Next we need to show that there is a {\em sufficiently large} open 
subset $\mathscr{U}$ of $\mathscr{B}$ where the inclusion   ${\mathscr F}_{\mathscr{U}}(\Theta^*(\DDiff(\mathcal{K})))\subset {\mathscr F}_{\mathscr{U}}(\Theta^*(\DDiff(\mathcal{G})))$ holds.
Observe first that since $ \mathscr{B}$ is normal, maybe after removing a closed subscheme $\mathscr{Y}\subset \mathscr{B}$ of codimension at least  two,  the open subset $\mathscr{U}=\mathscr{B}\setminus \mathscr{Y}$ is smooth. 
 This means that 
the restriction of $\Theta$ to  $\mathscr{U}$   is a finite type morphism of smooth schemes.   

Consider an open  cover, $\{\mathscr{U}_{\lambda}\}_{\lambda\in \Lambda}$ of $\mathscr{U}$, so that 
the restriction 
    $V\stackrel{\Theta|_{\mathscr{U}_{\lambda}}}{\longleftarrow} \mathscr{U}_{\lambda}$ factorizes  
as a composition of a smooth morphism,  say $\mathscr{Z}=V\times {\mathbb A}^{t}_k\stackrel{\varphi}{\longrightarrow} V$, followed by a  closed  immersion of  smooth schemes, say $\mathscr{U}_{\lambda}\stackrel{i_\lambda}{\hookrightarrow} \mathscr{Z}$. Then  ${\mathscr F}_\mathscr{Z}(\varphi^*(\DDiff(\mathcal{K})))\subset {\mathscr F}_\mathscr{Z}(\varphi^*(\DDiff(\mathcal{G})))$,  and since 
$\mathscr{U}_{\lambda}\subset \mathscr{Z}$ is a   smooth closed subscheme,   
Proposition \ref{restriction} applies to guarantee  that  we still have the inclusion ${\mathscr F}_{\mathscr{U}_{\lambda}}(\Theta^*(\DDiff(\mathcal{K}))|_{\mathscr{U}_{\lambda}})\subset{\mathscr F}_{\mathscr{U}_{\lambda}}(\Theta^*(\DDiff(\mathcal{G}))|_{\mathscr{U}_{\lambda}})$.

Also by  Proposition \ref{restriction},    $\Theta^*(\DDiff(\mathcal{K}))|_{\mathscr{U}_{\lambda}}$ is a differential Rees algebra, since it 
is the restriction of the differential Rees algebra $\varphi^*(\DDiff({\mathcal K}))$ 
to a smooth  closed subscheme ${\mathscr{U}_{\lambda}}$. Hence, one has that 
$$\Sing\Theta^*(\DDiff({\mathcal K}))|_{\mathscr{U}_{\lambda}}=V(I_N\mathcal{O}_{\mathscr{U}_{\lambda}})$$
(see \ref{difsingular}), and as a consequence, $N_1,\ldots, N_s\geq N$ in 
(\ref{factorizacion}). Under these hypotheses, Proposition \ref{CasoHyper}  ensures that 
   $\mathcal{O}_{{\mathscr{U}_{\lambda}},\mathbb{H}_i}[J_N{\mathcal O}_{{\mathscr{U}_{\lambda}},\mathbb{H}_i}W^{N}]\subset \overline{{\mathcal O}_{{\mathscr{U}_{\lambda},\mathbb{H}_i}}[I(\mathbb{H}_i)^{N_i}{\mathcal O}_{\mathscr{U}_{\lambda},\mathbb{H}_i}W^{N}]}$. By Lemma \ref{almostIntegral} this means that 
 $J_N{\mathcal O}_{\mathscr{U}_{\lambda},\mathbb{H}_i} \subset \overline{I(\mathbb{H}_i)^{N_i}{\mathcal O}_{\mathscr{U}_{\lambda},\mathbb{H}_i} }=\overline{I_N{\mathcal O}_{{\mathscr{U}}_{\lambda},\mathbb{H}_i}}$.   Since $\mathscr{Y}$ is a closed subscheme of codimension greater than or equal to two, $\mathcal{O}_{{\mathscr{U}_{\lambda},\mathbb{H}_i}}=\mathcal{O}_{\mathscr{B},\mathbb{H}_i}$ and so $J_N{\mathcal O}_{\mathscr{B},\mathbb{H}_i} \subset \overline{I_N{\mathcal O}_{\mathscr{B},\mathbb{H}_i}}$. By Remark \ref{IntegCriterBlowup},  it follows that  ${J_N} \subset\overline{I_N}$ in $V$. \qed
\end{Paragraph}

\section{Other equivalence relations}
\label{other_weak_equivalence}

The notion of weak equivalence in Definition \ref{Defweak} relies on that of local sequences introduced in Definition \ref{local_sequenceA}. If we change the notion of local sequence it is plausible that the equivalence relation will change as well. Here 
there 
are two examples of alternative notions of local sequences also useful for the problem of constructive resolution. 

\begin{Paragraph}{\bf Restricting the class of smooth morphisms.}\label{solo_mult_restr} 
Assume that in Definition \ref{local_sequenceA}  we only consider smooth maps of the following type:
\renewcommand{\labelenumi}{\arabic{enumi})}
\begin{enumerate}
\item \label{producto} Projection on the first coordinate, $V_1=V\times \mathbb{A}^n_k \stackrel{\varphi}{\longrightarrow} V$.
\item \label{restriccion}Restriction to a Zariski's open subset  $V_1$ of $V$, $V_1\stackrel{\varphi}{\longrightarrow} V$.
\end{enumerate}

So, we may  consider local sequences over $V$ defined in terms of monoidal transformations and smooth morphism of type \ref{producto}) and \ref{restriccion}). In such case, for a given Rees algebra ${\mathcal G}$,  the definition of $\mathcal{G}$-local sequence has changed. 
As in Definition \ref{Defweak}   one can define an equivalence relation on ${\mathcal O}_V$-Rees algebras now considering this new notion of local sequence. A priori, the new 
equivalence relation is different from the latter. In order to distinguish both equivalence relations,  we denote by $\widetilde{\mathscr{C}}_V(\mathcal{G})$ the equivalence class of ${\mathcal G}$ according to the new notion of local sequence. 

\

We claim that  both equivalence classes coincide, i.e.,  $\widetilde{\mathscr{C}}_V(\mathcal{G})=\mathscr{C}_V(\mathcal{G})$ for any $\mathcal{G}$. 
To prove this claim it is enough to   review the proof of the Duality Theorem (see \ref{ProofTheo}), 
and that of Proposition \ref{CasoHyper}. Note that the smooth morphisms used in the proofs are, in fact, as in \ref{producto}) and \ref{restriccion}). This already ensures that both notions of weak equivalences coincide. Let us underline here that the smooth morphisms of type \ref{producto}) and \ref{restriccion}) are so specific that they do not require the existence of retractions in  Question  \ref{comparar_arriba_abajo}.
\end{Paragraph}

\begin{Paragraph}{\bf Restricting the class of  permissible monoidal transformations.} \label{less_monoidal_transformations} When we transform a Rees algebra by a permissible monoidal transformation, as in \ref{weaktransforms}, an exceptional hypersurface appears. So, when  addressing the    resolution of a Rees algebra, see \ref{weaktransforms}, some exceptional hypersurfaces will be introduced, and in order to get a constructive resolution we will require that they have normal crossings. Now we introduce a new structure adapted to this new information, and a new definition of permissible transformation which will be handy starting in Section \ref{Ejemplos_Invariantes}.

\begin{Definition}\label{def_ob_permi} A \emph{basic object} is a triple $(V,\mathcal{G},E)$ where
\begin{itemize}
\item V is a smooth scheme,
\item $\mathcal{G}$ is a $\mathcal{O}_V$-Rees algebra, and
\item $E=\{H_1,H_2,\ldots,H_r\}$ is a set of smooth hypersurfaces so that their union has normal crossings.
\end{itemize}
\end{Definition}

\begin{Definition}\label{trans_bo}
We say that a smooth closed subscheme $Y\subset \Sing\mathcal{G}$ is \emph{permissible} for $(V,\mathcal{G},E)$ if it satisfies the additional constraint of having normal crossings with the union of the hypersurfaces in $E$. A \emph{permissible transformation} for a basic object is the blow up, $V\leftarrow V_1$, at a permissible center for this basic object. The \emph{transform} of $(V,\mathcal{G},E)$ 
is defined as: $$(V_1,\mathcal{G}_1,E_1)$$ where 
$\mathcal{G}_1$ is the transform of $\mathcal{G}$, as in \ref{weaktransforms}, and $E_1=\{H_1,\ldots, H_r,H_{r+1}\}$. Here $H_i\in E_1$ denotes the strict
transform of $H_i\in E$ 
for $i=1, \ldots, r$, and $H_{r+1}$ is the exceptional hypersurface of 
the blow up.

\end{Definition}

\begin{Paragraph}{\bf Local sequences for basic objects.} \label{def_local_sequence_bo} 
Analogously to Definition \ref{local_sequenceB}, we will define \emph{local sequences for basic objects}. Fix a basic object $(V,\mathcal{G},E)$. We say that a sequence
\begin{equation}\label{local_sequence_bo}
\xymatrix@R=2pt@C=25pt{
(V,\mathcal{G}, E)=(V_0,\mathcal{G}_0, E_0)  & \ar[l]_{\hspace{35pt}\pi_0}
(V_1,\mathcal{G}_1,E_1) & \ar[l]_{\hspace{20pt}\pi_1} \cdots & \ar[l]_{\hspace{-20pt}\pi_{m-1}} 
(V_m,\mathcal{G}_m,E_m)
}
\end{equation}
is a $(V,\mathcal{G},E)$-\emph{local sequence} if, 
for $i=0,1,\ldots,m-1$,  each
$\pi_i $ 
 is either a permissible monoidal transformation for 
$(V_i,{\mathcal G}_i,E_i)$ (and then $(V_{i+1},{\mathcal G}_{i+1}, E_{i+1})$ is the transform of $(V_i,{\mathcal G}_i,E_i)$  
in the sense of  \ref{trans_bo}),
or a smooth morphism (and then ${\mathcal G}_{i+1}$ and $E_{i+1}$ are, respectively, the pull-backs of 
${\mathcal G}_i$ and $E_i$ in $V_{i+1}$).
\end{Paragraph}
Therefore, in this new context, the condition of normal crossings is imposed  in the notion of local sequence.  And as before, we can define a new equivalence relation, now on basic objects in the obvious way. Specifically, we will say that $B=(V,\mathcal{G},E)$ and $B'=(V,\mathcal{K},E)$ are \emph{equivalent basic objects} if $\mathcal{G}$ and $\mathcal{K}$ are equivalent with this 
new (more restrictive) definition of local sequence (since we take into account the hypersurfaces in $E$ in the notion of permissible center). Once again, we claim that two basic objects, $B=(V,\mathcal{G},E)$ and $B'=(V,\mathcal{K},E)$, are equivalent if and only if $\mathcal{G}$ and $\mathcal{K}$ are weakly equivalent as Definition \ref{Defweak}.  

\

To prove the claim, observe that,  
on the one hand, by \ref{solo_mult_restr}, it is enough to consider projections and restrictions to open sets instead of general smooth morphisms, since we already know that this does not change the equivalence relation (see \ref{solo_mult_restr}).  

\

To finish, again we review the proofs of Theorem \ref{Canonical} (see \ref{proof_canonical}), and Proposition \ref{CasoHyper} and adapt the unique point where the permissible transformations are required. Namely, we focus on the local sequence (\ref{sucesion_sin_cruces_normales}) of the proof of Proposition \ref{CasoHyper}. This proof consists in comparing the weighted transforms of two Rees algebras, $\mathcal{G}$ and $\mathcal{H}$, after certain local sequences. So now, we impose that these local sequences be permissible for the basic objects $(V,\mathcal{G},E)$ and $(V,\mathcal{H},E)$ with the common set of hypersurfaces $E=\{H_{1},H_{2},\ldots, H_{r}\}$. Thus, we need the centers to have normal crossings with the transforms of $E$, in the following sequence \begin{equation}\label{secuencia_otra_we}
\begin{xy} 
(0,6)*+{V}; (12,6)*+{U}; (35,6)*+{U\times \mathbb{A}^1_k=V'}; (70,6)*+{V'_1}; (90,6)*+{V'_2};(108,6)*+{\ldots}; (128,6)*+{V'_{m-1}};(148,6)*+{V'_m};
(35,0)*+{Y_0=\mathbb{H}|_{U}\times \{x\}}; 
(70,0)*+{Y_1=\mathbb{H}_1^{\prime}\cap H_{r+1}}; 
(128,0)*+{Y_{m-1}=\mathbb{H}_{m-1}^{\prime}\cap H_{r+m-1}};
{\ar (9,6);(3,6) \POS?*^+!D{\scriptstyle \varphi_1}} ; 
{\ar (21,6);(15,6)\POS?*^+!D{\scriptstyle\varphi_2} };
{\ar (67,6);(47,6)\POS?*^+!D{\scriptstyle\rho_0}}; 
{\ar (86,6);(73,6)\POS?*^+!D{\scriptstyle\rho_1}}; 
{\ar (104,6);(94,6)}; 
{\ar (122,6);(112,6)};
{\ar (144,6);(134,6)\POS?*^+!D{\scriptstyle\rho_{m-1}} };
\end{xy}\end{equation}

Here $\varphi_1$ and $\varphi_2$ are the smooth morphisms of \ref{solo_mult_restr}, and, for $i=0,\ldots m-1$, $\rho_i$ is the blow up of $V'_i$ at $Y_i$, and $H_{r+i+1}$ is its exceptional hypersurface. 

\

We observe that if $\mathbb{H}|_U$ and $E$ have normal crossings then $\mathbb{H}|_U\times \mathbb{A}^1_k$ and $Y_0$ have normal crossings with $E'=\{ H_1|_U\times \mathbb{A}^1_k,H_{2}|_U \times \mathbb{A}^1_k,\ldots, H_{r}|_U\times \mathbb{A}^1_k\}$. So, $\pi_0$ is permissible in this new sense. The following monoidal transformations are also permissible. Note that 
 all the centers are the intersection of the strict transform of $\mathbb{H}$ and the new exceptional divisor.  So the local sequence (\ref{secuencia_otra_we}) is in fact permissible for both, $(V,\mathcal{G},E)$ and $(V,\mathcal{H},E)$, and the arguments in the proof of Proposition \ref{CasoHyper} still hold.

\

Otherwise, if $\mathbb{H}|_U$ and $E$ do not have normal crossings, we can replace $U$ with say $\widetilde{U}=U \setminus F$, where $F$ is closed of codimension at least two, so that $\mathbb{H}|_{\widetilde{U}}$ and $E|_{\widetilde{U}}$ have normal crossings. Therefore, placing a smooth morphism of \ref{solo_mult_restr} among the local sequence, we proceed as in the previous case. 

\end{Paragraph}

\part{Applications}
\label{applications}

\section{The Hilbert-Samuel function and Differential Rees Algebras}
\label{HilbertSamuel}

In the following lines we recall the main ideas behind Hironaka's proof of resolution of singularities, and show some consequences of Theorem \ref{Canonical} when applied in this context. 

\

Let 
$X$ be a variety over a  perfect field. Let $\max \operatorname{HS}_X$ be the maximum value of the Hilbert-Samuel function of $X$, say $\operatorname{HS}_X$, and define the closed set 
$$\vMax \operatorname{HS}_X:=\{x\in X\,: \operatorname{HS}_X(x)=\max \operatorname{HS}_X\}.$$  

Hironaka showed (see \cite{Hironaka77}) that, locally, in an \'etale neighborhood of each  point $x\in \vMax \operatorname{HS}_X$, there is  an immersion of $X$ in a smooth scheme, say $X\subset V$, together with a  pair $(J,b)$ on $V$,    so that $\Sing(J,b)=\vMax \operatorname{HS}_X$. Moreover, he showed that  the pair $(J,b)$ can be defined so that its  resolution  leads to a lowering of $\max \operatorname{HS}_X$.  More precisely, locally,   in an \'etale neighborhood of a point $x\in \vMax\operatorname{HS}_X$,    there are functions $f_1,\ldots, f_s\in {\mathcal O}_{V,x}$ such that 
\begin{equation}
\label{HSi}
\vMax \operatorname{HS}_X=\bigcap_{j=1}^s \vMax\operatorname{HS}_{\{f_j=0\}}.  
\end{equation}
Here $\vMax\operatorname{HS}_{\{f_j=0\}}$ denotes the maximum multiplicity locus of ${\{f_j=0\}}$, for 
 $j=1,\ldots,s$. 
Furthermore,  equality (\ref{HSi}) is preserved under blowing ups at smooth centers contained in $\vMax\operatorname{HS}_X$, and its transforms, until the maximum   value of the Hilbert-Samuel function drops. In other words, let 
\begin{equation}
\label{HilSam}
\xymatrix@R=2pt@C=20pt{
V=V_0 & V_1 \ar[l]_{\hspace{20pt} \rho_{0}} &    \ldots  \ar[l]_{\rho_{1}} & V_{n} \ar[l]_{\rho_{n-1}}\\
X=X_0 & X_1 & \ldots & X_n,
}
\end{equation} 
be a sequence of blow ups at centers $Y_i\subset \vMax\operatorname{HS}_{X_i}$, where $X_i$ denotes the 
strict transform of $X_{i-1}$ in $V_i$, for $i=1,\ldots,n$, and assume that 
$$ { {\max}\,}\operatorname{HS}_{X_0}= {{\max}\,}\operatorname{HS}_{X_1}=\ldots= {{\max}\,}\operatorname{HS}_{X_n}.$$
Then 
\begin{equation}
\label{HSSequence}
\vMax\operatorname{HS}_{X_i}=\bigcap_{j=1}^s \vMax\operatorname{HS}_{\{f_{j,i}=0\}}
\end{equation}
where $f_{j,i}$ denotes a strict transform of $f_{j,{i-1}}$ in $V_i$ for $i=1,\ldots,n$  
(here $f_{j,0}=f_j$). It can   be checked that the equality (\ref{HSi}) is also preserved by pull-backs of smooth morphisms. 

\

Observe that    the description  of the  closed sets  $\vMax\operatorname{HS}_{X_i}$ in (\ref{HilSam}) is similar to that of the tree of closed sets determined by Rees algebras and local sequences on smooth schemes.  More precisely, the equalities  in (\ref{HSi}) 
and (\ref{HSSequence}) can be interpreted as follows. Let $n_j$ be the maximum multiplicity of $\{f_j=0\}$.   Then setting ${\mathcal G}={\mathcal O}_{V}[f_1W^{n_1},\ldots,f_sW^{n_s}]$,  one has that
$$\vMax\operatorname{HS}_X=\Sing{\mathcal G}=\Sing \DDiff({\mathcal G}),$$
and this equality is preserved by local sequences if $\max \operatorname{HS}_X$ does not drop. Therefore, the problem of 
lowering the maximum value of the Hilbert-Samuel function of $X$, in a neighborhood of a point 
$x\in \vMax\operatorname{HS}_X$, is equivalent to that of finding a resolution of the ${\mathcal O}_V$-Rees algebra ${\mathcal G}$. 
 Note, in addition, that as the singular locus of $\mathcal G$ is a prescribed closed set, namely a Hilbert-Samuel  stratum, it is clear that such algebra is well defined up to weak equivalence. 
 
\

Now, assume that $X$ is globally embedded in some smooth scheme $V$. 
 Theorem \ref{Canonical} ensures that there is a canonical way to assign a Rees algebra to 
$\vMax\operatorname{HS}_X$, namely,
 $\overline{\DDiff({\mathcal G})}$, locally, in an \'etale neighborhood of each point 
$x\in \vMax\operatorname{HS}_X$. This implies that  these locally defined Rees algebras, patch so as to define a sheaf of Rees algebras  ${\mathcal D}$ on  $V$ (at least in \'etale topology). By construction, $\Sing {\mathcal D}=\vMax\operatorname{HS}_X$, and this equality is preserved by permissible monoidal transformations at centers 
$Y\subset \vMax\operatorname{HS}_{X_i}\subset V_i$,  and  pull-back by smooth morphisms.     

\

When $V$ is smooth over a field of characteristic zero,  there is a theorem of resolution of Rees algebras (see Section \ref{Ejemplos_Invariantes}).  Therefore, a local resolution of each $\overline{\DDiff({\mathcal G})}$ defines a global resolution of ${\mathcal D}$. This leads to an improvement of the Hilbert-Samuel function  of $X$, and, ultimately,  
to a resolution of singularities of $X$.   This trivializes the local-global problem in resolution 
of singularities regarding to the assignation of basic objects to $\vMax\operatorname{HS}_{X}$.

\ 

On the other hand, it can also be  proved  that 
a constructive resolution of $X$ does not depend on the immersion.

\section{Invariants and Resolution of singularities}
\label{Ejemplos_Invariantes}

As indicated in Section \ref{HilbertSamuel}, Hironaka reduces the problem of resolution of singularities  over a perfect field, and hence that of log-resolution of ideals,  to that of resolution of basic objects (see Section \ref{other_weak_equivalence} for the definitions of basic object and that of permissible transformation of basic objects).  Then he shows  that basic objects can be resolved at least over fields of characteristic zero (cf. \cite{Hironaka64}).

\begin{Definition} A \emph{resolution of a basic object}, $(V,(J,b),E)$, is a finite sequence of permissible transformations, 
\begin{equation}\label{secuencia_de_resolucion}
\xymatrix@R=2pt@C=20pt{
(V,{\mathcal G},E)& \ar[l]_{\rho_0}(V_1,{\mathcal G}_1,E_1) & \ar[l]_{\hspace{25pt}\rho_1} \dots & \ar[l]_{\hspace{-20pt}\rho_{n-1}} (V_n,{\mathcal G}_n,E_n),}
\end{equation} such that  $\Sing{\mathcal G}_n=\emptyset$.
\end{Definition}

Hironaka's proof is existential. In order to achieve constructive resolution we need to go one step further: 
we have to resolve basic objects in a way that two weakly  equivalent basic objects undergo the same resolution. This motivates the study of  \emph{invariants}  (see Definition  \ref{DefinicionInvariante}  below): they  will be used to determine the permissible centers in a resolution sequence as (\ref{secuencia_de_resolucion}). In fact, a constructive (or algorithmic) resolution of singularities is defined 
by fixing a suitable {\em  invariant  function of resolution of basic objects} 
that indicates the centers to blow up in a resolution process (see \ref{Problem_2}).

\begin{Definition}\label{DefinicionInvariante} Let $V$ be a smooth scheme over a perfect field $k$,  and  let $\mathscr{R}$ be the set of all (finitely generated) ${\mathcal O}_V$-Rees algebras. 
Recall that we can assign a closed set to  each Rees algebra 
${\mathcal G}\in \mathscr{R}$, namely $\Sing {\mathcal G}$. Suppose  we assign a value to 
each   ${\mathcal G}\in \mathscr{R}$, at each  
   point $x\in \Sing {\mathcal G}$, and denote it by $\mu_{\mathcal G}(x)$. 
We will say that  $\mu_{\mathcal G}(x)$ is an {\em invariant} if for any 
Rees algebra ${\mathcal K}$ weakly equivalent to ${\mathcal G}$, one has that 
$\mu_{\mathcal G}(x)=\mu_{\mathcal K}(x)$. 
\end{Definition}

\begin{Paragraph}
\label{Problema_Local_Global} {\bf The local-global problem.} The goal of this section   is to present an expository  account  of the use of the Canonicity Principle in questions related to resolution of singularities, particularly when it comes to the problem of the globalization of local invariants. A first step in this direction appears already in the previous section, where the globalization of the Hilbert-Samuel-invariant is discussed for a scheme $X$ included in a smooth scheme $V$ over a perfect field $k$. 
Essentially, the point is: given a basic object $(V, \mathcal G, E)$
we want to construct a resolution. Such a resolution can be achieved, but only over fields of characteristic zero. In this context, the globalization of local invariants also appears in the construction of this  resolution (and hence on that of log-resolution of ideals). 
\end{Paragraph} 

\

The material will be presented in chronological order. We begin   with  part (A) (see \ref{HironakasOrder})    describing  the invariants that appear in \cite{Hironaka74}: 
Hironaka's order function and the $\tau$-invariant. The goal in part (B) (see \ref{FirstSatellite}) is twofold: on the one hand,   
we describe some invariants used in constructive resolution (see B1); on the other, we give some hints about how 
these invariants are used to achieve resolution of basic objects in characteristic zero. Finally,  parts (C) and (D)  
are dedicated to  invariants that grow from a form of  elimination (see \ref{EliminationOrder} and \ref{HOrdFunctions}).

\

\begin{Paragraph}\label{HironakasOrder} {\bf  (A) Hironaka's main invariants.} \cite{Hironaka74}, \cite{Giraud1975}, \cite{Hironaka70}, \cite{Hironaka70Certain},  \cite{Hironaka77}, \cite{Oda1973} 
\end{Paragraph}
\noindent {\bf Hironaka's order function.} Let $\mathcal{G}=\bigoplus_n I_nW^n$  be a Rees algebra on   a smooth 
scheme $V$ over a perfect field $k$.   Recall 
that 
Hironaka's order function is  defined as: 
\begin{equation}\label{defeord}
\begin{array}{rrcl}
\ord^{}_{{\mathcal G}}: & \Sing {\mathcal G}& \longrightarrow &  \mathbb{Q}\geq 1\\
 & x & \longmapsto &  \ord^{}_{{\mathcal G}}(x)=\inf_{n\geq 1} \left\{\frac{\nu_x(I_n)}{n}\right\},
\end{array}
\end{equation}
(see \ref{orderRees}). The statements in  
\ref{SingOrdIntegral}, \ref{difsingular} and Theorem \ref{CanonicalChoice} guarantee that 
Hironaka's order function is an invariant (see also \cite[Theorem 10.10]{Hironaka74}).

\

\noindent {\bf Hironaka's 
 $\tau$-invariant for Rees algebras.}  Let $\mathcal{G}=\bigoplus_n I_nW^n$  be as before, and let
  $x\in \Sing \mathcal{G}$ be a closed point  with residue field $k'$. Fix a regular system of parameters,  
$\{z_1,\dots,z_d\}\subset {\calo}_{V,x}$, and consider the
graded $k'$-algebra associated to its maximal ideal $m_x$,
$\operatorname{Gr}_{m_x}{\mathcal O}_{V,x}$. This graded ring  is isomorphic to a
polynomial ring in $d$-variables with coefficients in $k'$, i.e., $k'[Z_1,\dots , Z_d]$, where $Z_i$ denotes the initial form of $z_i$ in $m_x/m_x^2$.  Note that 
 $\operatorname{Gr}_{m_x}{\mathcal O}_{V,x}$ is the coordinate ring
associated to the tangent space of $V$ at $x$, namely
$\Spec (\operatorname{Gr}_{m_x}{\mathcal O}_{V,x})={\mathbb T}_{V,x}$.
The {\em initial ideal} or {\em tangent ideal} of ${\mathcal G}$ at
$x$, $\operatorname{In}_{x}({\mathcal G})$, is the ideal of
$\operatorname{Gr}_{m_x}{\mathcal O}_{V,x}$ generated by the elements
$\operatorname{In}_x(I_n)$ for all $n\geq 1$. Observe that
$\operatorname{In}_{x}({\mathcal G})$ is zero unless $\ord_{\mathcal G}(x)=1$. The zero set of the tangent ideal in $\Spec\,(\operatorname{Gr}_{m_x}{\mathcal O}_{V,x})$ is the {\em tangent cone} of
${\mathcal G}$ at $x$, ${\mathcal C}_{{\mathcal G},x}$. The {\em $\tau$-invariant of ${\mathcal G}$ at $x$} is the minimum number of variables needed to describe
$\operatorname{In}_{x}({\mathcal G})$. This in turns is the codimension of the largest linear
subspace  ${\mathcal L}_{{\mathcal G},x}\subset {\mathcal
C}_{{\mathcal G},x}$ such that $u+v\in {\mathcal C}_{{\mathcal
G},x}$ for all $u\in {\mathcal C}_{{\mathcal G},x}$ and $v\in
{\mathcal L}_{{\mathcal G},x}$. The $\tau$-invariant of ${\mathcal G}$ at $x$ is denoted by  $\tau_{{\mathcal G},x}$.  
The inclusion ${\mathcal G}\subset \DDiff({\mathcal G})$ defines an inclusion 
${\mathcal
C}_{\DDiff({\mathcal G}),x} \subset {\mathcal
C}_{{\mathcal G},x}$,
and in fact, 
$${\mathcal C}_{\DDiff({\mathcal G}),x}={\mathcal L}_{\DDiff({\mathcal G}),x}={\mathcal L}_{{\mathcal
G},x}.$$ 
Note, in particular, that $\mathcal{G}$, $\overline{\mathcal G}$, and $\DDiff({\mathcal G})$ have the same $\tau$-invariant at all singular point (see for instance \cite[Remark 4.5, Theorem 5.2]{B}), and therefore 
it is an invariant. 

Some of these invariants have also been studied by H. Kawanue and K. Matsuki in the frame of their ``Idealistic filtration program'' (see  \cite{kaw} and \cite{kaw-mat}).

\

\noindent{\bf Key point.} The $\tau$-invariant bounds the local codimension of 
the singular locus of a Rees algebra say ${\mathcal G}$, 
at a given point $x\in \Sing{\mathcal G}$. In particular, 
if $Y\subset \Sing{\mathcal G}$ is a permissible center, then 
$\operatorname{codim}_x\,Y\geq \tau_{{\mathcal G},x}$: note that   $ {\mathbb T}_{Y,x}\subset  {\mathbb T}_{V,x}$, is a linear
subspace, and  that moreover,  ${\mathbb T}_{Y,x}\subset {\mathcal
L}_{{\mathcal G},x}$ for all $x\in Y\subset \Sing{\mathcal
G}$, thus   $\operatorname{codim}_x\,
Y\geq \tau_{{\mathcal G},x}$ (cf. \cite[Theorem 6.5]{BrV}). We say that  {\em  $\mathcal{G}$  is of \em codimensional type} $e$ at $x\in \Sing{\mathcal G}$   if
$\tau_{\mathcal{G},x}\geq e$.  We say that {\em ${\mathcal G}$ is of codimensional type 
$\geq e$} if $\tau_{\mathcal{G},x}\geq e$   for all $x\in \Sing \mathcal{G}$. In this case we will write $\tau_{\mathcal G}\geq e$.

\begin{Paragraph}\label{FirstSatellite} {\bf (B) Invariants that are used in constructive resolution.} \cite{Villa89}, \cite[4.11, 4.15]{EncVil97:Tirol}   
\end{Paragraph}

\noindent{\bf A resolution strategy.} Let $V$ be a smooth $d$-dimensional scheme over a perfect field $k$. Going back to the discussion in the previous paragraph, observe that if  ${\mathcal G}$ is an ${\mathcal O}_{V}$-Rees algebra of 
 codimensional type   $d$, then $\Sing \mathcal G$ consists of finitely many points. In such case, our arguments will show that a resolution can be obtained by simply blowing up these points.  More generally, let $e\geq 1$.  Then, if the codimensional type of ${\mathcal G}$ is $\geq e$ at some $e$-codimensional   component $Y$ of $\Sing {\mathcal G}$, it can 
be shown that $Y$ is smooth and a natural center to blow up (see \cite[Lemma 13.2]{BrV}).  The   point is that if  ${\mathcal G}$ is of codimensional type $\geq e$ at a given point $x\in \Sing {\mathcal G}$, but the 
 codimension of $\Sing {\mathcal G}$ is larger than $e$ at $x$, then it is possible to {\em associate} to ${\mathcal G}$ 
another Rees algebra, say $\widetilde{\mathcal G}$, of codimensional type larger than $e$, and 
     so that 
\begin{equation}
\label{Aumento_tau}
{\mathcal G}\subset \widetilde{\mathcal G}.
\end{equation} 
 This motivates the 
use of induction on the codimensional type of a Rees algebra as a strategy for resolution of basic objects. To this end, we make use of the so called  {\em satellite functions}. 

\

Satellite functions will be presented in (B1) below, where it will also be shown that they  are invariants. 
In this point we find it convenient to use the language of pairs instead of that of Rees algebras, since we think 
that it clarifies the exposition. 

\

Satellite functions play a key role in the inductive strategy for resolution of basic objects in characteristic zero. This is indicated 
in (B2), specially in 
Proposition \ref{prxybo1}. As we will see,  the Canonicity Principle will ensure that the Rees algebra $\widetilde{\mathcal G}$ from (\ref{Aumento_tau})  is unique up to weak equivalence. This settles the local-global problem in constructive resolution (see \ref{Problema_Local_Global}).

\

\noindent{\bf (B1) Satellite functions.}

\

\noindent {\bf The first satellite function \cite{Villa89}, \cite[4.11]{EncVil97:Tirol}.} As indicated before, for the clarity of the exposition,   we will write basic objects in 
terms of pairs instead of Rees algebras (see \ref{Resol_Rees}). 
  Let $(V,(J,b),E)$ be a basic object with $E=\{H_1,\ldots,H_r\}$, and consider \emph{any} local sequence as in Definition \ref{def_local_sequence_bo},
\begin{equation}\label{Atransfuno}
\xymatrix@R=2pt@C=20pt{
(V,(J,b),E)=(V_0,(J_0,b),E_0) &  \ar[l] (V_1,(J_1,b), E_1) & \ar[l] \cdots
& \ar[l] (V_m,(J_m,b),E_m).
}
\end{equation}
Let $\{H_{r+1}, \dots , H_{r+{m'}}\}(\subset E_m)$ with $m'\leq m$ denote the exceptional hypersurfaces introduced by the steps that are permissible monoidal transformations  (i. e., by the steps not given by smooth morphisms). We may assume, for simplicity, that these hypersurfaces are irreducible. Then for $i=1,\ldots,m$ there is a well defined factorization of the sheaves of ideals $J_i \subset \mathcal{O}_{V_i}$, say:
\begin{equation}\label{eqfcword}
J_i=I(H_{r_1})^{b_1}I(H_{r_2})^{b_2}\cdots I(H_{r_{i'}})^{b_{i'}}\cdot
\tilde{J}_i
\end{equation} so that $\tilde{J}_i$ does not vanish along
$H_{r_j}$ for $j=1, \ldots  i'$. Define  $\vword^{}_{(J_i,b)}$ (or simply $\vword^{}_{i}$):
\begin{equation}\label{eqdword}
\begin{array}{rrcl}
\vword^{}_i: & \Sing(J_i,b) & \longrightarrow   & \mathbb{Q}\\ 
 & x & \longmapsto & \vword_i(x)=\frac{\nu_{x}(\tilde{J}_i)}{b}\  (=\ord_{(\tilde{J}_i,b)}\,(x)),
\end{array}
\end{equation}
where $\nu_{x}(\tilde{J}_i)$ denotes the order of
$\tilde{J}_i$ at $\calo_{V_i,x}$.  As we will show 
below,  these functions  derive from Hironaka's order function, 
and hence are invariants.  

\

\noindent {\bf The second satellite function: the inductive function t \cite{Villa89}, \cite[4.15]{EncVil97:Tirol}.}  
 Consider \emph{any} local sequence as in Definition \ref{def_local_sequence_bo}, where now each $V_i\leftarrow V_{i+1}$ is defined with center $Y_i
\subset \vMax \vword_i$, 
\begin{equation}
\label{esol6}
\xymatrix@R=2pt@C=20pt{
(V,(J,b),E)& \ar[l]_{\rho_0} (V_1,(J_1,b),E_1) & \ar[l]_{\hspace{35pt}\rho_1} \cdots
&\ar[l]_{\hspace{-30pt}\rho_{m-1}} (V_m,(J_m,b),E_m),}
\end{equation}
Then,
\begin{equation}\label{eqdgw}
 \max \vword \geq \max
\vword_1\geq \dots \geq \max \vword_m.
\end{equation}

We now define a function $t_m$, only under the assumption that $\max
\vword_m>0$.
Set $l \leq m$   such that
\begin{equation}\label{bowo} \max \vword \geq \ldots \geq \max \vword_{l-1}> \max \vword_{l}= \max \vword_{l+1}\dots= \max \vword_{m},
\end{equation} and write:
\begin{equation}\label{ldde}
E_m=E_m^+ \sqcup E_m^- \text{ (disjoint union)},
\end{equation}
where
$E_m^-$ are the strict transforms of hypersurfaces in $E_{l}$.
Define
\begin{equation}
\label{Def_function_t}
\begin{array}{rrcl}
t_m: & \Sing(J_m,b) & \longrightarrow & {\mathbb Q} \times {\mathbb N} \\
 & x & \longmapsto & t_m(x)=(\vword_m(x),\sharp \{H_i\in E_m^- \,: x\in H_i\})\end{array}
 \end{equation}
where $\mathbb{Q}\times \mathbb{N}$ is a set ordered lexicographically, and $\sharp S$ denotes the total number of elements of a set $S$. We underline that:
\renewcommand{\labelenumi}{\roman{enumi})}
\begin{enumerate}
\item  If each step $(V_i,(J_i,b),E_i) \leftarrow (V_{i+1},(J_{i+1},b),E_{i+1})$ in
(\ref{esol6}) is defined with center $Y_i \subset \vMax t_i $, then
\begin{equation}\label{eqdgt}
 \max t \geq \max t_1\geq \dots \geq \max t_m.
\end{equation}
\item If $\max t_m=( \frac{d}{b}, a)$, then $\max
\vword_m=\frac{d}{b}$. Clearly $\vMax t_m \subset \vMax \vword_m$.
\end{enumerate}
Recall that the functions $t_i$ are defined only if $\max \vword_i>0$. We say that a sequence of transformations is $t$-{\em permissible} when  $Y_i \subset \vMax t_i $ for all $i$. Similarly, a sequence is  $w$-{\em permissible} when $Y_i \subset \vMax \vword_i $ for all $i$.

\

We will show 
next that  these functions  derive from Hironaka's order function, 
and hence are invariants.

\

\noindent {\bf Satellite functions derive from Hironaka's order function.} Let us draw attention here on the fact that the function  $\vword $ from (\ref{eqdword}), and the factorization in
(\ref{eqfcword}), grow from  Hironaka's order function.
Fix $H_{r+i}$ as in (\ref{eqfcword}). Assume, for simplicity, that all steps in sequence    (\ref{Atransfuno}) are permissible monoidal transformations  
with centers $Y_{i-1}\subset \Sing (J_{i-1},b)$, for $i=0,\ldots,m-1$. Then   define the a  function $\operatorname{exp}_i$ along the points in $\Sing(J_i,b)$  by setting 
\begin{equation}
\label{exponentes}
\operatorname{exp}_i(x)=
\left\{\begin{array}{ll}
\frac{b_i}{b}=\frac{\ord_{Y_{i-1}}J_{i-1}}{b}-1\ & \operatorname{ if }  x\in H_{r+i}\cap \Sing(J_i,b); \\
 & \\
0 & \operatorname{ otherwise}. 
\end{array}\right.
\end{equation}
Since $Y_{i-1}\subset \Sing (J_{i-1},b)$, one has that $b_{i}\geq 0$. So, we can express each rational number $\operatorname{exp}_i(x)$ in terms of the functions
$\ord_{(J_{j},b)}$, for $j<i$. More precisely, in terms of the functions $\ord_{(J_{j},b)}$  evaluated at the generic points, say $y_{j}$, of the centers
 $Y_{j} ( \subset V_{j})$ of the monoidal transformation. Finally note that 
by induction on the integer $i$,
$$\vword^d_{(J_i,b)}(x)=\ord_{(J_i,b)}(x)- \operatorname{exp}_1(x)-\operatorname{exp}_2(x)-\dots - \operatorname{exp}_i(x).$$
Thus the  satellite functions  derive from Hironaka's order functions and hence are invariants (they 
take the same value for any basic object $(V,(I,c),E)$ weakly equivalent to $(V,(J,b),E)$).

\

\noindent{\bf (B2) The role of the satellite functions in the inductive strategy for resolution of basic objects.}

\

\noindent We return to the language of Rees algebras in this part. Recall that we are interested in finding a resolution of a given basic object $(V,{\mathcal G}, E)$. 

\

\noindent{\bf The key to  the induction on the codimensional type.} 
A basic object $(V,{\mathcal G}, E)$ is said to be   {\em simple}  if
$\ord_{\mathcal{G}}(x)=1$ for all $x\in \Sing \mathcal{G}$.
Observe that if $(V,{\mathcal G}, E)$ is simple, then $\tau_{{\mathcal G},x}\geq 1$ for all $x\in \Sing {\mathcal G}$, and its transform, by any permissible transformation, is  simple again. 

\

\noindent  {\bf Why simple basic objects?}
Simple basic objects  play a central role    because, in characteristic zero,    their resolution can be addressed in an inductive manner. Traditionally the approach to resolve simple basic objects  by induction makes use of restriction to smooth hypersurfaces (of maximal contact). Here we describe an 
  alternative form of induction for  resolution of simple basic objects making use of the notion of the codimensional 
type of a Rees algebra (thus following the resolution strategy sketched in \ref{FirstSatellite}). 
 At the same time it shows that in this reduction simple basic objects appear only up to equivalence. This already highlights the importance of the notion of invariant from Definition \ref{DefinicionInvariante}  in constructive resolution of singularities.

\begin{Proposition} \label{prxybo1} \cite[Theorems 12.7 and 12.9]{BrV} Let $(V,{\mathcal G}, E)$ be a basic 
object with $V$ smooth over a perfect field $k$. Consider a   $t$-permissible local sequence (see \ref{FirstSatellite}), 
\begin{equation}
\label{sequencet}
\xymatrix@R=2pt@C=25pt{
(V,{\mathcal G},E) & \ar[l] (V_1,{\mathcal G}_1,E_1) & \ar[l] \cdots
& \ar[l] (V_m,{\mathcal G}_m,E_m),
}
\end{equation}
and assume that  $\max  \vword_m>0$. Let  $l$ be the smallest index so that $\max t_{l}=\max t_m$ (see \ref{FirstSatellite}). Then: 

(A) There is a simple ${\mathcal O}_{V_l}$-Rees algebra $\widetilde{\mathcal{G}}$   (or say there is a basic object  $(V_l,\widetilde{ \mathcal{G}},\widetilde{E})$), with the following property: Any {\bf local sequence} starting on $(V_l,\widetilde{\mathcal{G}},\widetilde{E})$, say
\begin{equation}\notag
\xymatrix@R=2pt@C=20pt{
(V_l, \widetilde{\mathcal{G}},\widetilde{E})& \ar[l] (\widetilde{V}_{l+1},  \widetilde{\mathcal{G}}_1,\widetilde{E}_{1}) & \ar[l]\cdots &\ar[l] (\widetilde{V}_{l+S},  \widetilde{\mathcal{G}}_S ,\widetilde{E}_S),
}
\end{equation}
induces a $t$-permissible local sequence starting on  $(V_l, \mathcal{G}_l,E_l)$ (also enlarging the first $l$-steps of sequence (\ref{sequencet})), say:
\begin{equation}
\label{reesollll6bis}
\xymatrix@R=2pt@C=20pt{
(V_l, \mathcal{G}_l,E_l)& \ar[l] (\widetilde{V}_{l+1}, \widetilde{\mathcal{G}}_{l+1},\widetilde{E}_{l+1}) & \ar[l] \cdots
&\ar[l] (\widetilde{V}_{l+S},\widetilde{\mathcal{G}}_{l+S},\widetilde{E}_{l+S}),
}
\end{equation}
with the following condition on the functions $t_j$ defined for this last sequence (\ref{reesollll6bis}):

\ a) $\vMax t_{l+k}=\Sing( \widetilde{\mathcal{G}}_k)$ for $k=0,1, \ldots, S-1$;

\ b) $\max  t_l=\max  t_{l+1}=\cdots = \max  t_{l+S-1} \geq \max  t_{l+S}$;

\ c) $\max  t_{l+S-1} = \max  t_{l+S}$ 
if and only if 
$\Sing (\widetilde{\mathcal{G}}_S) \neq \emptyset$ , in which case $\vMax t_{l+S}=\Sing (\widetilde{\mathcal{G}}_S)$;

(B) (Canonicity) If an ${\mathcal O}_{V_l}$-Rees algebra $\mathcal{K}$    also fulfills (A), then $\mathcal{K}$ and $\widetilde{\mathcal{G}}$ are weakly equivalent.
\end{Proposition}

\noindent Note that     the    claim in (B)   already follows from the conditions imposed in (A).

\

 We emphasize here that   Proposition 
\ref{prxybo1} is valid for basic objects $(V,{\mathcal G},E)$ with 
$V$ smooth over a perfect field $k$. Proposition \ref{prxybo1} is needed for the proof of       Theorem  \ref{Teorema_Monomial} below, which is also valid for smooth schemes over perfect fields. Then: 
\begin{itemize}
\item 
When the characteristic of the base field is zero, using   the output of  Theorem \ref{Teorema_Monomial} 
  it can be shown that resolution of basic objects  follows from a combinatorial argument. 

\item When the characteristic 
  of the base field is positive, Theorem \ref{Teorema_Monomial} says that some form of simplification of the singularities can be achieved. 
  
\end{itemize}
These ideas will be made clearer in the next point dedicated to elimination.

\begin{Paragraph}\label{EliminationOrder} {\bf Elimination.} \cite{hpositive}, \cite{BrV}  
\end{Paragraph}

As was previously indicated,   the resolution of basic objects over fields of characteristic zero can be addressed  by using   induction on the codimensional type (see \ref{HironakasOrder} and \ref{FirstSatellite}), and it is here where the theory 
of {\em elimination} is used. A first step within this approach is given in  
Proposition \ref{prxybo1}, which, in particular, states that there is a canonical way to attach a 
Rees algebra of codimensional type $\geq 1$ to a given a Rees algebra of codimensional type $\geq 0$.  

\

In the following lines, we describe the main ideas of elimination, and state Theorem \ref{Teorema_Monomial} 
as a consequence of Proposition \ref{prxybo1}.  We  will also  indicate why Theorem \ref{Teorema_Monomial}  implies 
resolution of basic objects in characteristic zero, and  the obstruction we find when trying to obtain a similar statement in 
positive characteristic.

\

In what follows, all schemes are suppose to be of finite type over a perfect field $k$. 
 Let $\beta: V^{(d)}\to V^{(d-e)}$ be a smooth morphism of smooth schemes of 
dimensions $d$ and $(d-e)$ respectively, with $0\leq e\leq d$. For each closed point  $z\in V^{(d)}$ denote by 
    $d\beta_z: \mathbb{T}_{V^{(d)},z} \to
\mathbb{T}_{V^{(d-e)},\beta(z)}$   the linear (surjective) map induced on the corresponding tangent spaces.  Let  $\mathcal{G}^{(d)}=\bigoplus_n I_nW^n$ be an ${\mathcal O}_{V^{(d)}}$-Rees algebra, and assume that  
$x\in\Sing \mathcal{G}^{(d)}$ is a  closed point with $\tau_{{\mathcal G}^{(d)},x}\geq e$. We say that $\beta: V^{(d)}\to V^{(d-e)}$ is {\em transversal} to $\mathcal{G}^{(d)}$ at $x$, if
the subspaces ${\mathcal L}_{\mathcal{G}^{(d)},x}$ from \ref{HironakasOrder},  and  $\operatorname{ker}(d\beta_x)$ are in general position in   $\mathbb{T}_{V^{(d)},x}$.

\

For a given ${\mathcal O}_{V^{(d)}}$-Rees algebra ${\mathcal G}^{(d)}$ and a point $x\in \Sing {\mathcal G}^{(d)}$ 
with $\tau_{{\mathcal G}^{(d)},x}\geq e$, it is not hard to construct a smooth morphism $\beta: V^{(d)}\to V^{(d-e)}$ 
transversal to   ${\mathcal G}^{(n)}$ at $x$, at least in an \'etale neighborhood of $x$; moreover, this construction  can be done  in a neighborhood of $x$ (see \cite[\S 8]{BrV}). 

\

If  $\beta: V^{(d)}\to V^{(d-e)}$  is transversal to   ${\mathcal G}^{(d)}$ at $x$, then it can be shown that
$\Sing {\mathcal G}^{(d)}$ and $\beta(\Sing {\mathcal G}^{(d)})$ are homeomorphic, and that a  closed subset $Y\subset \Sing {\mathcal G}^{(d)}$ is smooth if and 
only if $\beta(Y)\subset \beta(\Sing {\mathcal G}^{(d)})$ is smooth (see \cite[8.4]{BrV}).  

\

Transversality is 
preserved by permissible transformations: as indicated above, if $Y\subset \Sing {\mathcal G}^{(d)}$ is a smooth center then 
$\beta(Y)$ is a smooth center, and there is a commutative diagram of permissible transformations and smooth projections:
\begin{equation}
\label{diagrama_transversality}
\begin{xy}
(0,16)*+{V^{(d)}}; (25,16)*+{V_1^{(d)}}; 
(0,0)*+{V^{(d-e)} }; (25,0)*+{ V^{(d-e)}_1,}; (12,8)*+{\circlearrowleft};%
{\ar (20,16);(5,16)_{\rho^{(d)}}};
{\ar (0,13);(0,3)^\beta};
{\ar (18,0);(6,0)_{\rho^{(d-e)}}};
{\ar (23,13);(23,3)^{\beta_1}}
\end{xy}
\end{equation}
 where $\beta_1$ is transversal to the transform ${\mathcal G}_1^{(d)}$ of ${\mathcal G}^{(d)}$ in 
$V_1^{(d)}$ in a neighborhood of any point dominating $x$   (see \cite[\S 9]{BrV}). 

\

In \cite{BrV} it is shown that,  given a differential Rees algebra ${\mathcal G}^{(n)}$,  and a transversal smooth projection  $\beta: V^{(d)}\to V^{(d-e)}$ as before, 
it is possible to define an {\em elimination algebra} ${\mathcal G}^{(d-e)}\subset {\mathcal O}_{{\mathcal O}^{(d-e)}}[W]$. 
This elimination algebra has  
  the following property: given a  sequence of permissible transformations, there are commutative diagrams 
of transversal projections and  transforms of Rees algebras, 
\begin{equation}
\label{diagrama_compatibilidad}
\xymatrix@R=0pt@C=30pt{
\mathcal{G}=\mathcal{G}_0 & \mathcal{G}^{(d)}_1 &  & \mathcal{G}^{(d)}_m \\
V^{(d)}=V^{(d)}_0     \ar[ddd]^\beta     & V^{(d)}_1   \ar[l]_{\rho_0}  \ar[ddd]^{\beta_1}   &  \cdots \ar[l]_{\rho_1} &  V^{(d)}_m \ar[l]_{\rho_{m-1}}  \ar[ddd]^{\beta_m}\\
\\
\\
V^{(d-e)}=V^{(d-e)}_0           & V^{(d-e)}_1   \ar[l]_{\overline{\rho}_0}     & \cdots \ar[l]_{\overline{\rho}_1} &  V^{(d-e)}_m \ar[l]_{\overline{\rho}_{m-1}} \\
\mathcal{G}^{(d-e)}=\mathcal{G}^{(d-e)}_0 & \mathcal{G}^{(d-e)}_1
& & \mathcal{G}^{(d-e)}_m
}
\end{equation}
such that 
\begin{equation}
\label{Inclusion_Estricta}
\beta_i(\Sing {\mathcal G}^{(d)}_i)\subset \Sing {\mathcal G}^{(d-e)}_i
\end{equation}
 for $i=0,\ldots,m$.

\

  Given a differential Rees algebra ${\mathcal G}^{(d)}$  of codimensional type  $\geq e$, and an elimination algebra ${\mathcal G}^{(d-e)}$ 
on some $(d-e)$-dimensional smooth scheme, 
 the function  
{\em $\ord^{(d-e)}_{{\mathcal G}^{(d)}}$} is defined as: 
\begin{equation}\label{orden_e}
\begin{array}{rrcl}
\ord^{(d-e)}_{{\mathcal G}^{(d)}}: &  \Sing{\mathcal G}^{(d)} & \longrightarrow & {\mathbb Q}_{\geq 0}\\ 
 & z & \longmapsto & \ord_{{\mathcal G}^{(d-e)}}(\beta(z)),
\end{array}
\end{equation} 
where  $\ord_{{\mathcal G}^{(d-e)}}$ is the usual Hironaka's order function for a Rees algebra as in \ref{orderRees}.  Note that the elimination algebra provides   information of local nature, and that the choice 
of the local projection for its construction is not unique. 
A fundamental theorem for elimination algebras  is that this information does not depend on the  choice of projection, 
and that it is in fact an invariant (see \cite[Theorem 10.1]{BrV}). 

\

In the same way as in \ref{FirstSatellite}, satellite functions of $\ord^{(d-e)}_{{\mathcal G}^{(d)}}$ can be 
defined. Altogether, these invariants lead to a proof of Proposition \ref{prxybo1}, and hence to the following 
result: 

\begin{Theorem} \label{Teorema_Monomial}  \cite[Part 5]{BrV} \cite[Corollary 6.15]{positive} 
 Let $V^{(d)}$ be a smooth scheme over a perfect field $k$, and  let ${\mathcal G}^{(d)}$ be a differential ${\mathcal O}_{V^{(d)}}$-Rees algebra of codimensional type $\geq e$.   Let $\beta: V^{(d)}\to V^{(d-e)}$ 
 be a transversal projection in a neighborhood of a point $x\in \Sing {\mathcal G}^{(d)}$,  and let ${\mathcal G}^{(d-e)} 
 \subset {\mathcal O}_{V^{(d-e)}}[W]$ 
 be an elimination algebra. Then 
  a sequence of permissible transformations as (\ref{diagrama_compatibilidad}) can be defined, so that  up to integral closure,  
$${\mathcal G}^{(d-e)}_m={\mathcal O}_{V^{(d-e)}_m}[{\mathcal M}W^n],$$
where ${\mathcal M}$ is a locally principal ideal supported on the exceptional divisor of 
$$V^{(d-e)}_0\longleftarrow  V^{(d-e)}_m.$$  
\end{Theorem}

When $\operatorname{char }(k)=0$, the inclusion (\ref{Inclusion_Estricta}) is an equality, and, as a consequence, sequence (\ref{diagrama_compatibilidad}) can be enlarged so as to 
obtain a resolution of ${\mathcal G}^{(d)}$ using arguments of combinatorial nature.   On the other hand, if  $\operatorname{char }(k)=p>0$, the inclusion in (\ref{Inclusion_Estricta}) may be strict. Still, formula 
indicates that the singularities of ${\mathcal G}^{(d)}$ can be, somehow, simplified. In particular, in \cite{BVV} it is shown how sequence (\ref{diagrama_compatibilidad}) can be enlarged 
so as to  obtain  resolution of surfaces in positive characteristic.

\

\noindent {\bf On the representability of the singular locus in fewer variables.} 

\

Let ${\mathcal O}_{V^{(d)}}$ be a smooth scheme over a perfect field $k$, and let ${\mathcal G}^{(d)}$ be an ${\mathcal O}_{V^{(d)}}$-Rees algebra. Fix a transversal smooth projection as before, $\beta: V^{(d)}\to V^{(d-e)}$, in a neighborhood of some point 
$x\in \Sing {\mathcal G}^{(d)}$ where the codimensional type is $\geq e\geq 1$.   If $V^{(d)}\stackrel{\varphi}{\leftarrow} V^{(d)}\times {\mathbb A}_k^r$ 
  is the multiplication by an $r$-dimensional affine space, then there is a commutative diagram 
of smooth morphisms 
\begin{equation}
\label{transversal_smooth}
\xymatrix@R=20pt@C=25pt{
 V^{(d)} \ar[d] ^{\beta}  & \ar[l]_{\varphi^{(d)}} V^{(d)} \times {\mathbb A}_k^r \ar[d]^{\beta\times  id_{{\mathbb A}^r_k}}    \\  
V^{(d-e)} &  V^{(d-e)} \times {\mathbb A}_k^r  \ar[l]_{\varphi^{(d-e)}} 
.}
\end{equation}
One can check  that $\beta\times  id_{{\mathbb A}^r_k}: V^{(d)} \times {\mathbb A}_k^r \longrightarrow  V^{(d-e)} \times {\mathbb A}_k^r$  is transversal to $(\varphi^{{(d)}})^*({\mathcal G}^{(d)})$  at any point of $(\varphi^{(d)})^{-1}(x)$. 

\

Therefore  transversality is preserved by both, 
permissible and smooth morphisms as in (\ref{transversal_smooth}). 
 In this context we introduce the following definition:  we will say that $\Sing {\mathcal G}^{(d)}$ is {\em $\beta$-representable  in 
$(d-e)$-variables}, if there is an ${\mathcal O}_{V^{(d-e)}}$-algebra, say ${\mathcal T}^{(d-e)}$, 
 so that given any local  sequence over $V^{(d)}$ (or over $V^{(d-e)}$) inducing commutative  diagrams as 
 (\ref{diagrama_transversality})   and (\ref{transversal_smooth}), 
\begin{equation}
\label{Comm_Repre}
\xymatrix@R=0pt@C=30pt{
\mathcal{G}=\mathcal{G}_0 & \mathcal{G}^{(d)}_1 &  & \mathcal{G}^{(d)}_m \\
V^{(d)}=V^{(d)}_0     \ar[ddd]^\beta     & V^{(d)}_1   \ar[l]_{\pi_0}  \ar[ddd]^{\beta_1}   &  \cdots \ar[l]_{\pi_1} &  V^{(d)}_m \ar[l]_{\pi_{m-1}}  \ar[ddd]^{\beta_m}\\
\\
\\
V^{(d-e)}=V^{(d-e)}_0           & V^{(d-e)}_1   \ar[l]_{\overline{\pi}_0}     & \cdots \ar[l]_{\overline{\pi}_1} &  V^{(d-e)}_m \ar[l]_{\overline{\pi}_{m-1}} \\
\mathcal{T}^{(d-e)}=\mathcal{T}^{(d-e)}_0 & \mathcal{T}^{(d-e)}_1
& & \mathcal{T}^{(d-e)}_m
}
\end{equation}
one has that $\beta_i(\Sing  {\mathcal G}^{(d)}_i)=\Sing {\mathcal T}^{(d-e)}_i$  for $i=0,1,\ldots,m$ (here 
$\mathcal{G}_i^{(d)}$  and $\mathcal{T}^{(d-e)}_i$  denote  the transforms (or pull-backs) of $\mathcal{G}_{i-1}$  and 
$\mathcal{T}^{(d-e)}_{i-1}$  by $\pi_{i-1}$  and $\overline{\pi}_{i-1}$ respectively).

\

Suppose that ${\mathcal G}^{(d)}$ is a differential Rees algebra  of codimensional type $\geq e$. Let $\beta: V^{(d)}\longrightarrow  V^{(d-e)}$ be a transversal smooth 
projection, and let ${\mathcal G}^{(d-e)}$ be the corresponding elimination algebra. 
Then 
\begin{equation}
\label{Iguales}
\beta(\Sing {\mathcal G}^{(d)})=\Sing {\mathcal G}^{(d-e)}.
\end{equation} 
Moreover, if 
 $Y\subset \Sing {\mathcal G}^{(d)}$ is a permissible center, then $\beta(Y)\subset \Sing {\mathcal G}^{(d-e)}$ is permissible and there is a commutative diagram of smooth projections, transforms of Rees algebras,  and elimination algebras: 
\begin{equation}
\label{diagrama_eliminacion_com}
\xymatrix@R=0pt@C=30pt{
\mathcal{G} & \mathcal{G}^{(d)}_1  \\
V^{(d)}     \ar[ddd]^\beta     & V^{(d)}_1   \ar[l]_{\rho}  \ar[ddd]^{\beta_1}   \\
\\
\\
V^{(d-e)}           & V^{(d-e)}_1   \ar[l]_{\overline{\rho}}      \\
\mathcal{G}^{(d-e)}  & \mathcal{G}^{(d-e)}_1 
}
\end{equation}
 
If the characteristic of the base field is zero, then $\beta_1(\Sing {\mathcal G}^{(d)}_1)=\Sing {\mathcal G}^{(d-e)}_1$, and 
  it can be checked that $\Sing {\mathcal G}^{(d)}$ is $\beta$-representable in $(d-e)$-variables. 

\

In positive characteristic,  there is an inclusion, 
$\beta_1(\Sing {\mathcal G}^{(d)}_1)\subset \Sing {\mathcal G}^{(d-e)}_1$, which may be strict. However, this containment is strong enough so as to make the statement of 
Theorem \ref{Teorema_Monomial} valid over perfect fields of arbitrary characteristic.

\begin{Paragraph}\label{HOrdFunctions} {\bf  $\operatorname{H-ord}$-functions.} \cite{BVVV}
\end{Paragraph}

In the latter section we have discussed about the functions 
$\operatorname{ord}^{(d-e)}_{{\mathcal G}^{(d)}}$, defined for algebras 
${\mathcal G}^{(d)}$ of codimensional type $\geq e$. There are other functions, also defined for algebras of codimensional type $\geq e$ which are particularly relevant in positive characteristic. They are the  so called  {\em $\operatorname{H-ord}$-functions}. 

\

More precisely, let ${\mathcal G}^{(d)}$ be a Rees algebra on a $d$-dimensional smooth scheme $V^{(d)}$, of codimensional type $\geq e >1$ at  $x\in \Sing {\mathcal G}^{(d)}$. Then, in a neighborhood of $x$,  it is possible to construct 
a smooth local projection to some $(d-e)$-dimensional smooth scheme  $V^{(d-e)}$,  transversal to ${\mathcal G}^{(d)}$ at $x$, and an elimination algebra ${\mathcal G}^{(d-e)}$, say
$$\begin{array}{rrcl}
\beta: & V^{(d)} &   \longrightarrow  &   V^{(d-e)}\\
 & {\mathcal G}^{(d)} & &  {\mathcal G}^{(d-e)},
\end{array}$$
so that, locally, in an \'etale neighborhood 
of $x$, there are $e$ sections of $\beta$, say, 
$z_1,\ldots, z_e$, and $e$ elements $f_1W^{n_1},\ldots, f_eW^{n_e}
\in {\mathcal G}^{(d)}$,  with 
$$\begin{array}{l}
f_{1}=z_1^{n_1}+a_1^{(1)}z_1^{n_1-1}+\ldots+a_{n_1}^{(1)}\in {\mathcal O}_{V^{(d-e)}}[z_1]\\
\vdots \\
f_{e}=z_e^{n_e}+a_1^{(e)}z_e^{n_e-1}+\ldots+a_{n_e}^{(e)}\in {\mathcal O}_{V^{(d-e)}}[z_e]
\end{array}$$
and so  that, up to weak equivalence, 
$${\mathcal G}^{(d)}={\mathcal O}_V[f_{1}W^{n_1},\ldots, f_{e}W^{n_e}]\odot 
{\mathcal G}^{(d-e)},$$
where, as before, ${\mathcal G}^{(d-e)}$ is an elimination algebra (cf. \cite[Proposition 6.3]{BVVV}). 
With this notation, define:
$$H-\operatorname{ord}^{(d-e)}_{{\mathcal G}^{(d)}}(x):=\operatorname{min}_{\tiny{
\begin{array}{c}
1\leq i\leq e\\
1\leq j\leq n_i
\end{array}}} \left\{\frac{\nu_{\beta(x)}(a_{j}^{(i)})}{j}, 
\operatorname{ord}_{{\mathcal G}^{(d)}}^{(d-e)}(\beta(x))\right\}.$$

We mention here two properties of the function $H-\operatorname{ord}^{(d-e)}_{{\mathcal G}^{(d)}}$:

1) The value $H-\operatorname{ord}^{(d-e)}_{{\mathcal G}^{(d)}}(x)$
 is independent 
of the choice of the sections, and it does not depend 
on the choice of the smooth (transversal) projection either. Moreover, it is 
  an invariant (see \cite[Theorem 6.12]{BVVV}). 

2) $H-\operatorname{ord}^{(d-e)}_{{\mathcal G}^{(d)}}= \operatorname{ord}^{(d-e)}_{{\mathcal G}^{(d)}}$ as functions along 
$\Sing({\mathcal G}^{(d)})$ every time when ${\mathcal G}^{(d)}$ is $\beta$-representable. In particular, the equality always holds in characteristic zero.

\section{On the inductive representation of the multiplicity} 
\label{Non_Maximal}

Let $V$ be a smooth $d$-dimensional scheme over a perfect field $k$, and let $X\subset V$ be a hypersurface. A natural strategy to 
approach a resolution of the singularities of $X$ is to  lower its maximum multiplicity, say  $b$. Let ${\mathcal G}={\mathcal O}_V[I(X)W^b]$. 
Then  a resolution of $(V,{\mathcal G}, E=\{\emptyset\})$, leads, by composition, to a proper and birational morphism $V\leftarrow V^{\prime}$ where the maximum multiplicity of the  strict 
transform of $X$, say $X^{\prime}$, drops below $b$. 

\

According to the discussion in Section \ref{Ejemplos_Invariantes} (see \ref{EliminationOrder}), if the codimensional type of ${\mathcal G}$ is $\geq e$ at a given point $x\in\Sing {\mathcal G}$, then, locally, in a neighborhood of $x$, $\Sing {\mathcal G}$  can be described using  $(d-e)$ variables, and in fact, it can be projected 
(bijectively) to  a smooth $(d-e)$-dimensional scheme.  When the characteristic is zero,  there is an elimination algebra defined on a 
$(d-e)$-dimensional smooth scheme whose resolution leads to a resolution of 
${\mathcal G}$, and hence, to a lowering of the maximum multiplicity of 
$X$ (see (\ref{Iguales}) and (\ref{diagrama_eliminacion_com})).

\

As the following example shows,   there are hypersurfaces in positive characteristic whose maximum multiplicity locus cannot be represented 
by the singular locus of any Rees algebra in lower dimensions 
(see \ref{EliminationOrder} for the notion of representability,  especially diagram (\ref{Comm_Repre})). In fact we find an obstruction already by looking at  plane curves. 

\

Let $k$ be a perfect field,  let $V=\Spec k[z,x]$, and consider the 
curve $C:={\mathbb V}(z^2+x^3)$, whose maximum multiplicity is 
2 at the origin. Let ${\mathcal G}={\mathcal O}_V[(z^2+x^3)W^2]$ and let  
${\mathcal K}={\mathcal O}_V[(z^2+x^5)W^2]$. If the characteristic of $k$ is different from 2, then 
$$\DDiff({\mathcal K})={\mathcal O}_V[zW, x^4W, x^5W^2]\subset \DDiff({\mathcal G})={\mathcal O}_V[zW, x^2W, x^3W^2].$$
So, by Theorem \ref{Canonical}, ${\mathscr F}_V({\mathcal G})\subset {\mathscr F}_V({\mathcal K})$. 

\

Now suppose that  $k$ is a perfect field of characteristic 
2.  Observe that    the codimensional type of ${\mathcal G}$ is greater 
or equal to 1, and that       $\operatorname{H-ord}^{(1)}_{\mathcal G}$   equals $3/2$ at $(0,0)=\Sing {\mathcal G}$ (see 
\ref{HOrdFunctions} (2)). In particular, if $\Sing {\mathcal G}$ were representable by a Rees algebra 
in one variable, say ${\mathcal A}$, then the order of ${\mathcal A}$ 
at the single point of its singular locus would be $3/2$. But then, up to integral closure,  
${\mathcal A}={\mathcal O}_{V}[x^3W^2]\subset {\mathcal G}$, which 
would imply that $zW\in {\mathcal G}$, and thus  ${\mathcal G}={\mathcal O}_V[zW, x^3W^2]$. Hence 
${\mathcal K}\subset {\mathcal G}$ and  
 ${\mathscr F}_V({\mathcal G})\subset {\mathscr F}_V({\mathcal K})$. 
 We will see that this inclusion does not hold.

\

To get a contradiction, suppose that ${\mathscr F}_V({\mathcal G})\subset {\mathscr F}_V({\mathcal K})$. Then, by the  Theorem \ref{Canonical}, $\DDiff({\mathcal K}) \subset 
\overline{\DDiff({\mathcal G})}$. Thus $(z^2+x^5)W^2\in \overline{\DDiff({\mathcal G})}$,  and hence $(z^2+x^5)W^2-(z^2+x^3)W^2=(x^5+x^3)W^2=x^3(x^2+1)W^2=x^3(x+1)^2W^2\in \overline{\DDiff({\mathcal G})}$.  Set ${\mathcal H}={\mathcal O}_V[x^3(x+1)^2W^2, (z^2+x^3)W^2]$. Observe that $${\mathcal G}\subset {\mathcal H} \subset \overline{\DDiff({\mathcal G})}$$ 
and  the three Rees algebras  are weakly equivalent. 

\

\noindent {\bf Step 0.}
Set $U=V\setminus \{x=1\}$. So, up to integral closure,  ${\mathcal H}|_U={\mathcal O}_U[zW, x^3W^2]\subset \overline{\DDiff({\mathcal G})}|_U$. Set ${\mathcal G}={\mathcal G}|_U$ and ${\mathcal H}={\mathcal H}|_U$.

\

\noindent {\bf Step 1.} Consider the product with an affine line, 
$$\xymatrix@R=2pt@C=20pt{
U & \ar[l]_{\hspace{-20pt}\varphi_0} U_1=U\times \mathbb{A}_k^1 
}$$
and let $${\mathcal G}_1=\varphi_0^*({\mathcal G})\ \ \operatorname{ and } 
{\mathcal H}_1=\varphi_0^*({\mathcal H}).$$ 

\

\noindent{\bf Step 2.} Let $t$ be a regular parameter at some closed point of the affine line. Blow up  the point $V(\langle z, x, t\rangle)$, 
$$\xymatrix@R2pt@C=20pt{
U_1 & \ar[l]_{\rho_1} U_2.
}$$

Denote by $E_1$ the exceptional divisor, and by ${\mathcal G}_2$ and 
${\mathcal H}_2$ the  transforms of ${\mathcal G}_1$ and ${\mathcal H}_1$, respectively. 
 Set $z_1=\frac{z}{t}$ and  $x_1=\frac{x}{t}$. Then 
at ${\mathcal O}_{U_2,t}={\mathcal O}_{U_1}[z_1,x_1,t]$, 
$${\mathcal G}_2={\mathcal O}_{U_2,t}[(z_1^2+tx_1^3)W^2]\ \ \operatorname{ and } \ \ {\mathcal H}_2={\mathcal O}_{U_2,t}[z_1W, tx^3_1W^2].$$

\

\noindent{\bf Step 3.} Blow up  the point $V(\langle z_1, x_1, t\rangle)$, 
$$\xymatrix@R=2pt@C=20pt{
U_{2,t} & \ar[l]_{\hspace{3pt}\rho_2} U_3.
}$$
Denote by $E_2$ the exceptional divisor, and by ${\mathcal G}_3$ and 
${\mathcal K}_3$ the  transforms of ${\mathcal G}_2$ and ${\mathcal K}_2$, respectively. 
 Set $z_2=\frac{z_1}{t}$ and  $x_2=\frac{x_1}{t}$. Then 
at ${\mathcal O}_{U_3,t}={\mathcal O}_{U_{2,t}}\left[z_2,x_2,t\right]$, 
$${\mathcal G}_3={\mathcal O}_{U_3,t}[(z_2^2+t^2x_2^3)W^2]\ \ \operatorname{ and }\ \  {\mathcal H}_3={\mathcal O}_{U_3,t}[z_2W, t^2x_2^3W^2].$$

\

\noindent {\bf Step 4.} Now  choose some non-zero element $\alpha\in k$.  Then 
$$z^2_2+t^2x_2^3=z^2_2+t^2(\alpha+x_2+\alpha)^3=z^2_2+t^2(\alpha^{
3/2
})^2+t^2(x_2+\alpha)u=(z_2+t\alpha^{
3/2
})^2+t^2(x_2+\alpha)u,$$
for some non-zero element $u\in {\mathcal O}_{U_3}$.  Set $z_3=z_2+t\alpha^{
3/2
}$ and $x_3=x_2+\alpha$. 
Then in $U_3^{\prime}=U_{3,t}\setminus\left( \{x_2=0\} \cup \{E_1\cap E_2\}\right)$, 
$${\mathcal G}_3={\mathcal O}_{U_3^{\prime}}[(z_3^2+t^2x_3)W^2]$$
and
$${\mathcal H}_3={\mathcal O}_{U_3^{\prime}} [(z_3+t\alpha^{
3/2
})W,t^2W^2]=
{\mathcal O}_{U_3^{\prime}}[(z_3+t\alpha^{
3/2
})W, tW]={\mathcal O}_{U_3^{\prime}}[z_3W, tW].$$

\

\noindent{\bf Step 5.} Multiply by an affine line, 
$$\xymatrix@R=2pt@C=20pt{
U_3^{\prime} & \ar[l]_{\hspace{-20pt}\varphi_3} U_4=U_3^{\prime}\times {\mathbb A}^1_k
}$$
and let 
$${\mathcal G}_4=\varphi_3^*({\mathcal G}_3), \  \ \operatorname{ and }\ \ {\mathcal H}_4=\varphi_3^*({\mathcal K}_3). $$

\

\noindent{\bf Step 6.} Let $s$ be a regular parameter  at some closed point of the new affine line. Blow up the point $V(\langle z_3, t, x_3, s\rangle)$, 
$$\xymatrix@R=2pt@C=20pt{
U_4 &\ar[l]_{\hspace{3pt}\rho_4} U_5.
} $$
Denote by $E_3$ the new exceptional divisor, and let ${\mathcal G}_5$ and ${\mathcal H}_5$ be the  transforms of ${\mathcal G}_4$ and ${\mathcal H}_4$ respectively.  Set $z_4=\frac{z_3}{s}$,  $x_4=\frac{x_3}{s}$ and 
$t_4=\frac{t}{s}$. Then, at ${\mathcal O}_{U_5,s}={\mathcal O}_{U_4^{\prime}}\left[z_4, x_4, t_4,s\right]$, 
$${\mathcal G}_5={\mathcal O}_{U_5,s}[(z_4^2+ st^2_4x_4)W^2] \ \ \operatorname{ and } \ \  {\mathcal H}_5={\mathcal O}_{U_5,s}[z_4W,t_4W].$$

\

\noindent{\bf Step 7.} Blow up the point $V(\langle z_4, t_4, x_4, s\rangle)$, 
$$\xymatrix@R=2pt@C=20pt{U_{5,s} & \ar[l]_{\hspace{3pt}\rho_5} U_6}$$
Denote by $E_4$ the new exceptional divisor, and let ${\mathcal G}_6$ and ${\mathcal H}_6$ be the weighted transforms of  ${\mathcal G}_5$ and ${\mathcal H}_5$ respectively.  Set $z_5=\frac{z_4}{s}$,  $x_5=\frac{x_4}{s}$ and 
$t_5=\frac{t_4}{s}$. Then, at ${\mathcal O}_{U_6,s}={\mathcal O}_{U_5,s}\left[z_5, x_5, t_5,s\right]$, 
$${\mathcal G}_6={\mathcal O}_{{\mathcal O}_{U_6,s}}[(z_5^2+ t^2_5x_5s^2)W^2] \ \ \operatorname{ and } \ \  {\mathcal H}_6={\mathcal O}_{{\mathcal O}_{U_6,s}}[z_5W,t_5W].$$
Thus, 
$$V(\langle z_5,s\rangle )\subset \Sing{\mathcal G}_6$$
but 
$$V(\langle z_5,s\rangle )\nsubseteq\Sing{\mathcal H}_6,$$
contradicting the assumption that ${\mathcal G}$ and ${\mathcal H}$ were weakly equivalent
(i.e., that $\mathscr{F}_V(\mathcal{G})=\mathscr{F}_V(\mathcal{H})$). Thus ${\mathscr F}_V({\mathcal G})\not \subset {\mathscr F}_V({\mathcal K})$.


\begin{bibdiv}

\begin{biblist}



\bib{Aab}{book}{
	author={Abhyankar, S. S.},
	title={Ramification theoretic methods in Algebraic Geometry},
	series={Annals of Mathematics Studies},
	volume={43},
	publisher={Princeton University Press},
	address={Princeton New Jersey},
	date={1959},
}


\bib{AbM}{article}{
	author={Abhyankar, S. S.},
	author={Moh, T. T.},
	title={Newton-Puiseux expansion and generalized Tschirnhausen transformation},
	partial={
		part={I},
		journal={J. Reine Angew. Math.},
		volume={260},
		date={1973},
		pages={47--83},
	},
	partial={
		part={II},
		journal={J. Reine Angew. Math.},
		volume={261},
		date={1973},
		pages={29--54},
	},
}


\bib{AMBL}{article}{
	author={\'Alvarez-Montaner, J.},
	author={Blickle, M.},
	author={Lyubeznik, G.},
	title={Generators of $D$-modules in positive characteristic},
	journal={Math. Res. Lett.},
	volume={12},
	date={2005},
	number={4},
	pages={459--473},
}


\bib{B}{article}{
	author={Benito, A.},
	title={The tau invariant and elimination},
	journal={J. Algebra},
	volume={324},
	date={2010},
	number={8},
	pages={1903--1920},	
} 


\bib{BVV}{article}{
	author={Benito, A.},
	author={Villamayor U., O. E.},
	title={Techniques for the study of singularities with applications to resolution of 2-dimensional schemes},
	journal={Math. Ann.},
	eprint={DOI: 10.1007/s00208-011-0709-5},
	date={2011},
}     


\bib{BVVV}{article}{
	author={Benito, A.},
	author={Villamayor U., O. E.},
	title={On elimination of variables in the study of singularities in positive characteristic},
	status={Preprint},
	eprint={arXiv:1103.3462 [math.AG]},
}


\bib{BM}{article}{
	author={Bierstone, E.},
	author={Milman, P.},
	title={Canonical desingularization in characteristic zero by blowing up the maxima strata of a local invariant},
	journal={Invent. Math.},
	volume={128},
	date={1997},
	number={2},
	pages={207--302},	
}


\bib{BlE2011}{article}{
	author={Blanco, R.},
	author={Encinas, S.},
	title={Coefficient and elimination algebras in resolution of singularities},
	journal={Asian J. Math.},
	volume={15},
	date={2011},
	number={2},
	pages={251--271},
}


\bib{BMS}{article}{
	author={Blickle, M.},
	author={Musta\c{t}\v{a}, M.},
	author={Smith, K. E.},
	title={Discreteness and rationality of $F$-thresholds},
	booktitle={Special volume in honor of Melvin Hochster},
	journal={Michigan Math. J.},
	volume={57},
	date={2008},
	pages={43--61},
}


\bib{BMS2}{article}{
	author={Blickle, M.},
	author={Musta\c{t}\v{a}, M.},
	author={Smith, K. E.},
	title={$F$-thresholds of hypersurfaces},
	journal={Trans. Amer. Math. Soc.},
	volume={361},
	date={2009},
	number={12},
	pages={6549--6565},
} 


\bib{BrV}{article}{
	author={Bravo, A.},
	author={Villamayor U., O. E.},
	title={Singularities in positive characteristic, stratification and simplification of the singular locus},
	journal={Adv. in Math.},
	volume={224},
	date={2010},
	number={4},
	pages={1349-1418},
}








\bib{CJS}{article}{
	author={Cossart, V.},
	author={Jannsen, U.},
	author={Saito, S.},
	title={Canonical embedded and non-embedded resolution of singularities for excellent two-dimensional schemes},
	status={Preprint},
	eprint={arXiv:0905.2191v1 [math.AG]},
}


\bib{CP1}{article}{
	author={Cossart, V.},
	author={Piltant, O.},
	title={Resolution of singularities of threefolds in positive characteristic},
	partial={
		part={I},
		subtitle={Reduction to local uniformization on Artin-Schreier and purely inseparable coverings},
		journal={J. Algebra},
		volume={320},
		date={2008},
		number={3},
		pages={1051--1082},
	},
}

\bib{CP2}{article}{
	author={Cossart, V.},	
	author={Piltant, O.},
	title={Resolution of singularities of threefolds in positive characteristic},
	partial={
		part={II},
		journal={J. Algebra},
		volume={321},
		date={2009},
		number={7},
		pages={1836--1976},
	}
}




\bib{Cut09}{article}{
	author={Cutkosky, S. D.},
	title={Resolution of singularities for 3-folds in positive characteristic},
	journal={Amer. J. Math.},
	volume={131},
	date={2009},
	number={1},
	pages={59--127},
}

\bib{Cut11}{article}{
	author={Cutkosky, S. D.},
	title={A skeleton key to Abhyankar's proof of embedded resolution of characteristic $p$ surfaces.},	
	journal={Asian. J. Math.},
	volume={15},
	date={2011},
	number={3},
	pages={369-416},
}




\bib{EV}{article}{
	author={Encinas, S.},
	author={Villamayor U., O. E.},
	title={Rees algebras and resolution of singularities},
	book={
		title={Proceedings del XVI-Coloquio Latinoamericano de Algebra (Colonia, Uruguay, 2005) },
		series={Biblioteca de la Revista Matematica Iberoamericana}, 
		address={Madrid}
		date={2007},
		publisher={Rev. Mat. Iberoamericana},
	},
	pages={63-85},
}



\bib{EncVil97:Tirol}{article}{
	author={Encinas, S.},
	author={Villamayor U., O. E.},
	title={A Course on Constructive Desingularization and Equivariance},
	pages={147-227},
	book={
		title={Resolution of Singularities},
		subtitle={A research textbook  in tribute to Oscar Zariski},
		editor={Hauser, H.}
		editor={Lipman, J.}
		editor={Oort, F.},
		editor={Quir\'os, A.},
		series={Progr. Math.},
		volume={181},
		address={Birkh\"auser, Basel},
		date={2000},
	},
}



\bib{EncVillaIberoamericana}{article}{
	author={Encinas, S.},
	author={Villamayor U., O. E.},
	title={A new proof of desingularization over fields of characteristic zero},
	journal={Rev. Mat. Iberoamericana},
	volume={19},
	date={2003},
	number={2},
	pages={339--353},
}


\bib{mariluz}{thesis}{
	author={Garcia-Escamilla, M. L.},
	title={Rees algebras, differential operators, and applications to resolution of
singularities},
	status={in preparation},
}




\bib{Giraud1975}{article}{
	author={Giraud, J.},
	title={Contact maximal en caract\'eristique positive},
	journal={Ann. Sci. \'Ecole. Norm. Sup. 4\`eme s\'erie},
	volume={8},
	date={1975},
	number={2},
	pages={201--234},
} 


\bib{EGAIV}{book}{
	author={Grothendieck, A.}, 
	title={\'El\'ements de G\'eom\'etrie Alg\'ebrique},
	subtitle={IV \'Etude locale des sch\'emas et des morphismes de sch\'emas},
	volume={20,24,28,32},
	publisher={Inst. Hautes \'Etudes Sci. Publ. Math.},
	series={Publications Math\'ematiques},
	date={1964--1967},
}


\bib{HH17obstacles}{article}{
	author={Hauser, H.},
	title={Seventeen obstacles for resolution of singularities},
	book={
		title={Singularities},
		subtitle={The Brieskorn Anniversary Volume},
		series={Progr. Math.},
		volume={162},
		publisher={Brikh\"auser Verlag},
		address={Basel},
		date={1998},
	},	
	pages={289--313},
}

\bib{HauserBulletin}{article}{
	author={Hauser, H.},
	title={On the problem of resolution of singularities in positive characteristic 
(or: a proof we are still waiting for)},
	journal={Bull. Amer. Math. Soc. (N.S.)},
	volume={47},
	date={2010},
	number={1},
	pages={1--30},
} 

\bib{Hironaka64}{article}{
	author={Hironaka, H.},
	title={Resolution of singularities of an algebraic variety over a field of characteristic zero},
	partial={
		part={I},
		journal={Ann. of Math.},
		volume={79},
		date={1964},
		number={1},
		pages={109--203},
	},
	partial={
		part={II},
		journal={Ann. of Math.},
		volume={79},
		date={1964},
		number={2},
		pages={205--326},
	},
	series={Second Series},
} 


\bib{Hironaka70}{article}{
	author={Hironaka, H.},
	title={Additive groups associated with points of a projective space},
	journal={Ann. of Math.},
	volume={92},
	date={1970},
	number={2},
	pages={327--334},
} 


\bib{Hironaka70Certain}{article}{
	author={Hironaka, H.},
	title={Certain numerical characters of singularities},
	journal={J. Math. Kyoto Univ.},
	volume={10},
	date={1970},
	number={1},
	pages={151--187},
} 


\bib{Hironaka74}{book}{
	author={Hironaka, H.},
	title={Introduction to the theory of infinitely near singular points},
	series={Memorias de Matem\'atica del Instituto ``Jorge Juan''},
	volume={28},
	publisher={Consejo Superior de Investigaciones Cient\'ificas},
	address={Madrid},
	date={1974},
}


\bib{Hironaka77}{article}{
	author={Hironaka, H.},
        title={Idealistic exponents of a singularity},
        conference={
        		title={J. J. Sylvester Sympos.},
        		address={Baltimore, Md},
        		date={1976},
	},
	book={
		title={Algebraic Geometry},
		publisher={Johns Hopkins University Press},
		address={Baltimore, Md.},
		date={1977},
		series={The Johns Hopkins centennial lectures},
            },
             pages={52--125},
}


\bib{Hironaka03}{article}{
	author={Hironaka, H.},
	title={Theory of infinitely near singular points},
	journal={J. Korean Math. Soc.},
	volume={40},
	date={2003},
	number={5},
	pages={901--920},
} 


\bib{Hironaka05}{article}{
	author={Hironaka, H.},
	title={Three key theorems on infinitely near singularities},
	pages={87--126},
	book={
		title={Singularit\'es Franco-Japanaises},
		series={S\'emin. Congr.},
		volume={10},
		publisher={Soc. Math. France},
		adress={Paris},
		date={2005},
	},
}


\bib{IrenaHuneke}{book}{
	title={Integral Closure of Ideals, Rings, and Modules}, 
	author={Huneke, C.},
	author={Swanson, I.}, 
	series={London Mathematical Society Lecture Note Series},
	volume={336},
	date={2006}, 
	publisher={Cambridge University Press}, 
	address={Cambridge},
}


\bib{kaw}{article}{
	author={Kawanoue, H.},
	title={Toward resolution of singularities over a field of positive characteristic},
	partial={
		part={I},
		subtitle={Foundation; the language of the idealistic filtration},
		journal={Publ. Res. Inst. Math. Sci.},
		volume={43},
		date={2007},
		number={3},
		pages={819--909},
	}	
}

\bib{kaw-mat}{article}{
	author={Kawanoue, H.},
	author={Matsuki, K.},
	title={Toward resolution of singularities over a field of positive characteristic (the idealistic filtration program},
	partial={
		part={II},
		subtitle={Basic invariants associated to the idealistic filtration and their properties},
		journal={Publ. Res. Inst. Math. Sci.},
		volume={46},
		date={2010},
		number={2},
		pages={359--422},
	}	
}


\bib{Lipman}{article}{
	author={Lipman, J.},
	author={Teissier, B.},
	title={Pseudorational local rings and a theorem of Brian\c{c}on-Skoda about integral closures of ideals},
	journal={Michigan Math. J.},
	volume={28},
	date={1981},
	number={1},
	pages={97--116},
} 


\bib{Oda1973}{article}{
	author={Oda, T.},
	title={Hironaka's additive group scheme},
	pages={181--219},
	book={
		title={Number theory, algebraic geometry and commutative algebra, in honor of Yasuo Akizuki},
		date={1973},
		address={Kinokuniya, Tokyo},
	},
} 


\bib{Villa89}{article}{
	author={Villamayor U., O. E.},
	title={Constructiveness of Hironaka's resolution},
	journal={Ann. Sci. \'Ecole. Norm. Sup. 4\`eme s\'erie},
	volume={22},
	date={1989},
	number={1},
	pages={1--32},
}


\bib{Villa92}{article}{
	author={Villamayor U., O. E.},
	title={Patching local uniformizations},
	journal={Ann. Sci. \'Ecole. Norm. Sup. 4\`eme s\'erie},
	volume={25},
	date={1992},
	number={6},
	pages={629--677},	
} 



\bib{hpositive}{article}{
	author={Villamayor U., O. E.},
	title={Hypersurface singularities in positive characteristic},
	journal={Adv. Math.},
	volume={213},
	date={2007},
	number={2},
	pages={687--733},
} 


\bib{positive}{article}{
	author={Villamayor U., O. E.},
	title={Elimination with applications to singularities in positive characteristic},
	journal={Publ. Res. Inst. Math. Sci.},
	volume={44},
	date={2008},
	number={2},
	pages={661--697},
} 


\bib{integraldifferential}{article}{
	author={Villamayor U., O. E.},
	title={Rees algebras on smooth schemes: integral closure and higher differential operators},
	journal={Rev. Mat. Iberoamericana},
	volume={24},
	date={2008},
	number={1},
	pages={213--242},	
} 



\bib{WLL}{article}{
	author={W\l odarczyk, J.},
	title={Simple Hironaka resolution in characteristic
 zero},
	journal={Journal of the A.M.S.},
	volume={18},
	date={2005},
	number={4},
	pages={779--822},	
} 

\end{biblist} 
\end{bibdiv}
\end{document}